%% file: combinatorial.hybrid.tex
\documentclass[11pt,reqno,a4paper]{myart}

\usepackage{amsmath}
\usepackage{amssymb}
\usepackage{amsthm}
\usepackage[pagebackref,linktocpage]{hyperref}
\usepackage{enumitem}
\usepackage{tracefnt}
\usepackage[T1]{fontenc}
\usepackage{pifont}
\usepackage{bbding}
\usepackage{xspace}
\usepackage{ifthen}
\usepackage{multicol}
\usepackage{comment}
\usepackage{calc}

\makeatletter

\DeclareFontFamily{T1}{greek}{}
\DeclareFontShape{T1}{greek}{m}{n}{ <-> psyr }{}
\DeclareFontShape{T1}{greek}{o}{n}{ <-> [0.9] psyr }{}
\DeclareFontShape{T1}{greek}{m}{sl}{ <-> rpsyro }{}

\DeclareFontFamily{T1}{title}{}
\DeclareFontShape{T1}{title}{m}{n}{ <-> ptmr }{}
\DeclareFontShape{T1}{title}{m}{sc}{ <-> ptmrc }{}
\DeclareFontShape{T1}{title}{m}{it}{ <-> ptmri }{}
\DeclareFontShape{T1}{title}{b}{n}{ <-> ptmb }{}

\DeclareFontFamily{T1}{computer}{}
\DeclareFontShape{T1}{computer}{m}{n}{ <-> pcrr }{}

\newcommand{\aupalpha}{\textsf{a}}
\newcommand{\aupsigma}{{\usefont{T1}{greek}{o}{n}s}}
\newcommand{\aalpha}{\textsf{\textsl{a}}}

\renewcommand{\@makefntext}[1]{\noindent\makebox[1.0em][l]{\@makefnmark}#1}

\include{macros}

\newcommand{\hybridsymbol}{\raisebox{-0.5pt}{\text{\tiny\SparkleBold}}}
\newcommand{\hybridmax}{(\hybridsymbol_{\hspace{0.5pt}\mathrm{max}})}
\newcommand{\hybrid}{(\hybridsymbol)}
\newcommand{\hybridboc}{(\hybridsymbol_{\hspace{0.5pt}\mathrm{boc}})}
\newcommand{\hybridbociso}{(\hybridsymbol_{\hspace{0.5pt}\mathrm{boc}}^{\hspace{0.5pt}\cong})}
\newcommand{\pcombsymbol}{\raisebox{-0.5pt}{\text{\tiny\SixStar}}}
\newcommand{\pcomb}{(\pcombsymbol)}
\newcommand{\pcombmax}{(\pcombsymbol_{\mathrm{max}})}
\newcommand{\pcombboc}{(\pcombsymbol_{\mathrm{boc}})}
\newcommand{\pcombbociso}{(\pcombsymbol_{\mathrm{boc}}^{\cong})}

\newcommand{\aprinc}{(\raisebox{-0.5pt}{\text{\tiny\FiveStarLines}})}
\newcommand{\princA}{(\mathrm{A})}
\newcommand{\princAstar}{(\mathrm{A}^*)}
\newcommand{\pstars}{(\star_s)}
\newcommand{\pstarsplus}{(\raisebox{-1pt}{\text{\scriptsize\EightStarTaper}}_{\hspace{-1pt}s})}

\newcommand{\pdstar}{\pstar^\partial}
\newcommand{\alphaproper}{$\aupalpha$-proper\xspace}
\newcommand{\deecmp}{$\mathbb D$-complete\xspace}

\newcommand{\gcomp}{\Game_{\mathrm{cmp}}}
\renewcommand{\ggen}{\Game_{\mathrm{gen}}}

\newcommand{\poset}{\R}

\newcommand{\club}{\mathcal C}
\renewcommand{\gen}{\operatorname{gen}}
\newcommand{\genc}{\operatorname{gen}^+}
\newcommand{\phimin}{\phi_{\min}}
\newcommand{\psimin}{\psi_{\min}}

\newcommand{\psicvx}{\psi_{\mathrm{cvx}}}
\newcommand{\psicls}{\psi_{\mathrm{cls}}}
\newcommand{\powcnt}[1]{[#1]^{\leq\aleph_0}}
\renewcommand{\NS}{\mathrm{NS}}
\renewcommand{\downcl}[1]{{\downarrow\hspace{-3pt}#1}}

\renewcommand{\L}{\mathcal L}
\renewcommand{\spc}{\,\,}

\newtheorem{thm}{Theorem}[section]
\newtheorem*{thmo}{Theorem}
\newtheorem*{coro}{Corollary}
\newtheorem{@lem}{Lemma}[section]
\newtheorem{sublem}{Sublemma}[@lem]
\newtheorem{@cor@lem}{Corollary}[@lem]
\newtheorem{corthm}{Corollary}[thm]
\newtheorem{question}{Question}
\newtheorem*{questiono}{Question}
\newtheorem{prop}[@lem]{Proposition}
\newtheorem{claim}{Claim}

\theoremstyle{definition}
\newtheorem{defn}[@lem]{Definition}
\newtheorem{rem}[@lem]{Remark}
\newtheorem{example}[@lem]{Example}

\theoremstyle{remark}
\newtheorem{notn}[@lem]{Notation}
\newtheorem{term}{Terminology}

\newcommand{\altbullet}{\scriptsize\ding{117}}

\newlength{\label@width}
\newlength{\label@sep}
\newlength{\left@margin}
\newcounter{@temp}

\newenvironment{lem}[1][]{\setcounter{@temp}{0}\begin{@lem}[#1]}{\end{@lem}}

\newenvironment{cor}[2][0]{\addtocounter{@lem}{-#1}
\ifthenelse{#2 = 0}{}{\addtocounter{@cor@lem}{\value{@temp}}}\begin{@cor@lem}\addtocounter{@lem}{#1}}{\setcounter{@temp}{\value{@cor@lem}}\end{@cor@lem}}

\newenvironment{enumeq}{
 \settowidth{\label@width}{(\arabic{equation})}
 \settowidth{\label@sep}{\espc}
 \setlength{\left@margin}{\label@width+\label@sep}
 \begin{enumerate}[label=(\arabic*), ref=\arabic*,
   labelwidth=\label@width, labelsep=\label@sep, leftmargin=\left@margin]
     \setcounter{enumi}{\value{equation}}}{\setcounter{equation}{\value{enumi}}\end{enumerate}}

\newenvironment{principle}[1]{
\setlength{\label@sep}{\labelsep}
\begin{list}{#1}{\setlength{\rightmargin}{0pt}\settowidth{\labelwidth}{#1}\settowidth{\labelsep}{\quad}\setlength{\leftmargin}{\labelwidth+\labelsep}}
\item
\begin{minipage}[t]{\linewidth}
}{
\end{minipage}\\[-4.5pt]
\end{list}
}

\title[Combinatorial and hybrid principles for $\sigma$-directed
families]{\usefont{T1}{title}{m}{n}\LARGE Combinatorial and hybrid principles\\ for
    \aupsigma-directed families\\ of countable sets modulo finite}
\author[James Hirschorn]{\usefont{T1}{title}{m}{sc}
\large James Hirschorn}
\date{\usefont{T1}{title}{m}{it}May 4, 2008}
\address{\usefont{T1}{title}{m}{sc}{No affiliation, Thornhill, ON, CANADA}}
\thanks{This research began while the author was being supported by the Japanese
  Society for the Promotion of Science, Project No.~P04301. However, this research
  was for the most part unsupported.}
\email{\usefont{T1}{computer}{m}{n}{j\_hirschorn@yahoo.com}}
\urladdr{\usefont{T1}{computer}{m}{n}{http://homepage.univie.ac.at/James.Hirschorn/start.html}}

\begin{document}
\maketitle

\thispagestyle{empty}

\begin{abstract}{{\usefont{T1}{title}{b}{n}Abstract}}
%
%
We consider strong combinatorial principles for $\sigma$-directed 
families of countable
sets in the ordering by inclusion modulo finite, e.g.~$P$-ideals of countable
sets. We try for principles as strong as possible while remaining compatible with
$\ch$, and we also consider principles compatible with the existence of nonspecial
Aronszajn trees. 
The main thrust is towards abstract principles with game theoretic
formulations. 
Some of these principles are purely combinatorial, while the ultimate
principles are primarily combinatorial but also have aspects of forcing
axioms.
\newline\newline\footnotesize
\noindent{\usefont{T1}{title}{m}{it}
{2000 Mathematics Subject Classification}\tu.
\normalfont Primary 03E05; Secondary 03E35, 03E50, 06A06, 91A46.}
\newline
\noindent{\usefont{T1}{title}{m}{it}{Keywords.} \normalfont
  Sigma-directed, continuum hypothesis,
  special Aronszajn tree, hybrid principle.}
\vspace{12pt}
\end{abstract}

\newlength{\@contentsheight}
\setlength{\@contentsheight}{189.25pt}
\setlength{\columnseprule}{0.4pt}
\hspace*{176.8pt}\rule{0.4pt}{\@contentsheight}\vspace*{-\@contentsheight}
\addtocontents{toc}{\protect\begin{multicols}{2}}
\vspace*{\@contentsheight}\hrule\vspace*{-\@contentsheight}

{\footnotesize\tableofcontents}

\section{Introduction}
\label{sec:pstar-does-not}

Very often combinatorial phenomena occurring at the first uncountable ordinal can
be expressed by sentences in the second level of the L\'evy hierarchy over the
structure $(H_{\aleph_2},\in)$ or more generally the structure
$(H_{\aleph_2},\in,\NS)$, where $\NS$ denotes the ideal of nonstationary subsets of
$\oone$. For example, the continuum hypothesis ($\ch$) is $\Sigma_2$ (i.e.~of the form
$\ulc\exists x\,\forall y\spc\varphi(x,y)\urc$ where the formula $\varphi$ has no
unbounded quantifiers)\footnote{More precisely, it can be expressed as a $\Sigma_2$
  statement, e.g.~$\ulc$there exists a function on~$\oone$ 
  whose range includes all of the real numbers $\reals\urc$.} 
over $(H_{\aleph_2},\in)$ and Souslin's hypothesis ($\sh$) is $\Pi_2$ (i.e.~of the form
$\ulc\forall x\,\exists y\spc\varphi(x,y)\urc$) over $(H_{\aleph_2},\in)$. 
And there is a spectrum of axioms with axioms such as $V=L$ at one end, 
which tend to minimize the collection of $\Pi_2$ sentences that hold
over $(H_{\aleph_2},\in,\NS)$ (and thus maximize the $\Sigma_2$ sentences), and
axioms such as Martin's Maximum ($\mm$) at the other end, 
maximizing the $\Pi_2$ sentences
in the theory of $(H_{\aleph_2},\in,\NS)$. Similarly, there are $\Pi_2$-poor
(i.e.~$\Sigma_2$-rich) models like $L$ at one end of the spectrum and $\Pi_2$-rich
models such as the optimal $L(\reals)^{\pmax}$ at the other end. 

In the present article we are interested in the latter end of the spectrum where $\Pi_2$
sentences are maximized in the theory of $(H_{\aleph_2},\in,\NS)$. When one wants to
demonstrate $H_{\aleph_2}\models\phi$ for some $\Pi_2$ sentence $\phi$, for example in
the typical case of establishing a consistency result, there are a number of
options. One can try and show that some forcing axiom such as Martin's Axiom
$(\ma)$, or some more powerful forcing axiom such as the proper
forcing axiom $(\pfa)$ or $\mm$,
implies $H_{\aleph_2}\models\phi$ (i.e.~internal forcing). Or one can try to directly 
force the truth of $\phi$ over $H_{\aleph_2}$
with a forcing notion tailored to the sentence $\phi$ (i.e.~external forcing). 
Alternatively, the theory of $\pmax$ can be applied for this purpose. 

There is yet another approach, which is to find combinatorial statements so
strong that they entail numerous $\Pi_2$ statements over $(H_{\aleph_2},\in,\NS)$
including the  one of interest. There are many examples in the literature 
demonstrating the power of this method. What we are getting at is that
if we have some fixed strong combinatorial principle $(\mathrm{P})$, then a proof that
$(\mathrm{P})\to (H_{\aleph_2},\in,\NS)\models\nobreak\phi$
 is generally going to be a purely
combinatorial argument that is relatively short in length, whereas applying either
internal or external forcing requires the construction of a poset along with the
necessary density arguments which is often lengthier.

The advantage is even far more pronounced when we want to produce a model where
$H_{\aleph_2}\models\phi$, but at the same time making sure some other statement
$\psi$ holds. For
example, we shall need to obtain models of various $\Pi_2$ sentences over
$(H_{\aleph_2},\in,\NS)$ subject to the condition that $\ch$ holds. In this
situation, internal forcing may not be an option. In fact, to date there are only
a few known forcing axioms compatible with $\ch$ 
(\cite[\Section3]{MR924672},\cite[VII,~\Section2]{MR1623206},\cite[XVII,~\Section2]{MR1623206},\cite{math.LO/0003115})
one of which we shall use in
section~\ref{sec:main-theorems}. In many cases however, for example in
section~\ref{sec:main-theorems} when establishing
one of the strong combinatorial principles in conjunction with $\psi=\ulc\ch\lands$
there exists a nonspecial Aronszajn tree$\urc$, there is no known suitable
forcing axiom, and external forcing is the only viable option---indeed, 
to the author's knowledge no analogue of $\pmax$ has ever been 
successfully applied to produce models of $\ch$. 
What is more, the method of external forcing to
produce such a model of $\ch$ usually involves iterated forcing constructions that 
tend to be difficult and often long and tedious. 

To make an analogy with computer programming, 
external (typically iterated) forcing and $\pmax$ are like machine languages giving
maximum control over the resulting model, comparable to a low-level
programming language that gives maximum control over the computer. 
On the other
hand, internal forcing, and even more so applying some strong combinatorial
principle, are like high-level programming languages where there may be several 
layers of abstraction (perhaps one or two for internal forcing and two or three for combinatorial
principles) and the model can be constructed using a more human language and with
considerably less effort. In the latter case,
the underlying low-level arguments like iterated forcing or various $\pmax$ constructions can be
safely ignored because they have already been done to obtain the forcing axiom or
the combinatorial principle. This is analogous to the fact that high-level computer
languages are normally automatically compiled into low-level languages that the
computer can understand. 

The first combinatorial principle that we are aware of and that entails a
significant amount of the consistent $\Pi_2$ theory over $H_{\aleph_2}$ is a Ramsey
theoretic statement due to Todor\v cevi\'c written as
$\oone\to(\oone,(\oone\mathop{;}\fin\ \oone))^2$ 
in the 1983 paper~\cite{MR716846}. There are
surely many others with earlier dates, but we must draw the line somewhere as to what
constitutes a `strong' principle. We believe that the following  two much stronger
principles $\princA$ and $\princAstar$ below are due to Todor\v cevi\'c, 
although the first place in the literature where we
were able to find them is in the joint paper~\cite{MR1441232} 
appearing in 1997. 
(It is claimed in a recent article~\cite[Remark~2]{MR2221128} of Todor\v
cevi\'c that the dichotomy $\princA$ is just a restatement of his Ramsey
  theoretic principle appearing in the 1985
  paper~\cite{MR792822}, denoted there as $\pstar$, which is itself a
  slightly strengthened version of
  $\oone\to(\oone,(\oone\mathop{;}\fin\ \oone))^2$.
  However, this is either an error or
  a far-fetched exaggeration.) We need some notation and terminology before
  introducing these principles, both of which are dichotomies.

\begin{notn}
\label{notn:downcl}
We use the standard set theoretic notation $[A]^\lambda$,
$[A]^{\leq\lambda}$ and $[A]^{<\lambda}$ to
denote the collection of subsets $B\subseteq A$ of cardinality $|B|=\lambda$,
the collection of subsets $B\subseteq A$ of cardinality less than or equal to 
$\lambda$ (i.e.~$|B|\leq\lambda$), and the collection of subsets of
cardinality less than $\lambda$, respectively. 
We write $\Fin(A)$ for the set of all finite subsets of $A$. 

We use the standard order
theoretic notation $\downcl\H$ to denote the downwards closure of $\H$ in the
inclusion order, i.e.~$\downcl\H=\bigcup_{x\in\H}\power(x)$. Assume
$\H\subseteq\ideal$ for some set $\ideal$. When we want to take
the downwards closure with respect to some other quasi ordering $\altle$ of $\ideal$, 
we denote it by $\downcl(\H,\altle)$. 
\end{notn}

\begin{term}
\label{term:locally-in}
Let $\H$ be a family of subsets of some fixed set $S$. We say that a subset
$A\subseteq S$ is \emph{locally in $\H$}  if all of its countable subsets are
members of $\H$, symbolically $\powcnt A\subseteq\H$. 
We say that $B\subseteq S$ is \emph{orthogonal} to $\H$, expressed symbolically as
$B\perp\H$, if $B\cap x$ is finite for all $x\in\H$. 
\end{term}

\noindent By noting that for $A$ infinite, $A$ is locally in $\downcl\H$ is equivalent to
\begin{equation}
  \label{eq:44}
  [A]^{\aleph_0}\subseteq\downcl\H,
\end{equation}
while $B$ is orthogonal to $\H$ is equivalent to
\begin{equation}
  \label{eq:45}
  [B]^{\aleph_0}\cap{\downcl\H}=\emptyset,
\end{equation}
we can regard ``locally in'' and ``orthogonal to'' as dual notions. 

There is however
a nonequivalent dualization of orthogonality. 
We will consider the \emph{almost inclusion} quasi order $\subseteqfnt$ 
on any power set, where $x\subseteqfnt y$ if $x\setminus y$ is finite. 
An equivalent formulation of
equation~\eqref{eq:45} is $b\subseteqfnt y^\complement$ for all
$y\in\H$, for all $b\in[B]^{\aleph_0}$. 
This can be dualized as $a\subseteqfnt y$ for some
$y\in\H$, for all $a\in[A]^{\aleph_0}$, 
or equivalently $\powcnt A\subseteq\downcl(\H,\subseteqfnt)$. The latter
is not equivalent to equation~\eqref{eq:44}. We come back to this point just below.

\begin{term}
\label{term:directed}
Recall that a subset of $A$ of a quasi order $(Q,\leq)$ is \emph{directed} if every
pair of elements $a,b\in A$ has an upper bound $c\in A$, i.e.~$a,b\leq c$. More
generally, for a cardinal $\lambda$, $A$ is \emph{$\lambda$-directed} if every
subset of $A$ of cardinality less than $\lambda$ has an upper bound in $A$. 
Thus directed and $\aleph_0$-directed are the same notion. 
\emph{$\sigma$-directed} is the same thing as $\aleph_1$-directed.
\end{term}

\begin{principle}{$\princA$}
Let $S$ be a set. For every directed subfamily $\H$ of $(\powcnt S,\subseteqfnt)$
of cardinality at most~$\aleph_1$, either
\begin{enumerate}[leftmargin=*, label=(\arabic*), ref=\arabic*, widest=2, labelsep=\label@sep]
\item\label{item:47} $S$ has a countable decomposition into singletons and 
pieces locally in~$\downcl\H$, or
\item\label{item:46} there exists an uncountable subset of $S$ orthogonal to $\H$.
\end{enumerate}
\end{principle} 

\noindent The following dichotomy is the dual of $\princA$ in the sense that the
roles of the uncountable subset of $S$ and the countable decomposition have been
swapped. 

\begin{principle}{$\princAstar$}
Let $S$ be a set. For every directed subfamily $\H$ of $(\powcnt S,\subseteqfnt)$
of cardinality at most~$\aleph_1$, either
\begin{enumerate}[leftmargin=*, label=(\arabic*), ref=\arabic*, widest=2,
  labelsep=\label@sep]
\item\label{item:48} there exists an uncountable subset of $S$ locally in
  $\downcl\H$, or
\item\label{item:49} $S$ has a countable decomposition into countably
  many pieces  orthogonal to $\H$.  
\end{enumerate}
\end{principle}

\begin{rem}
\label{r-11}
In the dichotomies appearing in~\cite{MR1441232}, an $\aleph_1$-generated ideal $\ideal$
is specified rather than $\H$. If we were to use the weaker dual notion to
orthogonality mentioned above, where $\downcl\H$ is replaced with $\downcl(\H,\subseteqfnt)$,
then $\downcl(\H,\subseteqfnt)$ is the ideal generated by $\H\cup\Fin(S)$, and thus with
this modification our versions of the dichotomies become the same as the
originals, and in particular one does not need to
mention the singletons in the first alternative~\eqref{item:47} of the
principle $\princA$. 
In any case, while our statements $\princA$ and $\princAstar$ are formally
stronger than the originals, we will see in remark~\ref{r-12}
that they are in fact equivalent. 

The reason why we do not specify that $\H$ is closed under subsets 
rather than using the notation $\downcl\H$ will become apparent later. 
\end{rem}

The following result is discussed in~\cite{MR1441232}.

\begin{thmo}[Todor\v cevi\'c]
\label{u-11}
$\pfa$ implies $\princA$ and $\princAstar$.
\end{thmo}

\noindent Both of the dichotomies $\princA$ and $\princAstar$ are also known to negate the
continuum hypothesis (cf.~\cite{MR1441232}). 

The next major development in the domain of 
strong combinatorial principles took place in
1995 and was rather exciting. 
A combinatorial principle was formulated by Todor\v cevi\'c, namely
the case $\theta=\oone$ of the principle $\pstar$ below,
strong enough to imply many of the standard
consequences of $\ma_{\aleph_1}$ and $\pfa$ such as Souslin's
Hypothesis ($\sh$), \emph{all
$(\oone,\oone)$ gaps in $\pN\div\Fin$ are indestructible} and more;
and yet, as established by Abraham, consistent with $\ch$. 

\begin{principle}{$\pstar$}
Let $\theta$ be an ordinal.\footnotemark\ 
For every $\sigma$-directed subfamily $\H$ of $(\powcnt\theta,\subseteqfnt)$, either
\begin{enumerate}[leftmargin=*, label=(\arabic*), ref=\arabic*, widest=2, labelsep=\label@sep]
\item\label{item:50} there exists an uncountable $X\subseteq \theta$ locally in
  $\downcl\H$, or
\item\label{item:51} $\theta$ can be decomposed 
into countably many pieces orthogonal to $\H$.
\end{enumerate}
\end{principle}\footnotetext{Henceforth we shall work exclusively with
  families of sets of ordinals. Of course, the ordinal
  structure plays no role in this particular principle $\pstar$, and thus any set can be
  substituted for~$\theta$.} 

\noindent Remark~\ref{r-11} applies here as well. 
We have generalized the principle $\pstar$ so that it applies to 
arbitrary $\sigma$-directed families rather than just
$P$-ideals. In this case the strengthening is purely formal
(cf.~lemma~\ref{l-22}). The principle
$\pstar_\oone$, i.e.~$\pstar$ with $\theta=\oone$, was then generalized to arbitrary
$\theta$ by Todor\v cevi\'c in~\cite{MR1809418}. Let us note here that for any
ordinal of cardinality $\aleph_2$ or greater, this principle has large cardinal
strength (see e.g.~\cite{MR1809418},~\cite{H2}).\label{large-cardinal}
As we are mainly interested in the theory of $(H_{\aleph_2},\in,\NS)$, we
consider $\theta=\oone$ to be the most important case, and $\pstar_\oone$ has no
large cardinal strength as was demonstrated in~\cite{MR1441232}. 

\subsection{Description of results}
\label{sec:description-results}

The independence of Souslin's Hypothesis, 
and Jensen's theorem on its consistency with $\ch$, 
was a milestone in set theory; see for example the book~\cite{MR0384542} devoted to
this. Jensen moreover established that the stronger statement
\emph{all Aronszajn trees are special} is consistent with $\ch$. 
Now the proof 
in~\cite{MR1441232} of the consistency of $\pstar_\oone$ with $\ch$ and of the fact
that  $\pstar_\oone$ implies $\sh$, constituted the shortest proof to date of
Jensen's theorem $\con(\ch+\sh)$---and we should also note that it
relied heavily on the theory of Shelah~\cite{MR675955} 
developed in part for giving a more reasonable proof of
Jensen's theorem.  Hence this begged the following question.

\begin{questiono}[Abraham--Todor\v cevi\'c, 1995]
\label{q-3}
Does $\pstar$ imply that all Aronszajn trees are special\textup?
\end{questiono}

We give the expected negative answer here:

\begin{thm}
\label{u-1}
The conjunction of\/ $\pstar_\oone$ and\/ $\ch$ is consistent, 
relative to\/ $\zfc$, with the existence of a
nonspecial Aronszajn tree. The conjunction of\/ $\pstar$ and\/ $\ch$ is consistent
relative to a supercompact cardinal with the existence of a nonspecial Aronszajn
tree. 
\end{thm}

This can be interpreted as an indication that $\pstar$ is deficient for the study of
trees. For the good behaviour of trees, 
the optimal principles for $\sigma$-directed subfamilies of $\powcnt{\theta}$ 
seem to be the dual $\pdstar$ of
$\pstar$ and the club variation $\pstarc$ of~$\pstar$ described below. The
dualization of $\pstar$ for arbitrary ordinals involves subtleties beyond the scope
of this paper. Thus we just describe the dual of~$\pstar_\oone$.

\begin{principle}{$\pdstar_\oone$}
For every $\sigma$-directed family $\H$ of countable subsets of $\oone$, either
\begin{enumerate}[leftmargin=*, label=(\arabic*), ref=\arabic*, widest=2, labelsep=\label@sep]
\item there is a countable decomposition of $\oone$ into singletons
  and sets locally in
  $\downcl\H$, or
\item there is an uncountable subset of $\oone$ orthogonal to $\H$.
\end{enumerate}
\end{principle}
The proof that $\pstar\to\sh$ is easily adapted to
show $\pdstar_\oone\to\sh^+$, i.e.~$\pdstar_\oone$ implies all Aronszajn trees are special (this
follows immediately from lemmas~\ref{l-20} and~\ref{l-19} and 
corollary~\ref{l-21}). The
consistency of $\pdstar$ with $\ch$ was an open problem of Abraham--Todor\v
c\-evi\'c, and the theory developed in the present paper is applied in~\cite{Hir} to obtain a
positive solution. 

The principle $\pstarc$ is a variation of $\pstar$ due to Eisworth from 1997 
(see~\cite[Question~2.17]{MR1804704}), at least in
the most important case $\theta=\oone$, 
and its consistency with $\ch$ remains an open problem.
\begin{principle}{$\pstarc$}
Let $\theta$ be an ordinal with uncountable cofinality.
For every $\sigma$-directed subfamily $\H$ of $(\powcnt\theta,\subseteqfnt)$, either
\begin{enumerate}[leftmargin=*, label=(\arabic*), ref=\arabic*, widest=2, labelsep=\label@sep]
\item there is closed uncountable subset of $\theta$ locally in $\downcl\H$,
  or
\item there is a stationary subset of $\theta$ orthogonal to $\H$.
\end{enumerate}
\end{principle}
It is proved by the author in~\cite{H2} that $\pstarc\to\sh^+$, and our interest
there was that it also implies that $\sh^+$ holds in any random forcing extension. 
$\pstarc$ is shown to be a consequence of $\pfa$ in~\cite{H2}, and
this was first done by Eisworth for $\theta=\oone$. 

The present paper evolved from the joint work~\cite{NS} of the author with Abraham on the
problem of whether all Aronszajn trees being \emph{nearly special} implies that all
Aronszajn trees are in fact special. Let us say that an $\oone$-tree $T$ is nearly
special if every stationary subset $S\subseteq\oone$ has a stationary subset
$S'\subseteq S$ such that the subtree of $T$ formed by the levels in $S'$ is
special. After obtaining a negative answer to this tree problem in~\cite{NS}, it was
clear that the methods were relevant to the Abraham--Todor\v cevi\'c 
question~(cf.~page~\pageref{q-3}). 

In fact, the strong answer obtained in~\cite{NS}
(cf.~page~\pageref{thmo:AH}) naturally led to a new dichotomy.
\begin{principle}{$\pstarsplus$}
Let $\theta$ be a ordinal of uncountable cofinality. Then there exists a maximal
antichain $\A$ of $(\NS^+_\theta,\subseteq)$ such that for every $\sigma$-directed
subfamily $\H$ of~$(\powcnt\theta,\subseteqfnt)$, either
\begin{enumerate}[leftmargin=*, label=(\arabic*), ref=\arabic*, widest=2, labelsep=\label@sep]
\item\label{item:52} every $S\in A$ has an uncountable relatively closed $C\subseteq S$ locally 
in~$\downcl\H$, or
\item\label{item:53} there exists a stationary subset of $\theta$ orthogonal to $\H$,
\end{enumerate}
\end{principle}
where $\NS^+_\theta$ denotes the coideal of nonstationary subsets of $\theta$.
We will show that~$\pstarsplus_\oone$ 
implies that all Aronszajn trees are nearly special
in section~\ref{sec:nearly-spec-aronsz}. On the other hand, we shall see that
$\pstarsplus_\oone\nrightarrow\sh^+$, 
and thus the main result of~\cite{NS} is also established
here. Furthermore, in theorem~\ref{u-14} we obtain a model of $\ch$, that
satisfies both $\pstar_\oone$ and $\pstarsplus_\oone$ simultaneously, 
but also has a nonspecial Aronszajn tree in it. 

Thus this paper is to a certain extent a generalization of the results in~\cite{NS}
from the realm of $\oone$-trees to the realm of $\sigma$-directed families of
countable sets of ordinals quasi ordered by almost inclusion. This abstraction has
another advantage not discussed in the introduction. A comparison shows that the proofs of the
more general and thus stronger results for abstract families tend to be considerably
shorter and less involved than the corresponding results for the concrete objects,
in this case trees. On the other hand, it should be noted that the results here certainly
do \emph{not} supercede those in~\cite{NS}. In many situtations concrete objects are
preferable to abstract ones and can give a better view of the mathematical
argumentation. Moreover, in the paper~\cite{NS} the methods 
there have an interesting 
application to bases of Aronszajn trees, that does not appear to readily fit into the
present framework. 

We recently noticed (after $\pstarsplus$ was formulated) 
that the principle~$\pstarsplus$ is related to recent work in~\cite{EN}. Consider the
following weakening of $\pstarc_\oone$.
\begin{principle}{$\pstars$}
Every $\sigma$-directed subfamily $\H$ of $(\powcnt\oone,\subseteqfnt)$ has either
    \begin{enumerate}[leftmargin=*, label=(\arabic*), ref=\arabic*, widest=2, labelsep=\label@sep]
    \item  a stationary subset of $\oone$ locally in $\downcl\H$, or
    \item a stationary subset of $\oone$ orthogonal to $\H$.
    \end{enumerate}
\end{principle}
Thus $\pstarsplus_\oone\to\pstars$. 
This dichotomy is named $\mathrm {P}_{22}$ in~\cite{EN} and shown there to be consistent with
$\ch$.

The developments just discussed can be also be observed from the viewpoint described
in~\cite{MR1696852}. Since the actual definition is well beyond the scope of this paper, we roughly
describe a statement $\phi$ as \emph{$\Pi_2$-compact} (cf.~\cite{MR1696852}) 
if there is no pair $\psi_0$
and $\psi_1$ of $\Pi_2$ sentences, possibly with the unary symbol~$\NS$,
such that both $(H_{\aleph_2},\in,\NS)\models\ulc\psi_0\urc\lands\phi$ and
$(H_{\aleph_2},\in,\NS)\models\ulc\psi_1\urc\land\phi$ are consistent, yet
$(H_{\aleph_2},\in,\NS)\models\ulc\psi_0\lands\psi_1\urc\to\lnot\phi$ (thus we are
not relativizing $\phi$ to $H_{\aleph_2}$). Again speaking roughly, any
$\Pi_2$-compact sentence $\phi$, according to the formulation
in~\cite{MR1696852}, can be satisfied by some model $M\models\phi$
that also satisfies every $\Pi_2$
sentence $\psi$ over $(H_{\aleph_2},\in,\NS)$ that can be forced to be true
over $(H_{\aleph_2},\in,\NS)$, in conjunction with $\phi$ (in the
presence of some large cardinals). 

It is our understanding that the unanswered question of whether the continuum
hypothesis is $\Pi_2$-compact (e.g.~\cite[Ch.~10,11]{MR1713438}) 
is a substantial obstacle to continued progress on 
settling the continuum problem. This question is also very closely related to some
of the questions of Shelah in~\cite[\Section2]{MR1804704}. 
Finding stronger and stronger combinatorial
principles compatible with $\ch$ is highly relevant to determining the
$\Pi_2$-compactness of $\ch$. 
For example, if we found two combinatorial principles compatible
with $\ch$ whose conjunction negates $\ch$ then we will have proved it
non-$\Pi_2$-compact.

To our knowledge the $\Pi_2$-compactness of Souslin's hypothesis remains an open
problem, although results in~\cite{MR1696852} and~\cite{MR1683897} suggest a negative
answer. 
We have no idea whether the related question about the $\Pi_2$-compactness of the
existence of a nonspecial Aronszajn tree has ever been considered.
\begin{question}
\label{q-7}
Is $\ulc$\tu{there exists a nonspecial Aronszajn tree}$\urc$ a $\Pi_2$-compact
statement\tu? How about $\lnot\sh^+\lands \ch$\tu?
\end{question}

\noindent In striving for strong combinatorial principles such as $\pstarsplus$ compatible
$\lnot\sh^+\lands\ch$ we are working towards an answer to question~\ref{q-7}. 

We began by using games to facilitate the proof of the consistency of
$\pstarsplus$, and this was highly successful because it 
allowed us to produce e.g.~models of both~$\pstar$ and $\pstarsplus$ while only
arguing once for properties such as $\aupalpha$-properness, the properness isomorphism
and so forth. However, they ended up being absorbed into the theory
itself and led to new and stronger combinatorial principles based on
these games. This resulted in the
addition of a whole new layer of abstraction. 

Until this work,
proving each variation of the original principle $\pstar$ to be consistent
with $\ch$ required repeating long arguments with slight modifications
according to the particular variation. Moreover, within a given
variation, long arguments are again repeated, particularly between
$\aupalpha$-properness and the proper isomorphism condition. 
It is our hope that this additional abstraction 
will result in reduced redundancy in proving other
combinatorial principles compatible with~$\ch$. This has already been
realized in~\cite{Hir}. 

Here is one of the principles without
definitions to give some idea of the concept.\label{principle:hybrid}
\begin{principle}{$\hybrid$}
Let $(\F,\H)$ be a pair of subfamilies of $[\theta]^{\leq\aleph_0}$ for some ordinal
$\theta$, with $\F$ closed under finite reductions. 
If $\H$ is $\F$-extendable, and Complete has a forward 
nonlosing strategy in the
parameterized game $\ggen(\F,\H)$, $\psi$-globally for some $\psi\to\psimin$,
then there exists an uncountable
$X\subseteq\theta$ such that every proper initial segment of~$X$ is in~$\downcl(\H,\sqsubseteq)$. 
\end{principle}
The games all take as parameters 
a pair $(\F,\H)$ of families of countable sets of ordinals. These
principles all assert the existence of an uncountable set all of whose proper
initial segments are in
$\downcl(\H,\sqsubseteq)$. For example, we can obtain the principle $\pstar$ by
inputting $(\downcl\H,\downcl\H)$ into the appropriate game theoretic combinatorial principle,
and then verifying that a certain player in the corresponding class of games has
winning strategy. This drastically reduces the amount of work needed to obtain a
model of $\pstar$ and $\ch$ `from scratch'. 
Of course the consistency of the game theoretic principle must be established 
first, and this is done in section~\ref{sec:main-theorems}.

There are two main parameterized games involved. 
One of the games $\gcomp$ is purely
combinatorial and thus leads to purely combinatorial principles. 
There is one class $\poset$
of forcing notions associated with pairs $(\F,\H)$ of families of countable sets
of ordinals, introduced in section~\ref{sec:forcing-notion}. In the second
parameterized game~$\ggen$, the outcome of the game depends on a genericity condition 
in the poset $\poset(\F,\H)$. Thus the corresponding principles, e.g.~$\hybrid$
above, have aspects of both a combinatorial principle and a forcing axiom, and thus
are hybrid principles. These seem to us to still have much more the essence of a
combinatorial principle as opposed to a forcing axiom. 

A number of questions concerning these principles are posed in
section~\ref{sec:main-theorems}. We also make a comparison with the principle
$\princA$ in section~\ref{sec:p_aleph_1}. Specifically,
we found the ``$\aleph_1$'' in the two dichotomies $\princA$ and
$\princAstar$ for arbitrary $S$ a bit out of place in a combinatorial principle,
and more like something we would expect in a forcing axiom 
(of course in the most important case $S=\oone$ it is quite natural, 
and disappears when we express the
statement over $H_{\aleph_2}$). We suggest a kind of remedy for this in
section~\ref{sec:p_aleph_1}, 
and strengthen Todor\v cevi\'c's theorem on page~\pageref{u-11} in doing so.

\subsection{Credits and acknowledgements}
\label{sec:credits}

This paper evolved from~\cite{NS} which began during the author's
visit to Ben Gurion University in September~2006. 
The idea to use games in the iterated forcing
construction of~\cite{NS} was due to the author, 
but was undoubtedly influenced by the
author's readings of~\cite{math.LO/0003115}. 
The iteration scenario in the proof of
the consistency of $\pstarsplus_\oone$ 
with the existence of a nonspecial Aronszajn tree 
and $\ch$ is based on the scenario from~\cite{NS} which was 
joint work resulting from discussions with Abraham during the visit. 
We thank Uri Abraham for patiently going through the details of
Shelah's theory on preserving
nonspecialness with us, cf.~\cite[Chapter~IX]{MR1623206}, and for
providing an atmosphere conducive to mathematical research during the visit.
\section{Prerequisites}
\label{sec:prerequisites}
\subsection{Terminology}
\label{sec:terminology}

For a model $M$ and a poset $P\in M$, we write $\Gen(M,P)$ for the collection of
filters $G\subseteq P\cap M$ that are generic over $M$, 
i.e.~$D\cap G\cap M\neq\emptyset$ for every dense $D\subseteq P$ in $M$. 
For $p\in P\cap M$,
$\Gen(M,P,p)$ denotes those $G\in\Gen(M,P)$ with $p\in G$. A condition
$q\in P$ is \emph{generic over $M$}, also called
\emph{$(M,P)$-generic} if $q\forces\dot G_P\cap M\in\Gen(M,P)$. We
write $\gen(M,P)$ for the set of $(M,P)$-generic conditions.
As usual, a poset $P$ is proper 
if for every countable elementary $M\prec H_\kappa$,
for $\kappa$ some sufficiently large regular cardinal, with $P\in M$,
every $p\in P\cap M$ has an extension that is generic over $M$,
i.e.~$\gen(M,P,p)\neq\emptyset$. 

We let $\Genc(M,P)$ denote all
$G\in\Gen(M,P)$ with a common extension in~$P$, 
i.e.~some $q\in P$ with $q\geq p$ for all $p\in G$; 
and $\Genc(M,P,p)$ is defined similarly. Following~\cite{hbst}, we say
that a condition is \emph{completely $(M,P)$-generic} 
if $q\forces\dot
G_P\cap M\in\Genc(M,P)$. Note that for a separative poset $P$ this is
equivalent to $\{p\in M:q\geq p\}\in\Genc(M,P)$. We write $\genc(M,P)$
for the set of completely $(M,P)$-generic conditions, and
we say that that a poset~$P$ is \emph{completely proper} 
if for every $M$ as above, $\genc(M,P,p)\neq\emptyset$ for every $p\in
P\cap M$; 
or equivalently, $\Genc(M,P,p)\neq\emptyset$.
It is easily seen that $P$ is completely proper iff it is
proper and does not add new reals. 

A \emph{tower} of elementary submodels refers to a continuous
$\in$-chain $M_0\in\allowbreak M_1\in\cdots$ such that
$\{M_\xi:\xi\le\alpha\}\in M_{\alpha+1}$ for all $\alpha$. 
For any collection $\M$, define $\gen(\M,P,p)=\bigcap_{M\in\M}\gen(M,P,p)$. 
For a tower $\vec M$ of elementary submodels of some
$H_\kappa$, we write $M_0\prec M_1\prec\cdots$ to emphasize the fact that lower
members are elementarily included in higher ones. 
A poset~$P$ is \emph{$\aalpha$-proper} if every tower $\vec M$ of
countable height $\alpha+1$ consisting of
countable elementary submodels of~$H_\kappa$, with $\kappa$ sufficiently
large and regular and $P\in M_0$: every $p\in P\cap M_0$ has an
extension generic over all members of the tower, 
i.e.~$\gen(\{M_\xi:\xi\leq\alpha\},P,p)\neq\emptyset$. For a class
$\E$ of elementary submodels, 
we say that a poset $P$ is \emph{$\E$-$\aalpha$-proper} to indicate that
$\gen(\{M_\xi:\xi\le\alpha\},P,p)\ne\emptyset$ for all $p\in P\cap
M_0$, whenever $M_0\prec\cdots\prec M_\alpha$ is a tower of members of
$\E$ with $P\in M_0$. 

For a countable model $M$, we let
\begin{equation}
  \label{eq:4}
  \delta_M=\sup(\oone\cap M).
\end{equation}
When $M$ is either an elementary submodel of $H_{\aleph_1}$ or a
transitive model then we have $\delta_M=\oone\cap M$. 

A \emph{quasi order} is a pair $(Q,\le)$ where $\le$ is a reflexive
and transitive relation on $Q$. A subset $C\subseteq Q$ is called
\emph{convex} if for all $p\le q$ in $C$, $p\le r\le q$ implies $r\in C$.

For any binary relation $R\subseteq A\times B$ and $x\subseteq A$,
we write $R[x]$ for the image of $x$ under $R$, i.e.~$R[x]=\{y\in
B:R(x,y)\}$.

For any two sets $x$ and $y$ of ordinals, we write $x\sqsubseteq y$ to
indicate that $y$ \emph{end extends} $x$, i.e.~$x$ is an \emph{initial
  segment} of $y$. 

For a tree $T$, we write $T_\alpha$ for the $\alpha\Th$ level. 
An \emph{$\oone$-tree} is a tree of height~$\oone$ with all levels countable. An
\emph{Aronszajn tree} is an $\oone$-tree with no cofinal (i.e.~uncountable)
branches. A tree is called \emph{special} if it can be decomposed into
countably many antichains. Note that a special $\oone$-tree must be Aronszajn. 
For $R\subseteq\oone$ we write $T\restriction R$ for the restriction 
$\bigcup_{\alpha\in R}T_\alpha$ of $T$ to levels in $R$. We write
$\pred_T(t)=\{u\in T:u\lln T t\}$ for the set of predecessor of $t$.

We follow the standard set theoretic convection of writing $V$ for the
class of all sets. 
\subsection{Combinatorics}
\label{sec:combinatorics}

By an \emph{ideal of sets} we of course mean an ideal $\ideal$ 
in some lattice of sets, 
i.e.~$\ideal\subseteq\power(S)$ for some fixed set $S$ and $\ideal$ is closed
under subsets and pairwise unions of its members. We also say that $\ideal$ is an
\emph{ideal on $S$}. For a cardinal $\lambda$, a
\emph{$\lambda$-ideal} is a $\lambda$-complete ideal, i.e.~it is
closed under unions of cardinality less than $\lambda$. A
\emph{$\sigma$-ideal} means an $\aleph_1$-ideal, i.e.~closed under
countable unions.

\begin{defn}
Let $(Q,\leq)$ be a quasi order. A subset $A\subseteq Q$ is \emph{cofinal} in the
quasi ordering if every $q\in Q$ has an $a\geq q$ in $A$. 
To every quasi order, we associate $\J(Q,\leq)\subseteq\power(Q)$
consisting of all noncofinal subsets of $Q$. 
\end{defn}

\begin{lem}
\label{l-1}
Let $(Q,\leq)$ be a quasi order. Then $\J(Q,\leq)$ is a lower set. Moreover\textup:
\begin{enumerate}[leftmargin=*, label=\textup{(\alph*)}, ref=\alph*, widest=b]
\item\label{item:1} If $Q$ has no maximum elements then $Q\subseteq\J(Q,\leq)$, 
i.e.~$\J(Q,\leq)$ contains all of the singletons of $Q$.
\item\label{item:2} If $A\subseteq Q$ then $\J(A,\le)\subseteq\J(Q,\le)$.
\item\label{item:3} If $I\subseteq Q$ is directed, then $\J(I,\le)$ is a
  proper ideal on $Q$. More
  generally, if $I$ is $\lambda$-directed then $\J(I,\le)$ is a $\lambda$-ideal.
\end{enumerate}
\end{lem}
\begin{proof}[Proof of~\tu{\eqref{item:3}}]
Clearly $Q$ is cofinal in itself. Thus $\J(I)$ is a proper ideal, i.e.~$Q\notin\J(I)$.
Suppose $I$ is $\lambda$-directed, and $\A\subseteq\J(I)$ with $|\A|<\lambda$.
For each $A\in\A$, there exists $i_A\in I$ with no element of $A$
above it. By directedness, $\{i_A:A\in\A\}$ has an upper bound $j\in
I$. Then $j$ witnesses that $\bigcup\A\in\J(I)$. 
\end{proof}

We isolate the role of the property that a family of subsets of an ordinal~$\theta$,
has no countable decomposition of $\theta$ into orthogonal pieces. This has already
been done but in less generality in~\cite{MR1809418}.

\begin{prop}
\label{p-3}
Suppose $\F$ is a family of subsets of some ordinal\/ $\theta$.
Let $\alpha\leq\theta$ be the least ordinal with no countable decomposition into sets
orthogonal to $\F$. Then $\alpha$ has uncountable cofinality. 
\end{prop}

\begin{lem}
\label{l-2}
Let $\lambda$ be a cardinal, and
let $\F$ be a directed subset of the quasi order
$(\power(\theta),\subseteqfnt)$ with $\J(\F,\subseteqfnt)$ a $\sigma$-ideal.
Suppose the ordinal $\theta$ satisfies
\begin{enumerate}[leftmargin=*, label=\textup{(\roman*)}, ref=\roman*, widest=ii]
\item\label{item:4} $\theta$ has no decomposition into \tu(strictly\tu) 
less than $\lambda$ many sets orthogonal to~$\F$,
\item\label{item:5} every $\xi<\theta$ has a  decomposition into 
less than $\lambda$ many sets orthogonal to~$\F$.
\end{enumerate}
Then for every family $\X$ of cofinal subsets of $(\F,\subseteqfnt)$ with
$|\X|<\lambda$, and every $\xi<\theta$, there exists $\alpha\geq\xi$ in $\theta$
such that
\begin{equation}
  \label{eq:3}
  \{x\in X:\alpha\in x\}\textup{ is cofinal in $(\F,\subseteqfnt)$}\espc\textup{for all $X\in\X$}.
\end{equation}
\end{lem}
\begin{proof}
Let $\F$ be a directed subset of $(\power(\theta),\subseteqfnt)$ with $\J(\F)$ a
$\sigma$-ideal. 

\begin{sublem}
\label{sublem:cofinal}
For every $\subseteqfnt$-cofinal $X\subseteq\F$, 
if\/ $Y\subseteq\theta$ and $\{x\in X:\alpha\nobreak\in\nobreak x\}\in\J(\F)$ for all\/ $\alpha\in Y$,
then\/ $Y$ is orthogonal to $\F$. 
\end{sublem}
\begin{proof}
Suppose that some $Y\subseteq\theta$ is not
orthogonal to $\F$, say $z\subseteq Y$ is infinite and $z\subseteq y$ for some
$y\in\F$. Assume without loss of generality that $z$ is countable.
Since $X$ is cofinal while $\F$ is $\subseteqfnt$-directed, 
$\{x\in X:z\subseteqfnt x\}$ is cofinal. Therefore, it follows from the fact that $z$ is
countable while $\J(\F)$ is a $\sigma$-ideal that there exists a finite $s\subseteq z$ 
with $\{x\in X:z\setminus s\subseteq x\}$ cofinal, completing the proof.
\end{proof}

Now suppose $\theta$ and $\X$ are as in the hypothesis. Assume towards a
contradiction that for every $\alpha\geq\xi$, 
the set $\{x\in X_\alpha:\alpha\in x\}$ is noncofinal for some $X_\alpha\in\X$. 
For each $X\in\X$, put $Y_X=\{\alpha\geq\xi:X_\alpha=X\}$. 
Then by the sublemma, $\{Y_X:X\in\X\}$ is a decomposition of $\theta\setminus\xi$
into less than $\lambda$ many sets orthogonal to $\F$. 
But by~\eqref{item:5}, we have a decomposition of $\theta$ into less than $\lambda$
many orthogonal sets, contradicting~\eqref{item:4}. 
\end{proof}

Now we isolate the role played by the stronger property that the family has no
stationary orthogonal subset of $\theta$.

\begin{lem}
\label{l-23}
Suppose $S\subseteq\theta$ is stationary.
Let $\F$ be a directed subfamily of $(\power(\theta),\subseteqfnt)$ with
$\J(\F,\subseteqfnt)$ a $\sigma$-ideal, and with no stationary subset
of $S$ orthogonal to $\F$.
Then for every $M\prec H_{\theta^+}$
with $\F,S\in M$ and $|M|<\cof(\theta)$, 
for every cofinal $X\subseteq\F$ in $M$, 
$\{x\in X:\sup(\theta\cap M)\in x\}$ is cofinal. 
\end{lem}
\begin{proof}
Suppose $X\in M$ is a cofinal subset of $(\F,\subseteqfnt)$.
Let $Y\subseteq S$ be the set of all $\alpha<\theta$ 
such that $\{x\in X:\alpha\in x\}\in\J(\F,\subseteqfnt)$. 
Then $Y$ is orthogonal to $\F$, 
because sublemma~\ref{sublem:cofinal} applies here
as well.
Thus by assumption it is not stationary.
Therefore, as $\sup(\theta\cap M)\in\theta$ and $Y\in M$, we cannot have
$\sup(\theta\cap M)\in Y$. 
\end{proof}

We verify here that principles such as $\pstar$, $\pstarsplus$ and $\pstarc$ do not
become weaker when  the additional requirement that the families are $P$-ideals is
imposed. Notice however, that this argument does not apply to the abstract game
theoretic principles. 

\begin{lem}
\label{l-22}
Let $\H$ be a directed subfamily of $(\power(\oone),\subseteqfnt)$. 
Supposing that $X$ is a subset of\/
$\oone$ that is locally in the ideal $\<\H\cup\Fin(\oone)\>$ generated by $\H\cup\Fin(\oone)$, 
there exists a finite $s\subseteq X$ such that $X\setminus s$
is locally in $\downcl\H$. 
\end{lem} 
\begin{proof}
For each $\delta<\oone$, $X\cap\delta\in\<\H\cup\Fin(\oone)\>$. Hence there
exists a finite $\F_\delta\subseteq\H$ and a finite $s_\delta\subseteq\oone$ such that
$X\cap\delta\subseteq \bigcup\F_\delta\cup s_\delta$.
Since $\H$ is directed, $\F_\delta$ has a $\subseteqfnt$-bound $y_\delta\in\H$. 
Hence there is a finite $t_\delta\subseteq\oone$ 
such that $\bigcup\F_\delta\subseteq y_\delta\cup t_\delta$. Thus
\begin{equation}
  \label{eq:47}
  X\cap\delta\setminus(s_\delta\cup t_\delta)\subseteq y_\delta.
\end{equation}
Pressing down, there exists a stationary $S\subseteq\oone$ such that $s_\delta=s$
and $t_\delta=t$ for all $\delta\in S$. It now follows that $X\setminus(s\cup t)$ is
locally in $\downcl\H$. 
\end{proof}

\begin{rem}       
\label{r-12}
We only dealt with the case $\theta=\oone$. 
For general $\theta$ of uncountable cofinality the same argument
applies, so long as $\otp(X)=\oone$. 
This is of no concern for the principle $\pstar$, but for $\pstarsplus$
and $\pstarc$, where we want a closed uncountable set, we need at
least one limit of uncountable cofinality relative to~$X$. 
For example we could have $\otp(X)=\oone+1$, but the top point might
prevent $X$ from being locally in $\downcl\H$. 

Nevertheless, if weaken the principles $\pstarsplus$ or $\pstarc$ so
that $X$ is only required to be closed at limits of countable
cofinality relative to $X$ (and thus $X$ can have order type $\oone$),
we can prove that they are still equivalent to the principles stated
in section~\ref{sec:description-results}. This is left as an exercise
for the interested reader. 
\end{rem}

\subsection{Game Theory}
\label{sec:terminology-games}

We shall follow the convention that the first move of any game is move $0$, and that
the $k\Th$ move refers to move $k$. 
\emph{\textbf{Caution.}} There is
some potential for confusion, because this means for example
that the $1^{\mathrm{st}}$ move is move $1$, whereas the $0\Th$ move is really the
`first' move.

In the simplest $n$ player game, there is one winner and all of the other players
lose, and there are exactly $n$ possible outcomes. More generally, the game can
result in any possible ranking (not necessarily injective) of the $n$ players for a
total of
\begin{equation}
  \label{eq:37}
  \sum_{i=0}^{n-1}p(n,n-i)\cdot(n-i)!
\end{equation}
possible outcomes, where $p(n,k)$ is the number of partitions of $n$ into $k$ pieces. The first few
terms are $n!$, ${n \choose 2}(n-1)!$, $\left[{n \choose 3}+3{n\choose 4}\right](n-2)!$,\linebreak
$\left[{n\choose 4}+10{n\choose5}+15{n\choose 6}\right](n-3)!$, $\ldots$ .
We say that player $X$  \emph{wins}  if it is ranked strictly above all of
the other players. We say that it \emph{loses} if there is some other player ranked
strictly above him. Accordingly we distinguish between a \emph{winning
strategy} for player $X$ and a \emph{nonlosing strategy} for $X$. 
In the case of a two player game, 
there are $p(2,2)\cdot 2!+p(2,1)\cdot 1!=2+1=3$ possible outcomes, either player can
win and there can be a tie, i.e.~draw, where neither player wins or loses.

The following concept is useful for dealing with completeness systems.

\begin{defn}
\label{d-1}
Let $\game$ be a game of length $\delta$
and $X$ a player of $\game$. A \emph{forward strategy} for $X$
in the game $\game$ is a strategy $\Phi$ for $X$ such that for any position $P$ in
the game $\game$, where the game has not yet ended and it is $X$'s turn to play,
$\Phi(P)$ gives a move for $X$. A \emph{forward winning \tu[nonlosing\tu] strategy} for $X$ is a
forward strategy $\Phi$ for $X$ such that $X$ wins [does not lose]
the game so long as there exists
$\xi<\delta$ such that $X$ plays according to $\Phi$ on move $\eta$ for all
$\eta\geq\xi$.
\end{defn}

\noindent Thus the point of a forward strategy for $X$ is that $X$ can decide at any
point in the game to start playing according to the strategy. Of course a forward
(winning [nonlosing]) strategy is in particular a (winning [nonlosing]) strategy.

\begin{defn}
\label{d-9}
Let $\game(M,x,a_0,\dots,a_{n-1})$ be a parameterized game with a
fixed number of  players~$n$ (with respect to the parameters).
Suppose that the first parameter $M$ is taken from some family $\S$, the second
parameter is a subset of $M$ and
that the third parameter $\vec a=a_0,\dots,a_{n-1}$ is taken from some
family~$\T$.

We consider a function $F$ with domain $\S$ where
\begin{equation}
  \label{eq:38}
  F(M):\T\cap M\to\power(\power(M))\espc\text{for all $M\in\S$}, 
\end{equation}
i.e.~$F(M)(\vec a)\subseteq\power(M)$ for all $\vec a\in\T\cap M$.
We think of $F$ as describing 
a notion of suitability over $M$ for the second parameter. 
It is required to satisfy
\begin{equation}
  \label{eq:8}
  F(M)(\vec a)\neq\emptyset\espc\text{for all $M\in\S$ and all $\vec a\in\T\cap M$}.
\end{equation}

Suppose $X$ is some player in the parameterized game. 
For any
$\E\subseteq\S$, and any property $\Phi$ where $\Phi$ has $n+3$ variables,
we say that  $\Phi$ holds for $X$, \emph{$\E$-$F$-globally}, if 
$\Phi(X,M,x,\vec a)$ holds whenever
\begin{enumerate}[leftmargin=*, label=(\roman*), ref=\roman*, widest=iii]
\item\label{item:41} $M\in\E$,
\item\label{item:56} $x\in F(M)(\vec a)$,
\item\label{item:57} $\vec a\in\T\cap M$.
\end{enumerate}
Particularly, $X$ has a winning [nonlosing] strategy in the game
$\Game$, $\E$-$F$-globally, means that $X$ has a winning [nonlosing]
strategy in the game $\game(M,x,\vec a)$ whenever $M$, $x$ 
and $\vec a$ satisfy~\eqref{item:41}--\eqref{item:57}.

When we say that a property holds $F$-globally,
we mean $\S$-$F$-globally.
\end{defn}

We will consider parameterized games that are given by some
first order definition in the language of set theory. In some situations, e.g.~in
moving between forcing extensions, it is
important to distinguish between a parameterized game as
a mathematical object or a defined notion.\footnote{As it turned out,
  this never became an issue in the present article.}
The following notation is used to help accomplish this.

\begin{notn}
\label{notn:global}
Typically, when dealing with the parameterized game as a defined notion, we also
want the suitability function $F$ to be a defined notion, and vice versa. 
Thus we write a formula in place of $F$ to indicate that the parameterized
game, the function $F$ (and thus $\S$) and $\T$
are first order definable with no parameters. 
We may display objects to indicate that a particular variable is fixed, and also
when we want to specify some part as an object  (e.g.~$\E$ in example~\ref{x-5}). 
In the case when some of the parameters $b_0,\dots,b_{m-1}$ are fixed,
e.g.~some property holds globally for player $X$ in the game $\game(b_0,\dots,b_{m-1})$ 
(see e.g.~example~\ref{x-5}), 
equation~\eqref{eq:8} is only required for parameters $\vec a\in\T$ such that 
$\vec a=b_0,\dots,b_{m-1},a_m,\dots,a_{n-1}$ and for $M\ni\vec a$.
\end{notn}

\begin{example}
\label{x-4}
Let $\psi$ be a first order formula with $n+2$ free variables. 
When we say that a property $\Phi$ of the parameterized game $\game$ 
holds \emph{$\psi$-globally} for $X$ 
we are indicating that the parameterized game~$\game$ as well as $\S$ and $\T$
are definable by first order formulae with
no additional parameters and that the function $F$ is given by
\begin{equation}
  F(M)(\vec a)=\{x\subseteq M:\psi(M,x,\vec a)\}.\label{eq:41}
\end{equation}
Notice that equation~\eqref{eq:8} becomes
\begin{equation}
  \label{eq:39}
  \varphi_\S(M)\lands\varphi_\T(\vec a)\to \exists x\subseteq M\, \psi(M,x,\vec a)\espc
  \text{for all $M$ and all $\vec a\in M$}.
\end{equation}
\end{example}

\begin{example}
\label{x-5}
Letting $n=1$, suppose $\psi$ is a formula with $3$ free variables. 
Suppose $\E$ and $b$ are sets.
Then saying that $X$ has a winning strategy in the game $\game(b)$,
$\E$-$\psi$-globally, 
indicates that the parameterized game $\game(M,x,a)$ can be identified with a formula,
but only indicates that $X$ has a winning strategy when $a=b$.  
It also indicates that $F$ and $\T$ are definable with no parameters, 
but is specifying the object~$\E$. 
\end{example}

\begin{defn}
\label{definable:notn}
For two first order formulae $\varphi(v_0,\dots,v_n)$ and
$\psi(v_0,\dots,v_n)$, we write $\varphi\to\psi$ for its universal closure. 
We say that \emph{provably} $\varphi\to\psi$ to indicate that 
\begin{equation}
  \label{eq:32}
  \zfc\proves \ulc\forall v_0,\dots,v_n\spc
  \varphi(v_0,\dots,v_n)\to\psi(v_0,\dots,v_n)\urc.
\end{equation}
\end{defn}

Let us point out an obvious relationship.

\begin{prop}
\label{p-5}
Suppose $\varphi\to\psi$. 
Then if $X$ has a \tu[forward\tu]
\tu(winning\tu) strategy for $X$, $\E$-$\psi$-globally in a parameterized game\/ $\game$,
then $X$ also has a
\tu[forward\tu] \tu(winning\tu) strategy $\E$-$\varphi$-globally in\/ $\game$.
\end{prop}

\begin{defn}
\label{d-10}
Let $\game$ be some game of length $\delta$, with player $A$ playing first.
For a position $P$ in the game $\game$ with $A$ to play, we define the
restricted game $\game\restriction P$ as a game with the same players, playing in
the same order, of length $\delta-|P|$. The rules of the game are that in position
$Q$ of $\game\restriction P$ with $X$ to play, $m$ is a valid move for $X$ if $m$ is
a valid move for $X$ in the position $P\bigext Q$ of the game $\game$. 
And if $Q$ is the position at the end of the game $\game\restriction P$, 
then player $X$ wins [loses] in the play $Q$ 
if $X$ wins [loses] the play $P\bigext Q$ of the game $\game$. 
\end{defn}

The following lemma is used to obtain forward strategies. 

\begin{lem}
\label{l-16}
Assume that $\game(M,x,\vec a)$ is a parameterized game with a fixed number of players
where player $A$ always moves first, and
a fixed length $\delta$ with $\delta$ indecomposable, with respect to the
parameters. Suppose that
player $X$ has a winning \tu[nonlosing\tu] strategy 
in the parameterized game~$\game$, $\E$-$F$-globally. 
If every suitable triple $(M,x,\vec a)$ and every position $P$ of the
game $\game(M,x,\vec a)$ with $A$ to play, 
has an $x'\subseteq M$ and $\vec a'\in\T\cap M$ such that
\begin{samepage}
\begin{enumerate}[leftmargin=*, label=\tu{(\alph*)}, ref=\alph*, widest=b]
\item\label{item:13} $x'\in F(M)(\vec a')$,
\item\label{item:32} $\game(M,x',\vec a')=\game(M,x,\vec
  a)\restriction P$,
\end{enumerate}
\end{samepage}
then $X$ has a forward winning \tu[nonlosing\tu] strategy in\/~$\game$, $\E$-$F$-globally.
\end{lem}
\begin{proof}
We use the following notation in this proof.

\begin{notn}
\label{notn:inproof}
If $P$ is a position in some game $\game$, where $A$ moves first,
and it is player $X$'s turn to make its $\xi\Th$ move (note that this entails $|P|$ is
either $\xi$ or $\xi+1$ depending on whether or not $X=A$, resp.), 
then for each $\gamma\leq\xi$, we let $P_\gamma\sqsubseteq P$ 
be the position preceding $P$ where it is $A$'s turn to make its $\gamma\Th$ move. 
\end{notn}

Assuming the hypotheses, take a
suitable triple $(M,x,\vec a)$, and let $\Phi(M,x,\vec a)$ be a
winning [nonlosing] strategy for $X$ in the game $\game(M,x,\vec a)$. 
We define a forward strategy $\Phi'(M,x,\vec a)$ for $X$ as follows.
Take a position $P$ in the game $\game(M,x,\vec a)$ where it is $X$'s turn to make its $\xi\Th$ move.
For each $\gamma\leq\xi$,  let $x_{P,\gamma}$ and $\vec a_{P,\gamma}$ be the $x'$ and $\vec a'$ 
guaranteed by the hypothesis with $P:=P_\gamma$. 
Then by~\eqref{item:13} and~\eqref{item:32}, 
we can let $\gamma_{P}$ be the least ordinal $\gamma\leq\xi$ such that $P=P_\gamma\bigext Q^\gamma$ 
where $Q^\gamma$ is the result of $X$ playing according to the
strategy $\Phi(M,x_{P,\gamma},\vec a_{P,\gamma})$  in the game $\game(M,x_{P,\gamma},\vec a_{P,\gamma})$. And then
\begin{equation}
  \Phi'(M,x,\vec a)(P)=\Phi(M,x_{P,\gamma},\vec a_{P,\gamma_P})(Q^\gamma)\label{eq:19}
\end{equation}
defines a strategy for $X$ in the game $\game(M,x,\vec a)$ by~\eqref{item:32}, which is
moreover a forward strategy by its definition. 

Suppose that the game $\game(M,x,\vec a)$ has been played, where $X$ has
played according to $\Phi'(M,x,\vec a)$ on every move $\alpha\geq\xi$.
Assuming that $\xi$ is the least such ordinal, 
then by equation~\eqref{eq:19}, from its $\xi\Th$ move on, 
$X$ has played according to the strategy $\Phi(M,x_{P,\xi},\vec a_{P,\xi})$ 
in the game $\game(M,x_{P,\xi},\vec a_{P,\xi})$, where $P$ is the position in the game
$\game(M,x,\vec a)$ when it was $X$'s turn to make its $\xi\Th$ move. 
Therefore, $X$ wins [does not lose] in the game $\game(M,x_{P,\xi},\vec a_{P,\xi})$,
and thus $X$ wins
[does not lose] the game $\game(M,x,\vec a)$ by~\eqref{item:32}. 
\end{proof}

\begin{rem}
\label{r-7}
The only role of indecomposability is that it is entailed by~\eqref{item:32}.
\end{rem}
\section{General results}
\label{sec:general-results}
\subsection{The forcing notion}
\label{sec:forcing-notion}

\begin{defn}
\label{def:consistency}
For two families $\F$ and $\H$ of subsets of some ordinal, 
let $\poset(\F,\H)$ be the poset consisting of all pairs $p=(x_p,\X_p)$ where
 \begin{enumerate}[leftmargin=*, label=(\roman*), ref=\roman*, widest=iii]
   \item\label{item:34} $x_p\in\H$,
   \item\label{item:6} $\X_p$ is a nonempty countable family 
     of cofinal subsets of $(\F,\subseteqfnt)$, 
   \item\label{item:33} $\{y\in X:x_p\subseteq y\}$ is cofinal in $(\F,\subseteqfnt)$
     for all $X\in\X_p$,
 \save
 \end{enumerate}
 ordered by $q$ extends $p$ if
 \begin{enumerate}[leftmargin=*, label=(\roman*), ref=\roman*, widest=iii]
\restore
   \item $x_q\sqsupseteq x_p$ (i.e.~$x_q$ end extends $x_p$ with respect to the
     ordinal ordering),
   \item $\X_q\supseteq\X_p$.
 \end{enumerate}
We write $\poset(\H)$ for $\poset(\H,\H)$.
\end{defn}
\pagebreak

\begin{rem}
\label{r-5}
Note we can assume without loss of generality 
that $\H\subseteq\partial(\F)$, i.e.~%
\begin{equation}
  \label{eq:46}
  \poset(\F,\H)=\poset(\F,\H\cap\partial(\F)),
\end{equation}
where $\partial(\F)$ is the set of all $x$ such that
$\{y\in\F:x\subseteq y\}$ is cofinal in $(\F,\subseteqfnt)$.
Obviously $\partial(\F)\subseteq\downcl\F$.
\end{rem}

\begin{prop}
\label{p-1}
Let $(\F,\subseteqfnt)$ be $\lambda$-directed and $\H$ arbitrary. 
Then every $p\in\poset(\F,\H)$ and every $A\in[\F]^{<\lambda}$ has a\/ $y\in\F$
such that $x_p\subseteq y$ and $y$ is a $\subseteqfnt$-upper bound of $A$. 
\end{prop}
\begin{proof}
The ``nonempty'' in~\eqref{item:6} is needed here. Take any $X\in\X_p$.
There exists an upper bound $y'\in\F$ of $A$ because $\F$ is $\lambda$-directed,
and thus by~\eqref{item:33},
there exists $y\supseteqfnt y'$ with $x_p\subseteq y$ in $X$, as required. 
\end{proof}

\begin{defn}
\label{d-8}
Let $\F$ and $\H$ be families of subsets of some ordinal $\theta$. 
We say that $\H$ is \emph{$\F$-extendable} if for every $x\in \H$, 
every countable family $\X$ of $\subseteqfnt$-cofinal subsets of $\F$
such that $x\subseteq z$ for all $z\in X$ for all $X\in\X$,
and every $\xi<\theta$, there exists $y\subseteq\theta$ such that
\begin{enumerate}[leftmargin=*, label=(\roman*), ref=\roman*,
  widest=iii]
\item $y\sqsupseteq x$
\item $y\in\H$,
\item $y\setminus\xi\ne\emptyset$,
\item $\{z\in X:y\subseteq z\}$ is cofinal in $\F$ for all $X\in\X$.
\end{enumerate}
We just say that $\H$ is \emph{extendable} to indicate that it is $\H$-extendable. 
\end{defn}

\begin{rem}
\label{r-8}
The definition of $\H$ being $\F$-extendable actually depends on the choice of
$\theta$. We shall always implicitly assume that $\theta$ can be computed from $\F$,
i.e.~it is the supremum of the ordinals appearing in $\F$,
i.e.~$\theta=\sup\bigcup\F$.  
\end{rem}

\noindent The preceding definition was tailored for the following
density result that is crucial for obtaining the desired uncountable
$X\subseteq \theta$. 

\begin{prop}
\label{p-12}
Let $\F$ and\/ $\H$ be families of subsets of\/ $\theta$.
Suppose that $\H$ is $\F$-extendable. Then for every\/ $\xi<\theta$,
\begin{equation}
  \label{eq:2}
  \D_\xi=\{p\in\poset(\F,\H):\sup(x_p)\geq\xi\}\textup{ is dense}. 
\end{equation}
\end{prop}

For a filter $G\subseteq\poset(\F,\H)$, the generic object is
\begin{equation}
  \label{eq:33}
  X_G=\bigcup_{p\in G}x_p\subseteq\theta.
\end{equation}
We make the obvious observations.

\begin{prop}
\label{p-30}
$X_G$ is the union of a chain in $(\H,\sqsubseteq)$.
\end{prop}

\begin{prop}
\label{p-31}
Every proper initial segment $y\sqsubset X_G$ has an $x\sqsupseteq y$ in $\H$. 
I.e.~every proper initial segment of $X_G$ is in $\downcl{(\H,\sqsubseteq)}$.
\end{prop}
\subsection{The associated games}
\label{sec:associated-games}

Following is the natural game associated with our forcing notion.

\begin{defn}
\label{d-11}
Let $\F$ and $\H$ be families of subsets of an ordinal $\theta$, 
$y\subseteq\theta$
and $p\in\poset(\F,\H)$ with $x_p\subseteq\nobreak y$. 
Define the game $\gcomp(y,\F,\H,p)$ 
with players \emph{Extender} and \emph{Complete} of length~$\omega$.
Extender plays first and on move $0$ must play $p_0$ so that
\begin{itemize}[leftmargin=*]
\item $p_0$ extends $p$.
\end{itemize}
On the $k\Th$ move:
\begin{itemize}[leftmargin=*]
\item Extender plays $p_k\in\poset(\F,\H)$ satisfying
  \begin{enumerate}[leftmargin=*]
  \item\label{item:14}  $p_k$ extends $p_i$ for all $i=0,\dots,k-1$,
  \item\label{item:15}  $x_{p_k}\subseteq y\setminus\bigcup_{i=0}^{k-1} s_i$,
  \end{enumerate}
\item Complete plays a finite $s_k\subseteq y\setminus x_{p_k}$.
\end{itemize}
Complete wins if the sequence $p_k$ ($k<\omega$) has a common extension in
$\poset(\F,\H)$, and Extender wins otherwise.
\end{defn}

\begin{notn}
\label{notn:generate}
For a centered subset $C$ of some poset $P$, 
we let $\<C\>$ denote the filter on $P$
generated by $C$.
\end{notn}

\begin{defn}
\label{d-13}
We define a variation of the game $\gcomp$ as follows.
For some $M$, for families $\F,\H\in M$ of subsets of some ordinal, 
$y\subseteq M$
and $p\in\poset(\F,\H)\cap M$ with $x_p\subseteq y$, we define the game
$\ggen(M,y,\F,\H,p)$. It has the same
rules, but with the additional rule that Extender's $k\Th$ move, $p_k$ or
$(p_k,X_k)$, respectively, must satisfy
\begin{equation}
  \label{eq:35}
  p_k\in M.
\end{equation}
This game has three possible outcomes, determined as follows:
\begin{enumerate}[leftmargin=*, label=(\roman*), ref=\roman*, widest=iii]
\item\label{item:35} 
  Extender loses (i.e.~Complete wins) if
  $\<p_k:k<\omega\>\notin\Gen(M,\poset(\F,\H))$,
\item\label{item:36}
  the game is drawn (i.e.~a tie) if $\<p_k:k<\omega\>\in\Genc(M,\poset(\F,\H))$,
\item\label{item:7} Extender wins the game
if $\<p_k:k<\omega\>\in\Gen(M,\poset(\F,\H))$ 
  but~\eqref{item:36} fails.
\end{enumerate}
\end{defn}

The game $\ggen(M,y,\F,\H,p)$ is especially interesting for us 
because a draw in this game 
corresponds precisely with complete genericity.

\begin{prop}
\label{p-4}
Let\/ $p_k$ denote Extender's\/ $k\Th$ move in the game\/ 
$\ggen(M,\allowbreak y,\allowbreak\F,\allowbreak\H,\allowbreak p)$.
Then the game results in a draw iff\/ $\<p_k:k<\omega\>
\in\Genc(M,\allowbreak\poset(\F,\H),\allowbreak p)$.
\end{prop}

The following augmented game is used for preserving
the nonspecialness of trees. 

\begin{defn}
\label{d-12}
We define an augmented game $\ggen^*(M,y,\F,\H,p)$. 
Again Extender moves first with $p_0\ge p$ in $M$, 
but Extender additionally plays $X_k\subseteq\F$ on each move. 
The whole point is that $X_k$ is \emph{not} required to be in $M$.
On the $k\Th$ move:\pagebreak
\begin{itemize}[leftmargin=*, label=\altbullet]
\item Extender plays $(p_k,X_k)$ where $p_k\in \poset(\F,\H)\cap M$  
and $X_k\subseteq\F$ satisfy
  \begin{enumerate}[leftmargin=*]
  \item\label{item:22} $p_k$ extends $p_i$ for all $i=0,\dots,k-1$,
  \item\label{item:23} $x_{p_k}\subseteq y\setminus\bigcup_{i=0}^{k-1} s_i$,
  \item\label{item:28} $\{x\in X_i:x_{p_k}\subseteq x\}$ is cofinal in
    $(\F,\subseteqfnt)$ for all $i=0,\dots,k$,
  \end{enumerate}
\item Complete plays a finite $s_k\subseteq y\setminus x_{p_k}$. 
\end{itemize}
The possible outcomes are:
\begin{enumerate}[leftmargin=*, label=(\roman*), ref=\roman*, widest=iii]
\item\label{item:35a} 
  Complete wins if $\<p_k:k<\omega\>\notin\Gen(M,\poset(\F,\H))$.
\item\label{item:36a}
  The game is drawn if $\<p_k:k<\omega\>\in\Genc(M,\poset(\F,\H))$ and
  moreover $\{p_k:k<\omega\}$ has a common extension $q\in\poset(\F,\H)$ with 
\begin{itemize}[leftmargin=*, label=\altbullet]
\item $\{X_k:k<\omega\}\subseteq\X_q$,
\end{itemize} 
\item\label{item:7a} Complete loses the game 
if $\<p_k:k<\omega\>\in\Gen(M,\poset(\F,\H))$ 
  and~\eqref{item:36a} fails.
\end{enumerate}
\end{defn}

\begin{prop}
\label{p-6}
At the end of any of the three games $\gcomp$, $\ggen$ or $\ggen^*$,
\begin{equation}
  \label{eq:23}
  \bigcup_{k<\omega}x_{p_k}\subseteq y\setminus\bigcup_{k<\omega}s_k.
\end{equation}
\end{prop}

The augmented game $\ggen^*$ `includes' 
the game $\ggen$ in the following sense.

\begin{prop}
\label{p-7}
If $\bigl((p_0,X_0),s_0\bigr),\dots,\bigl((p_k,X_k),s_k\bigr)$ is a position in the
game $\ggen^*(y,\F,\H,p)$ then $(p_0,s_0),\dots,(p_k,s_k)$ is a position in the
game $\ggen(y,\F,\allowbreak\H,p)$. 
Conversely, if $(p_0,s_0),\dots,(p_k,s_k)$ is a position
in the game $\ggen(y,\F,\allowbreak\H,p)$ then
$\bigl((p_0,\F),s_0\bigr),\dots,\bigl((p_k,\F),s_k\bigr)$ is a position in the game
$\ggen^*(y,\F,\allowbreak\H,p)$. 
\end{prop}
\begin{proof}
$\X_p\neq\emptyset$ in definition~\ref{def:consistency}\eqref{item:6} 
is used for the converse.
\end{proof}

\begin{prop}
\label{p-25}
A nonlosing strategy for Complete in the game $\ggen^*(M,\allowbreak
y,\allowbreak \F,\allowbreak \H,\allowbreak p)$
yields a nonlosing strategy in the game $\ggen(M,y,\F,\H,p)$.
\end{prop}
\begin{proof}
At a position $(p_0,s_0),\dots,p_k$ (with Complete to move) 
in the game\linebreak $\ggen(y,\F,\H,\allowbreak p)$, 
by proposition~\ref{p-7}, $\bigl((p_0,\F),s_0\bigr),\dots,(p_k,\F)$ is a position in
the game $\ggen^*(y,\F,\H,p)$. Thus the strategy for Complete in the game
$\ggen^*(y,\allowbreak\F,\allowbreak\H,\allowbreak p)$ gives a move $s_k$. 
This defines a nonlosing strategy for Complete
in the game $\ggen(y,\F,\H,p)$, because a draw  our play of the game
$\ggen^*(y,\F,\H,p)$ results in  draw in the game $\ggen(y,\F,\H,p)$.
\end{proof}

We relate the `completeness' game to the latter `genericity' game. 

\begin{prop}
\label{p-13}
A \tu(forward\tu) winning strategy for Complete in the game $\gcomp(y,\F,\H,p)$ 
yields a \tu(forward\tu) nonlosing strategy for Complete in the game
$\ggen^*(M,y,\F,\H,p)$.
\end{prop}
\begin{proof}
At a position $\bigl((p_0,X_0),s_0\bigr),\dots,(p_k,X_k)$ in the game
$\ggen^*(M,y,\F,\H,p)$, let $\bar p_i=p_i\cup \bigcup_{j=0}^i X_j$ for
each $i=0,\dots,k$. Then each $\bar p_i\in\poset(\F,\H)$ by
rule~\eqref{item:28}, and thus $(\bar p_0,s_0),\dots,\bar p_k$ is a
position in the game $\gcomp(y,\F,\H,p)$. Thus if $s_k$ is played
according to Complete's winning strategy in the latter game, then
Complete wins the latter game. This means that $\{\bar p_0,\bar
p_1,\dots\}$ has a common extension $q$, and then $q$ extends  
$\{p_0,p_1,\dots\}$ and satisfies $\{X_0,X_1,\cdots\}\subseteq\X_q$,
yielding a draw in the former game.
\end{proof}

\begin{prop}
\label{p-9}
Assume $\F$ and $\H$ are families of subsets of\/ $\theta$, 
$p\in\poset(\F,\H)$ and $x_p\subseteq y$. Let\/ $t\in\Fin(\theta\setminus x_p)$. 
Then Complete has a winning strategy in the game\/ $\gcomp(y,\F,\H,p)$ 
iff it has a winning strategy in the game\/ 
$\gcomp(y\cup\nobreak t,\allowbreak\F,\allowbreak\H,\allowbreak p)$.
Similarly for the games\/ $\ggen$ and\/ $\ggen^*$.
\end{prop}
\begin{proof}
Assume Complete has a winning strategy in $\gcomp(y,\F,\H,p)$. 
Assume without loss of generality that $y\cap t=\emptyset$. 
Then a winning strategy for Complete in $\gcomp(y\cup t,\F,\H,p)$ 
is given by playing $s_k\cup t$ on move $k$ 
where $s_k$ is played according to the strategy for the former game. 
This is because they both give identical restrictions on Extender's choice of
moves according to rule~\eqref{item:15}. Conversely, 
if Complete plays $s_k'\setminus t$ in the former game, 
where $s_k'$ has been played in the latter game, then Extender has
less freedom to move in the former game. 
\end{proof}

We are interested in making finite extensions of the third parameter when dealing
with completeness systems; but unfortunately, 
the above approach does not seem to generalize to forward winning strategies. 

The argumentation in the preceding proof, 
i.e.~the fact that restricting Extender's moves is 
favourable for Complete in any of the $3$ games, does show the following.

\begin{prop}
\label{p-10}
Let $\Phi$ be a \tu(forward\tu) winning, resp.~nonlosing, strategy for Complete 
in the game $\gcomp(y,\F,\H,p)$ or the game $\ggen(M,y,\F,\H,p)$, respectively.
Then Complete wins, resp.~does not lose,
the game whenever it plays $s_k\supseteq\Phi(P_k)$ on every move
\tu(after some point\tu) in the game, where $P_k$ is the position after Extender
makes its\/ $k\Th$ move. Similarly, for the augmented game $\ggen^*$.
\end{prop}

Next we isolate the role played by the $\subseteqfnt$-cofinal 
subsets of the family~$\F$.

\begin{rem}
\label{r-13}
Henceforth, when we write $H_\kappa$ there is a tacit assumption that $\kappa$ 
is a sufficiently large regular cardinal for the argument at hand. 
This will always be in the context of some pair $(\F,\H)$ of
subfamilies of $\power(\theta)$. It will always be sufficiently large
as long as $\poset(\F,\H)\in H_\kappa$, e.g.~when
$\kappa\ge\bigl(2^{\max\{|\F|,|\H|\}}\bigr){}^+$. 
\end{rem}

\begin{defn}
We say that a family $\F$ of sets is \emph{closed under finite
  reductions} to indicate that it is 
closed under finite set subtraction, i.e.~$x\setminus s\in\F$ for all $x\in\F$ 
and all finite $s\subseteq x$. 
\end{defn}

\begin{lem}
\label{l-12}
Let $\F$ and $\H$ be subsets of $\power(\theta)$, with $\F$ closed under finite
reductions, and let $M\prec H_\kappa$ be a model containing $\F$ and $\H$ satisfying
\begin{equation}
  \label{eq:11}
  x\subseteq M\espc\textup{for all $x\in\F\cap M$.}
\end{equation}
Suppose $p\in\poset(\F,\H)\cap M$,
let $Q\subseteq\poset(\F,\H)$ be an element of $M$
and assume that $y\subseteq\theta$ is a $\subseteqfnt$-upper bound of $\F\cap M$
with $x_p\subseteq y$.
Then one of the following two alternatives must hold.\textup{
\begin{enumerate}[leftmargin=*, label=(\alph*), ref=\alph*, widest=b]
\item\label{item:26} \textit{There exists an extension $q$ of\/ $p$ in $Q\cap M$ with
    $x_q\subseteq y$.}
\item\label{item:27} 
  \textit{There exists a $\subseteqfnt$-cofinal $X\subseteq\F$ such that}
  \begin{enumerate}[leftmargin=*, label=(\arabic*), ref=\arabic*, widest=2]
  \item\label{item:29} \textit{$x_{p}\subseteq x$ for all $x\in X$,}
  \item\label{item:30} \textit{for no extension $q\in Q$ of\, 
      $p$ does there exist $z\in X$ satisfying $x_q\subseteq z$.}  
  \end{enumerate}
\end{enumerate}}
\end{lem}
\begin{proof}
Define $Y$ to be the set of all $x\in\F$ for which 
there is some  $y_x\supseteqfnt x$ in $\F$, with $x_{p}\subseteq y_x$,
such that no $q\geq p$ in $Q$ satisfies $x_q\subseteq y_x$. 
Clearly $Y\in M$. Take $x\in \F\cap M$. Taking any $Z\in\X_p$, 
by elementarity there exists
$z\in Z\cap M$ such that $x_{p}\subseteq z$ and $x\subseteqfnt z$. 
Since $y\supseteqfnt z$, 
\begin{equation}
\label{eq:28}
y_x=y\cap z\supseteqfnt x.
\end{equation}
And $y_x$ is the result of removing a finite subset from $z$. 
Thus $y_x\in M$ as $z\subseteq M$, and $y_x\in\F$ by the assumption on $\F$.

Assume that alternative~\eqref{item:26} fails.
Then as $y_x\subseteq y$ is in $M$, by
elementarity, $y_x$ witnesses that $x\in Y$. 
Therefore, $Y=\F$ by elementarity, and thus $X=\{y_x:\nobreak x\in\nobreak\F\}$ 
is $\subseteqfnt$-cofinal by~\eqref{eq:28}. 
And since the $y_x$'s are witnesses, alternative~\eqref{item:27} holds for $X$.
\end{proof}

\begin{cor}{0}
\label{l-9}
Let $\F$, $\H$, $M$, $p$, $Q$ and $y$ all be as specified in 
lemma~\tu{\ref{l-12}}.
Assume that $k+1$ moves have been made 
in either the game $\gcomp(y,\F,\H,p)$, or the
game $\ggen(M,y,\F,\H,p)$ with $y\subseteq M$, 
with Extender playing $p_i$ on its $i\Th$ move, and that each $p_i\in M$ 
\tu(in the former game\tu).
Then one of the following two alternatives must hold.\textup{
\begin{enumerate}[leftmargin=*, label=(\alph*), ref=\alph*, widest=b]
\item\label{item:16} \textit{Extender has a move with $p_{k+1}\in Q\cap M$,}
\item\label{item:17}
  \textit{There exists a $\subseteqfnt$-cofinal $X\subseteq\F$ such that}
  \begin{enumerate}[leftmargin=*, label=(\arabic*), ref=\arabic*, widest=2]
  \item\label{item:18} \textit{$x_{p_k}\subseteq x$ for all $x\in X$,}
  \item\label{item:19} \textit{for no extension $q\in Q$ of\, 
      $p_{k}$ does there exist $z\in X$ satisfying $x_q\subseteq z$.}  
  \end{enumerate}
\end{enumerate}}
\noindent Similarly for the augmented game  $\ggen^*(M,y,\F,\H,p)$.
\end{cor}
\begin{proof}
Let $s_0,\dots,s_k$ denote the moves played so far by Complete. 
Lemma~\ref{l-12} is applied with $p:=p_k$ 
and $y:=y\setminus\bigcup_{i=0}^k s_i$.
The second alternatives are identical, and thus if~\eqref{item:17}
fails, then there
is an extension $p_{k+1}\geq p_k$ in $Q\cap M$ with $x_{p_{k+1}}\subseteq
y\setminus\bigcup_{i=0}^k s_i$. 
Thus $p_{k+1}$ satisfies the
requirement~\eqref{item:15} of the game, as needed. 
\end{proof}

The main purpose the side condition $\X_p$ is to allow Extender to play 
inside a given dense subset of $\poset(\F,\H)$ in $M$.

\begin{cor}{0}
\label{o-2}
In the situation of corollary~\tu{\ref{l-9}},
if $Q$ is dense then Extender can always play $p_{k+1}\in Q\cap M$.
\end{cor}
\begin{proof}
Supposing towards a contradiction that Extender 
has no move with\linebreak $p_{k+1}\in Q\cap M$, by corollary~\ref{l-9}, 
there is a cofinal $X\subseteq\F$ as in alternative~\eqref{item:17}.
Then $\bar q=(x_{p_k},\X_{p_k}\cup\{X\})$ is a condition of $\poset(\F,\H)$ 
by~(\ref{item:27}\ref{item:29}),
and there exists $q\geq\bar q$ in $Q$ by density. But $x_q\subseteq z$ for cofinally many
$z\in X$ contradicting~(\ref{item:27}\ref{item:30}).
\end{proof}

\begin{cor}{0}
\label{o-4}
Let $\F$, $\H$, $M$, $p$, $Q$ and $y$ be as in lemma~\tu{\ref{l-12}}, with moreover
$M$ countable and $y\subseteq M$. Then Extender has nonlosing strategies in both of
the games $\ggen(M,y,\F,\H,p)$ and\/ $\ggen^*(M,y,\F,\H,p)$.
\end{cor}
\begin{proof}
Let $(D_k:k<\omega)$ enumerate all of the dense subsets of $\poset(\F,\H)$ in~$M$.
By corollary~\ref{o-2}, Extender can always make move $k$ with
\begin{equation}
\label{eq:40}
  p_{k}\in D_k\cap M.
\end{equation}
This describes a nonlosing strategy, because at the end of the game,
$\<p_k:k<\nobreak\omega\>\in\Gen(M,\poset(\F,\H))$. 
\end{proof}

\begin{defn}
\label{d-14}
Let $\psimin(M,y,\F,\H,p)$ be a formula expressing the conjunction of
\begin{enumerate}[leftmargin=*, label=(\roman*), widest=ii]
\item $x_p\subseteq y$,
\item $y$ is an upper bound of $(\F\cap M,\subseteqfnt)$. 
\end{enumerate}
\end{defn}

Later on we will use the fact that $\psimin(M,\cdot,\F,\H,p)$ 
defines a set that is second order definable over $M$.

\begin{defn}
\label{d-3}
Let $\phimin(y;\F,\H,p)$ be a second order formula expressing the conjunction of
\begin{enumerate}[leftmargin=*, label=(\roman*), widest=ii]
\item $x_p\subseteq y$; formally, $\forall \alpha\in x_p\spc y(\alpha)$,
\item $y$ is an upper bound of $(\F,\subseteqfnt)$; formally, 
$\forall x\in\F\spc y(\alpha)$ for all but finitely many $\alpha\in x$.
\end{enumerate}
\end{defn}

\begin{prop}
\label{p-11}
Suppose $M$ is a model of enough of 
$\zfc-\mathrm{P}$\footnote{$\zfc$ minus the Power Set axiom.} 
and $x\subseteq M$ for all $x\in\F$.
Then for all\/ $\F,\H,p\in M$ and all\/ $y\subseteq M$, 
$\psimin(M,y,\F,\H,p)\leftrightarrow M\models\phimin(y;\F,\H,p)$. 
\end{prop}

\begin{defn}
\label{d-15}
All three of the games considered, $\gcomp$,  $\ggen$ and $\ggen^*$, are
viewed as parameterized games of the form
$\game(M,x,a_0,a_1,a_2,a_3)$, as in
definition~\ref{d-9}, where $a_3$ is a ``dummy''
variable whose purpose is explained below. For example, 
$\gcomp(M,x,a_0,a_1,a_2,a_3)\equiv\gcomp(x,a_0,a_1,a_2)$ 
and $\ggen(M,\allowbreak x,\allowbreak a_0,\allowbreak a_1,\allowbreak
a_2,\allowbreak a_3)\equiv\ggen(M,x,a_0,a_1,a_2)$. 
Define $\S$ to be the class of all countable elementary submodels $M\prec H_\kappa$,  
with $\kappa$ a regular uncountable cardinal. For a given cardinal $\kappa$, we define $\S_\kappa\subseteq\S$ by 
$\S=\bigcup_{\mu\geq\kappa\text{ is regular}}\{M\prec
H_\mu:|M|=\aleph_0\}$.
 $\T$ is defined by $\varphi_\T(\F,\H,p,a_3)$ stating that $\ulc\F$ and $\H$ are families of
sets of ordinals, $p\in\poset(\F,\H)$ and $a_3=\poset(\F,\H)\urc$, 
with the provision that we may restrict $\T$ further when needed. 
Note that these games, as well as $\S$ and $\T$, 
are definable without any additional parameters. In this setting, we use a formula
$\psi$ to describe the suitability function $F$;
we suppress the last free
variable in $\psi$ since $a_3$ obviously plays no role in the definition of
$F$. The role played by $a_3$,
is that for any $M\in\S$, when $\vec a=(\F,\H,p,a_3)\in\T\cap M$ this implies that
$\poset(\F,\H)\in M$ and thus $M\prec H_\kappa$ for some sufficiently
large cardinal~$\kappa$ as in remark~\ref{r-13}. We could also use
$\S_\kappa$ below instead of $\S$, just as well. 

For $\E\subseteq\S$, $\E_\kappa=\E\cap\S_\kappa$.
Moreover, for $R\subseteq\theta$, we let 
\begin{equation}
  \label{eq:16}
  \E(R,\theta)=\{M\in\S:\sup(\theta\cap M)\in R\},
\end{equation}
and for $R\subseteq\oone$ we write $\E(R)$ for $\E(R,\oone)$. 
Thus $\E(R)=\{M\in\S:\delta_M\in R\}$
(cf.~\eqref{eq:4}).
Also $\E_\kappa(R,\theta)=\E(R,\theta)\cap\S_\kappa$. 
\end{defn}

\begin{example}
\label{x-7}
Let $R\subseteq\oone$ and suppose $\F,\H$ are families of sets of
ordinals. Suppose that $\psi(v_0,\dots,v_4)$ is a first order formula
such that for every $M\in\E(R)$ containing $\F$ and $\H$, and every
$p\in\poset(\F,\H)\cap M$, $\varphi_\T(\F,\H,p,a_3)$ implies there is
a $y\subseteq
M$ such that $\psi(M,y,\F,\H,p)$ holds, 
and thus equation~\eqref{eq:39} is satisfied. 
Then saying Complete has a winning strategy $\E(R)$-$\psi$-globally 
in the game $\gcomp(\F,\H)$,  means that it has a winning strategy in the game
$\gcomp(y,\F,\H,p)$ for all $M\in\S$ with $\F,\H,\poset(\F,\H)\in M$ and $\delta_M\in R$,
and all $p\in\poset(\F,\H)\cap M$. 
Alternatively, we could have omitted $a_3$ and equivalently referred
to ``$\E_\kappa(R)$-$\psi$-globally'' instead.
\end{example}

\begin{cor}[5]{1}
\label{o-6}
Restricting $\T$ in definition~\tu{\ref{d-15}} to only include
families $\F$ of countable sets of ordinals, 
Extender has nonlosing strategies in the games\/ $\ggen$
and\/ $\ggen^*$, $\psimin$-globally.
\end{cor}
\begin{proof}
First we have to show that equation~\eqref{eq:39} holds. But if
$\F,\H,p\in M$ and $\varphi_T(\F,\H,p,a_3)$, then in particular, $\F$
is a family of countable subsets of some ordinal $\theta$ and
$p\in\poset(\F,\H)$. 
Since members of $\F$ are countable, $x\subseteq M$ for all
$x\in\F\cap M$ and thus $\psimin(M,\theta\cap M,\F,\H,p)$ holds. 

Now we obtain a nonlosing strategy for Extender in both of the
games\linebreak\ 
$\ggen(M,\allowbreak\theta\cap M,\allowbreak\F,\allowbreak\H,\allowbreak p)$ 
and $\ggen^*(M,\theta\cap M,\F,\H,p)$ by
corollary~\ref{o-4}, since~\eqref{eq:11} holds. 
\end{proof}
\subsection{Complete properness}
\label{sec:complete-properness}

\begin{cor}[5]{0}
\label{o-5}
Let $\F$ and $\H$ be subfamilies of\/ $[\theta]^{\leq\aleph_0}$, 
with $\F$ closed under finite reductions.
Suppose that $\psi\to\psimin$,  and that $\psi$-globally, 
Extender has \emph{no} winning strategy 
in the parameterized game $\ggen(\F,\H)$.
Then\/ $\poset(\F,\H)$ is completely proper.
\end{cor}
\begin{proof}
Suppose $M\prec H_\kappa$ is a countable elementary submodel with $(\F,\H)\in M$.
Take $p\in\poset(\F,\H)\cap M$. 
The game $\ggen(M,\F,\H,y,p)$ is played 
with Extender playing $p_k$ on move $k$
according to a nonlosing strategy, which it has by corollary~\ref{o-6}
and proposition~\ref{p-5} since $\psi\to\psimin$. 
By the hypothesis that Extender's strategy in
the game $\ggen(M,\F,\H,y,p)$ is not a winning strategy, 
Complete can play in such a way that
the game does not result is a win for Extender.   
Thus the result is a drawn game, and hence $\<p_k:k<\omega\>\in\Genc(M,\poset(\F,\H),p)$
by proposition~\ref{p-4}, as required.
\end{proof}

The following weaker result gives a purely combinatorial characterization 
of complete properness, unlike corollary~\ref{o-5}.

\begin{cor}[5]{0}
\label{c-1}
Let $\F$ and $\H$ be subfamilies of\/ $[\theta]^{\leq\aleph_0}$, with $\F$ closed under
finite reductions.
Suppose that $\psi\to\psimin$, and that $\psi$-globally, 
Complete has a\/ winning strategy in
the parameterized game $\gcomp(\F,\H)$.
Then\/ $\poset(\F,\H)$ is completely proper.
\end{cor}
\begin{proof}
By propositions~\ref{p-25} and~\ref{p-13}, Complete has a nonlosing
strategy in the game $\ggen(\F,\H)$, $\psi$-globally, and in
particular, corollary~\ref{o-5} applies.
\end{proof}

\subsection{Completeness systems}
\label{sec:completeness-system}

Our formulation of \emph{completeness systems} differs slightly
from that observed in the literature. 
Completeness systems were invented by Shelah, 
and we use the same underlying ideas  
as in the original formulation in~\cite{MR675955}.

A full account of the
theory of completeness systems is given by Abraham in~\cite{hbst}. 
It is emphasized there that in order to apply the theory, a
$P_\alpha$-name $\dot Q_\alpha$ for 
a poset must be complete for some completeness
system that lies in the ground model. Then \emph{simple}---meaning simply
definable---completeness systems are introduced to achieve this. 
An alternative to completeness systems that has gained some popularity
was introduced in~\cite{MR1638230},
where the necessary combinatorial properties of $\dot Q_\alpha$ 
entailed by the completeness system are isolated. 
In this approach one directly verifies that the
name $\dot Q_\alpha$ itself satisfies the prerequisite  properties.

Our approach is less robust than in~\cite{hbst}; however, we do not know of any
examples not encompassed by our treatment,\footnote{We only deal with
  $\sigma$-complete systems, but our treatment could be adapted to
  more restrictive systems (e.g.~\cite[Ch.~VIII,~\Section4]{MR1623206}).} 
and it may have
some advantages, including we hope, conceptual simplicity. 
In our formulation, the fundamental notion is a second order formula 
rather than the system of filters it describes. 
Moreover, the completeness system
(i.e.~this formula) is good for exactly one class of posets. 
This captures every usage of the completeness systems that we have observed, 
although there may very well be  uses for undefinable (i.e.~nonsimple) 
completeness systems, or systems that work for more than one class of posets. 
A potential advantage of our approach is that, 
when the formula provably (in~$\zfc$) has the required properties, the
completeness system functions in arbitrary forcing extensions, 
whereas (as indicated in~\cite{hbst}) the approach of using 
a ground model system is only valid in forcing
extensions that do not add new countable subsets of the ground model. 

\begin{defn}
\label{d-6}
We say that a pair of formulae $\wp(v_0,\dots,v_{n-1})$ and
$\tau(v_0,\dots,\allowbreak v_{n-1})$ 
in the language of set theory \emph{describe} a poset over some model $N$ 
if $N\models\nobreak\ulc\forall x_0,\dots,x_{n-1}\spc\tau(\vec x)\to
\wp(\vec x)$ is a poset$\urc$.\footnote{$\tau$ is 
unnecessary but is used for presentational purposes.} 
If we are working with some ground model $V$ and  we say
that $(\wp,\tau)$ describes a poset, we mean that it describes it over $V$.
And the pair \emph{provably describes} a poset 
if $\zfc\proves\ulc\forall \vec x\spc\tau(\vec x)\to \wp(\vec x)$ 
is a poset$\urc$.
\end{defn}

\begin{example}
\label{x-2}
Let $\tau(v_0,v_1)$ express $\ulc v_0$ and $v_1$ consist of sets of ordinals$\urc$.
Then the pair $(\poset,\tau)$ provably describes a poset, where
$\poset$ is from definition~\ref{def:consistency}, i.e.~for $(\F,\H)$ satisfying
$\tau(\F,\H)$, the described poset is $\poset(\F,\H)$.
\end{example}

\begin{defn}
\label{d-2}
Suppose $(\wp,\tau)$ is a pair of formulae with $n$ free variables 
that (provably) describes a poset.
A (\emph{provable}) \emph{completeness system} for $(\wp,\tau)$ 
will refer to  a second order formula
$\varphi(Y_0,Y_1;\allowbreak v_0,\dots,v_n)$ 
for which (it is provable in $\zfc$ that): 
for every countable $M\prec H_\kappa$, 
where $\kappa$ is a sufficiently large regular cardinal, 
for all $\vec a\in M$ satisfying $\tau(\vec a)$, 
for every $p\in \wp(\vec a)^M$ (i.e.~$M\models p\in \wp(\vec a)$), 
the family of sets
\begin{equation}
  \label{eq:22}
  \G_Z=\{G\subseteq\Gen(M,\wp(\vec a),p):M\models\varphi(G,Z;\vec a,p)\}
  \espc\text{($Z\subseteq M$)}
\end{equation}
\begin{enumerate}[leftmargin=*, label=(\roman*), widest=ii]
\item\label{item:31} generates a proper filter on $\Gen(M,\wp(\vec a),p)$, 
i.e.~every finite intersection 
$\G_{Z_0}\cap\dots\cap\G_{Z_{n-1}}\subseteq\Gen(M,\wp(\vec a),p)$ is
nonempty,
\item\label{item:21} has a member that is a 
subset of $\Genc(M,\wp(\vec a),p)$, i.e.~there exists $Z\subseteq M$ such that
every element $G\in\G_Z$ has a common extension in~$\wp(\vec a)$.
\end{enumerate}
The completeness system is called \emph{$\sigma$-complete} 
if (it is provable in $\zfc$ that) for all $M$, $\vec a$ and $p$ as above, 
the filter generated by the family
from equation~\eqref{eq:22} is $\sigma$-complete.
\end{defn}

\begin{rem}
\label{r-4}
To avoid confusion, it should be noted that in the typical formulation
from the literature condition~\ref{item:21} is stipulated by stating
that the poset is complete for the given completeness system.
\end{rem}

We use $\aupalpha$-properness together with completeness systems in
what is now a standard method of forcing without adding reals. Since
we have made some adjustments to the usual terminology for
completeness systems, the following theorem needs to be taken
in the present context.

\begin{term}
Let $\vec P=(P_\xi,\dot Q_\xi:\xi<\mu)$ be an iterated forcing
construct. 

When we say that an \emph{iterand} $\dot Q_\xi$ of $\vec P$
satisfies some property $\Phi$ we of course mean that
$P_\xi\forces\Phi(\dot Q_\xi)$. 

As usual, when we say that $\vec P$ has \emph{countable supports} we
mean that $P_\delta=\varprojlim{}_{\xi<\delta}P_\xi$ for limit
$\delta$ of countable cofinality, and
$P_\delta=\varinjlim{}_{\xi<\delta}P_\xi$ for limit $\delta$ of uncountable
cofinality. This also determines $P_\mu$ for $\vec P$ when $\mu$ is a
limit, and of course $P_\mu=P_{\mu-1}\gdot\dot Q_{\mu-1}$ in case $\mu$
is a successor. 
\end{term}

\begin{thmo}[Shelah]
\label{Shelah:alpha}
Let $\vec P$ be an iterated forcing construction with countable
supports. 
Suppose that $\E\subseteq[H_\kappa]^{\aleph_0}$ is
stationary for some sufficiently large regular cardinal $\kappa$\tu;
and suppose for each $\xi$, $(\wp_\xi,\tau_\xi)$ provably describes a
poset and has a $\sigma$-complete completeness system $\varphi_\xi$.
If for each $\xi<\mu$, the iterand\/ $\dot Q_\xi$ is $\E$-$\aalpha$-proper, 
and $\dot Q_\xi=\wp_\xi(\dot{\vec a})$ and $\tau_\xi(\dot{\vec a})$ hold for some
parameter~$\dot {\vec a}$, then $P_{\length(\vec P)}$ does not add new reals.
\end{thmo}

\begin{defn}
\label{d-17}
We let \emph{$\mathbb D$-complete} 
denote the class of all posets $Q$
for which there exists $(\wp,\tau)$ provably describing a poset, such
that $(\wp,\tau)$ has a $\sigma$-complete completeness system and
$Q=\wp(\vec a)$ for some parameter $\vec a$ satisfying $\tau(\vec a)$. 
Then the forcing axiom $\ma(${\deecmp}$)$ is the statement
that for every \deecmp poset $Q$ and every family $\D$ of
cardinality $|\D|=\aleph_1$ consisting of dense subsets of $Q$, there
exists a filter $G\subseteq Q$ intersecting every member of $\D$. 
We define $\ma(\Phi$ and {\deecmp}$)$ analogously, where
$\Phi$ is some property of posets (e.g.~properness). 
\end{defn}

\begin{prop}
\label{p-2}
Every poset in \deecmp is completely proper.
\end{prop}

\begin{coro}[Shelah]
\label{Shelah:ma-alpha}
$\ma(\aupalpha$\tu{-proper and {\deecmp}}$)$ is consistent
with $\ch$ relative to the consistency of a supercompact cardinal.
\end{coro}

\begin{rem}
\label{r-2}
Although we have made the effort to distinguish when a formula $\psi$
\emph{provably} has some property, it will not actually have a direct
bearing on the topics in this paper. Provability is only needed for
the property of describing a poset.

In particular, although, as
already mentioned, a completeness system should be in the ground
model, simple completeness systems were designed with the following
property in mind (adapted to the present context).
\end{rem}

\begin{prop}
\label{l-26}
Suppose $(\wp,\tau)$ provably describes a poset.
The statement $\ulc\varphi$ is a completeness system for
$(\wp,\tau)\urc$ is absolute between transitive models \tu(of enough of
$\zfc$\tu) that have the same reals. The $\sigma$-completeness property
is similarly absolute.
\end{prop}
\begin{proof}
The point is that the larger model has 
no new isomorphism types of countable elementary
submodels. See e.g.~\cite{hbst}. 
\end{proof}

\begin{defn}
\label{d-16}
In the context of parameterized games, 
we may refer to a notion of suitability $F$ as \emph{describing} 
some type of family of subsets of $M$. Of course we may do the same
when $F$ is given by a formula $\psi$; moreover, in the latter case we
can say that $\psi$ \emph{provably describes} some family to indicate
that this fact is provable in $\zfc$.
\end{defn}

\begin{example}
\label{x-3}
We say that $\psi$ describes a $P$-filter, if for all $M\in\S$, for all
$\vec a\in\T\cap M$, $\{x\subseteq M:\psi(M,x,\vec a)\}$ is a $P$-filter on $M$.
\end{example}

\begin{prop}
\label{p-15}
$\psimin$ provably describes a $P$-filter \tu(cf. definition~\tu{\ref{d-15}}\tu). 
\end{prop}
\begin{proof}
For $M\in\S$ and $p\in\poset(\F,\H)\cap M$ 
the described family is $\{y\subseteq M:x_p\subseteq y$, 
and $x\subseteqfnt y$
for all $x\in\F\cap M\}$, which clearly forms a $P$-filter. 
\end{proof}

\begin{lem}
\label{l-7}
Let $\psi$ be a notion of suitability such that \tu(provably\tu)  $\psi\to\psimin$.
Then there is a \tu(provable\tu) $\sigma$-complete completeness system
for $\R(\F,\H)$ for all subfamilies $\F$ and $\H$ of\/ $[\theta]^{\leq\aleph_0}$ for some ordinal\/ $\theta$ 
with $\F$ closed under finite reductions for which $\psi$-globally,
Complete has a forward nonlosing strategy in the game $\ggen(\F,\H)$.
\end{lem}

\begin{rem}
\label{r-10}
Lemma~\ref{l-7} is asserting the existence of a completeness system 
for $(\poset,\tau)$ where 
$\tau(\F,\H)$ expresses $\ulc \F$ and $\H$ are families of countable sets of
ordinals with $\F$ closed under finite reductions such that
$\psi$-globally, Complete has a
forward nonlosing strategy for $\ggen(\F,\H)\urc$.
\end{rem}

\begin{proof}
We fix a definable method of coding
\begin{itemize}[leftmargin=*]
\item a subset $y$ of $\theta\cap M$,
\item for each $t\in\Fin(\theta)\cap M$, 
  a function $\Phi(t)$ with domain a subset of $M$ and
  codomain $\Fin(\theta)\cap M$,
\end{itemize}
by subsets $Z\subseteq M$. 
Then we let $\varphi(G,Z;\F,\H,p)$ be a formula expressing $\ulc$if both
\begin{enumerate}[leftmargin=*, label=(\alph*), ref=\alph*, widest=d]
\item\label{item:8} $\phimin(y;\F,\H,p)$ (cf. definition~\ref{d-3}),
\item\label{item:12} $\Phi(t)$ is a forward strategy for Complete in the game
  $\ggen(M,y\cup t,\F,\H,p)$ for all $t\in\Fin(\theta)$ such that 
  $\phimin(y\cup t;\F,\H,p)$,
\save
\end{enumerate}
then $G$ is the filter generated by $(p_k:k<\omega)$, 
where $(p_k:k<\omega)$ is some sequence satisfying:
there exists $m<\omega$ and $t\in\Fin(\theta)$ such that 
\begin{enumerate}[leftmargin=*, label=(\alph*), ref=\alph*, widest=d]
\restore
\item\label{item:20} $\phimin(y\cup t;\F,\H,p)$,
\item\label{item:9} the game $\ggen(M,y\cup t,\F,\H,p)$ is played 
and $(p_k,s_k)$ denotes move $k$,
\item\label{item:11} Complete plays $s_k\supseteq\Phi(t)(P_k)$ 
for all $k\geq m$,
 where $P_k$ is the position after Extender's $k\Th$ move$\urc$,
\end{enumerate}
where $y$ and $\Phi$ are the objects coded by $Z$.\footnote{More 
precisely, we fix a method of coding sequences of elements of
$M$ by subsets of $M$, and then $\varphi$ 
is of the form $\ulc\exists Y\, \rho(G,Z,Y;\F,\H,p)\urc$ 
where $\rho$ has no second order quantifiers and $Y$
is used to code the game played, i.e.~it codes $(p_k,s_k:k<\omega)$.} 

Assuming that $\vec a=(\F,\H)$ satisfies the hypotheses, 
then given $p\in\poset(\F,\H)^M$ 
we need to check that the family of subsets of $M$ given by
\begin{equation}
  \label{eq:12}
  \G_Z=\{G\in\Gen(M,\poset(\F,\H),p):M\models\varphi(G,Z;\F,\H,p)\}
  \espc\text{($Z\subseteq M$)}
\end{equation}
has the required properties. 
First we note that $\G_Z$ is a nonempty subset of $\Gen(M,\poset(\F,\H),p)$ 
whenever $M\models\varphi(G,Z;\F,\H,p)$.
This is because when clauses~\eqref{item:8} and~\eqref{item:12} are
true, letting
$y$ and $\Phi$ be the objects coded by $Z$, $y\supseteq x_p$
bounds $\F\cap M$ by proposition~\ref{p-11},
and $\Phi(\emptyset)$ is a (forward) strategy for Complete in the
game $\ggen(M,y,\F,\H,p)$. Thus setting $m=0$ and $t=\emptyset$, 
Extender can play according to a nonlosing strategy
for $\ggen(M,y,\F,\H,p)$ by corollary~\ref{o-6} 
and Complete can play valid moves $\Phi(\emptyset)(P_k)$, proving that a
sequence $(p_k:k<\omega)$ satisfying~\eqref{item:20}--\eqref{item:11} 
does exist,
and moreover $G$ is generic over $M$ because Extender does not lose.

Next we show that the family forms a $\sigma$-complete filter base, 
which in particular establishes~\ref{item:31}. 
We in fact establish the stronger property that its
upwards closure is a $\sigma$-complete filter. 
This is of course done by diagonalizing the coded objects.
Take $Z_n\subseteq M$ ($n<\omega$) and assume without loss of generality
that~\eqref{item:8} and~\eqref{item:12} are satisfied for all $n$. 
For each $n$, let $y_n$ and $\Phi_n$ be the objects coded by~$Z_n$.
Since $\psimin$ describes $\sigma$-$\supseteqfnt$-directed family,
using proposition~\ref{p-11} there exists $y_\omega\subseteq M$ such that  
$M\models\phimin(y_\omega;\F,\H,p)$ and $y_\omega\subseteqfnt y_n$ for all $n<\omega$.
Fix some enumeration $(u_j:j<\omega)$ of $\Fin(\theta)\cap M$.
For each $t\in\Fin(\theta)\cap M$ such that $M\models\phimin(y_\omega\cup t;\F,\H,p)$, let
\begin{multline}
  \label{eq:13}
  \Psi(t)(P)=\bigcup\bigl\{\Phi_i(u_j)(P):i,j< |P|\text{, }
  M\models\phimin(y_i\cup u_j;\F,\H,p)\text{, }\\
  P\text{ is a position of the game }\ggen(M,y_i\cup u_j,\F,\H,p)\bigr\}
\end{multline}
for every position $P$ in the game $\ggen(M,y_\omega\cup t,\F,\H,p)$ for which it is
Complete's turn to play (its $|P|-1\Th$ move), where an empty union is taken to be
the empty set, and then set
\begin{equation}
  \Phi_\omega(t)(P)=\Psi(t)(P)\cap (y_\omega\cup t).\label{eq:10}
\end{equation}
Then letting $Z_\omega\subseteq M$ code $y_\omega$ and $\Phi_\omega$, 
clearly~\eqref{item:8} and~\eqref{item:12} hold for $Z_\omega$.
We will show that $\G_{Z_\omega}\subseteq\bigcap_{n=0}^\infty\G_{Z_n}$. 
To see this, 
suppose $(p_k:k<\omega)$ is a sequence satisfying~\eqref{item:20}--\eqref{item:11}
for $Z_\omega$, witnessed by $m<\omega$, 
$t\in\Fin(\theta)\cap\nobreak M$ and Complete's moves $(s_k:k<\omega)$.
Given $n$, we need to show that $(p_k:k<\omega)$ 
satisfies~\eqref{item:20}--\eqref{item:11} for $Z_n$. 
Since $\psimin$ describes a filter,  
$M\models\phimin(y_n\cup\nobreak y_\omega\cup t;\F,\H,p)$. 
Then letting $t'\in\Fin(\theta)\cap M$ satisfy
\begin{equation}
y_n\cup t'=y_n\cup y_\omega\cup t,
\label{eq:14}
\end{equation}
$t'$ witnesses that~\eqref{item:20} holds.
Put $m'=\max\{m,n,j\}$ where $u_j=t'$.
Then satisfaction of the conditions for $Z_n$ 
is witnessed by $m'$ and $t'$, 
since~\eqref{eq:14} implies that~\eqref{item:9} holds and
for each $k\geq m$, $\Psi(t)(P_k)\cap(y_n\cup t')$ is a valid move in the game
$\ggen(M,y_n\cup t',\F,\H,p)$, where $P_k$ is the position in the game
$\ggen(M,y_\omega\cup t,\F,\H,p)$ after Extender's $k\Th$ move 
(and thus $|P_k|=k+1$), 
because by~\eqref{eq:10}, $\Psi(t)(P_k)\cap(y_\omega\cup t)\subseteq s_k$.
Now $\Psi(t)(P_k)\supseteq \Phi_n(t')(P_k)$  for all $k\geq\max\{n,j\}$
proving~\eqref{item:11} for $Z_n$ with $s_k:=\Psi(t)(P_k)\cap(y_n\cup t')$ 
for $k=m,m+1,\dots$.

It remains to verify~\ref{item:21}. Indeed the requirement of
equation~\eqref{eq:39} guarantees a $y\subseteq M$ satisfying 
$\psi(M,y,\F,\H,p)$; and by assumption, 
Complete has a forward nonlosing strategy $\Phi(t)$ in the game 
$\ggen(M,y\cup t,\F,\H,p)$ for all $t\in\Fin(\theta)\cap M$ satisfying 
$\psi(M,y\cup t,\F,\H,p)$. Let $Z\subseteq M$ code $y$ and $\Phi$.
Then~\eqref{item:8} and~\eqref{item:12} hold since $\psi\to\psimin$, 
and thus every member $G\in\G_Z$ is generated by
$(p_k:k<\omega)$ resulting from the game $\ggen(M,y\cup t,\F,\H,p)$ 
being played for some~$t$. Since Extender does not lose as
$\<p_k:k<\omega\>=G\in\Gen(M,\poset(\F,\H))$,
and since Complete plays supersets of the forward
nonlosing strategy $\Phi(t)$ for all but finitely many moves by~\eqref{item:11}, 
the game results in a draw by proposition~\ref{p-10}. 
Therefore $G=\<p_k:k<\omega\>\in\Genc(M,\poset(\F,\H),p)$ by proposition~\ref{p-4}.
\end{proof}

\subsection[Upwards boundedly order closed families]{Upwards boundedly order closed subfamilies}
\label{sec:upper-prer-subf}

In~\cite{Hir1} we introduced the following weakening of order closedness.

\begin{defn}
\label{d-4}
We say that a subset $A$ of a poset $P$ is 
\emph{upwards boundedly order closed} 
if for every nonempty $B\subseteq A$ with an upper bound $p\in P$
(i.e.~$b\leq p$ for all $b\in B$): 
if $B$ has a supremum $\spr B$ in $P$, then $\spr B\in A$. 

In the present context of subfamilies $\H\subseteq\power(\theta)$,
we say that $\H$ is \emph{upwards boundedly order closed} to indicate
that it is so in the tree $(\power(\theta),\sqsubseteq)$ of initial
segments. 
\end{defn}

There is a simple criterion for it.

\begin{lem}
\label{l-5}
Every convex \tu(cf.~\tu{\Section\ref{sec:terminology}}\tu) 
subset of a poset is upwards boundedly order closed.
\end{lem}
\begin{proof}
Let $(P,\le)$ be a poset, and let $C\subseteq P$ be convex. 
Take a nonempty $B\subseteq C$, 
say $b\in B$, with an upper bound $p\in P$. Suppose $B$ has a supremum
$\bigvee B$ in $P$. Then $b\le \bigvee B\le p$ implies $\bigvee B\in
C$ by convexity.
\end{proof}

\begin{cor}{0}
\label{p-26}
If\/ $\H$ is a convex subfamily of\/
$(\power(\theta),\sqsubseteq)$ then $\H$ is upwards boundedly order closed.
\end{cor}

Applying the definition in the present context gives:

\begin{prop}
\label{p-17}
$\H$ is upwards boundedly order closed iff every nonempty subfamily
$\K\subseteq\H$ with a\/  $y\in\H$, such that $x\sqsubseteq y$ for all $x\in\K$,
satisfies $\bigcup\K\in\H$. 
\end{prop}

\noindent This endows the poset $\poset(\F,\H)$ with the following crucial property.

\begin{prop}
\label{p-18}
Let $\F,\H$ be families of subsets of\/ $\theta$, 
with $\H$ upwards boundedly order closed. 
If a family $Q\subseteq\poset(\F,\H)$ has a common extension in
$\poset(\F,\H)$, then it has a common extension $q$ such that
\begin{equation}
  \label{eq:6}
  x_q=\bigcup_{p\in Q}x_p.
\end{equation}
\end{prop}

In the present context we consider a generalization.

\begin{defn}
\label{d-19}
Let $\H\subseteq\power(\theta)$ and $R\subseteq\theta$.
We say that $\H$ is \emph{upwards boundedly order closed beyond $R$}
if $\bigcup\K\in\H$ 
whenever $\K\subseteq\H$ is a nonempty subfamily with $y\in\H$ such
that $x\sqsubseteq y$ for all $x\in\K$ and such that
\begin{equation}
  \label{eq:27}
  \sup\Bigl(\bigcup\K\Bigr)\notin R.
\end{equation}
\end{defn}

\begin{prop}
\label{p-8}
Let $\F,\H$ be families of subsets of\/ $\theta$, 
with $\H$ upwards boundedly order closed beyond $R\subseteq\theta$. 
If a family $Q\subseteq\poset(\F,\H)$ has a common extension in
$\poset(\F,\H)$ and
\begin{equation}
  \label{eq:51}
  \sup\Biggl(\bigcup_{p\in Q}x_p\Biggr)\notin R,
\end{equation}
then it has a common extension $q$ such that $x_q=\bigcup_{p\in Q}x_p$.
\end{prop}

This property of the family $\H$ 
allows us to obtain $\aupalpha$-properness for the poset
$\poset(\F,\H)$. It will also
be used in the next section (\Section\ref{sec:isomorphic-models}) to obtain a
strong chain condition for the poset. 

\begin{claim}
\label{c-2}
Let $\F$ and $\H$ be subfamilies of\/ $[\theta]^{\leq\aleph_0}$, 
with $\F$ closed under finite reductions 
and $\H$ upwards boundedly order closed.
Suppose $\F,\H\in M\prec H_\kappa$ is countable,
$p\in\poset(\F,\H)\cap M$ and $y\subseteq M$ is a
$\subseteqfnt$-bound of $\F\cap M$ with $x_p\subseteq y$.
If Extender does not have a winning strategy in the game $\ggen(M,y,\F,\H,p)$, then 
there exists $q\in\gen^+(M,\poset(\F,\H),p)$ such that $x_q\subseteq y$. 
\end{claim}
\begin{proof}
The game $\ggen(M,y,\F,\H,p)$ is played with 
Extender playing according to a nonlosing strategy by
corollary~\ref{o-6}. Since Extender's strategy  is not a winning strategy,
Complete can play so that Extender does not win, and hence the game is
drawn. Then Extender's sequence of moves
$(p_k:k<\omega)$ generates a completely $(M,\poset(\F,\H),p)$-generic
filter $G$, say
with extension $q\in\poset(\F,\H)$, by proposition~\ref{p-4}. 
And by proposition~\ref{p-18},
we may assume that $x_q=\bigcup_{\bar p\in G}x_{\bar p}
=\bigcup_{k<\omega}x_{p_k}\subseteq y$ as needed.
\end{proof}

\begin{claim}
\label{a-1}
Let $\F$ and $\H$ be subfamilies of\/ $[\theta]^{\leq\aleph_0}$, 
with $\F$ closed under finite reductions 
and $\H$ upwards boundedly order closed beyond $R\subseteq\theta$, 
such that $\H$ is $\F$-extendable. 
Suppose $\F,\H\in M\prec H_\kappa$ is countable     
with $\sup(\theta\cap M)\notin R$, 
$p\in\poset(\F,\H)\cap M$ and $y\subseteq M$ is a
$\subseteqfnt$-bound of $\F\cap M$ with $x_p\subseteq y$.
If Extender has no winning strategy in the game $\ggen(M,y,\F,\H,p)$, then 
there exists $q\in\gen^+(M,\poset(\F,\H),p)$ such that $x_q\subseteq y$. 
\end{claim}
\begin{proof}
Set $\delta=\sup(\theta\cap M)$. 
We proceed as in the proof of claim~\ref{c-2}, but now since $\H$ is
$\F$-extendable, we also have by proposition~\ref{p-12}
that for every $\xi<\theta$ in $M$, $\D_\xi\in
M$ (cf.~equation~\eqref{eq:2}) is dense, and
thus $x_{p_k}\in\D_\xi$ for some $k$ since Extender did not lose.
This entails that $\sup\bigl(\bigcup_{k<\omega}x_{p_k}\bigr)
=\delta$, and since $\delta\notin R$ by assumption, we use
proposition~\ref{p-8} to obtain the desired extension $q$ with
$x_q\subseteq y$.
\end{proof}

\begin{defn}
\label{d-5}
We say that a suitability function $F$  is
\emph{coherent} if for all $M\in N$ in $\S$, and all $\vec a\in\T\cap M$,
every $y\in F(N)(\vec a)$ has an $x\subseteq y\cap M$ in $F(M)(\vec a)\cap N$.
When the function is given by a formula $\psi$ we say that $\psi$ is
\emph{provably coherent}, as usual to specify that coherence is
provable in~$\zfc$.
\end{defn}

\begin{prop}
\label{p-21}
If we restrict $\T$ in definition~\tu{\ref{d-15}} 
so that $(\F,\H,p)\in\T$ only if
$(\F,\subseteqfnt)$ is a $\sigma$-directed family of countable sets of
ordinals, then 
$\psimin$ is provably coherent.
\end{prop}
\begin{proof}
Suppose $M\in N$, $\F,\H,p\in M$ and $\psimin(N,y,\F,\H,p)$. 
Since $\F$ is $\sigma$-directed and $M$ is countable, 
by proposition~\ref{p-1},
there exists $x_p\subseteq y'\in\F$ bounding $\F\cap M$, 
and by elementarity we can find such a $y'\in N$.
Now $y'\subseteqfnt y$, and $y'\subseteq N$ as $y'$ is countable, and
thus $y'\cap y\in N$. Hence $y'\cap y\cap M\subseteq y$ is in~$N$ 
and clearly $\psimin(M,y'\cap y\cap M,\allowbreak\F,\H,p)$ holds. 
\end{proof}

\begin{defn}
\label{d-21}
Let $\psi$ be a formula that is to be used in definition~\ref{d-15}.
We say that $\psi$ \emph{respects $\ggen$} if for all $M\in\S$, all
$\F,\H\in M$ and all $p\in\poset(\F,\H)\cap M$, if $\psi(M,y,\F,\H,p)$ and
  $(p_0,s_0),\dots,(p_k,s_k)$ is a position in the game
  $\ggen(M,\allowbreak y,\F,\H,p)$, then 
  $\psi\bigl(M,y\setminus\bigcup_{i=0}^k s_i,\F,\H,p_k\bigr)$.
We may also specify that $\psi$ \emph{provably respects $\ggen$}.
\end{defn}

\begin{prop}
\label{p-33}
For any pair of families $(\F,\H)$, $\psimin$ provably respects $\ggen$.
\end{prop}
\begin{proof}
Immediate from the definitions.
\end{proof}
\subsubsection{$\aalpha$-properness}
\label{sec:aupalpha-properness}

\begin{lem}
\label{l-4}
Let\/ $\F$ be a subfamily of\/ $\powcnt\theta$
closed under finite reductions, 
and let\/ $\H$ be an upwards boundedly order closed subset 
of\/ $([\theta]^{\leq\aleph_0},\sqsubseteq)$.
Suppose that $\psi\to\psimin$ is coherent and respects $\ggen$. 
If Complete has a nonlosing strategy  for\/ $\ggen(\F,\H)$, $\E$-$\psi$-globally,
then $\poset(\F,\H)$ is\/ $\E$-$\aalpha$-proper.
\end{lem}
\begin{proof}
The proof is by induction on $\alpha<\oone$. 
The induction hypothesis is: 
\begin{principle}{$\ih\beta$}
For every tower $M_0\in M_1\in\cdots$ of members of $\E$, that are
also elementary submodels of $H_\kappa$,
of height $\beta+1$ with $\F,\H\in M_0$,
every $p\in\poset(\F,\H)\cap M_0$ and every $y\subseteq M_\beta$ such that
$\psi(M_\beta,y,\F,\H,p)$,
there exists $q\in\gen^+(\{M_\xi:\xi\leq\nobreak\beta\},\poset(\F,\H),p)$
with $x_q\subseteq y$.
\end{principle}
\noindent Assume that $\ih\beta$ holds for all $\beta<\alpha$. 

Suppose $M_0\in M_1\in\cdots$ is a
tower in $\E$ of elementary submodels of $H_\kappa$ of height $\alpha+1$, 
with $\F,\H\in M_0$, and take $p\in\poset(\F,\H)\cap M_0$. Suppose we are given
$y\subseteq M_\alpha$ satisfying $\psi(M_\alpha,y,\F,\H,p)$.
In the case $\alpha=0$, $\ih0$ follows immediately from claim~\ref{c-2} 
because $\psi\to\psimin$ and Complete has a nonlosing strategy in the
game $\ggen(M_\alpha,y,\F,\H,p)$.

Consider the case $\alpha=\beta+1$ is a successor. 
Since $\psi$ is coherent, there exists $y'\subseteq y\cap M_\beta$  in $M_\alpha$
satisfying $\psi(M_\beta,y',\F,\H,p)$.
Now applying $\ih\beta$ in the model $M_\alpha$ with $y:=y'$, 
we obtain $q'\in\gen^+(\{M_\xi:\xi\leq\beta\},\poset(\F,\H),p)\cap M_\alpha$ with
$x_{q'}\subseteq y$. Then we use claim~\ref{c-2} with $M:=M_\alpha$ 
and $p:=q'$ to extend $q'$ to $q$
completely generic over $M_\alpha$ with $x_q\subseteq y$.

Assume now that $\alpha$ is a limit, say
$(\eta_k:k<\omega)$ is a strictly increasing sequence cofinal in~$\alpha$.
The game $\ggen(M_\alpha,y,\F,\H,p)$ is played with Complete using its 
nonlosing strategy.
After $k+1$ moves $(p_0,s_0),\dots,(p_{k},s_{k})$
have been played, assume that $p_{k}\in M_{{\eta_{k-1}}+1}$ 
(taking $\eta_{-1}+1=0$)
is completely generic over $\{M_\xi:\xi\leq\eta_{k-1}\}$ 
with $x_{p_k}\subseteq y$. 
Since $\psi$ respects $\ggen$ and $(p_0,s_0),\dots,(p_k,s_k)$ is a
valid position in the game $\ggen(M_{\eta_k},y,\F,\H,p)$,
$\psi\bigl(M_{\eta_k},y\setminus\bigcup_{i=0}^k s_i,\F,\H,p_k\bigr)$ holds.
And by coherence, there exists $y'\subseteq
y\setminus\bigcup_{i=0}^k s_i\cap M_{\eta_k}$ in
$M_{\eta_k+1}$ satisfying $\psi(M_{\eta_k},y',\F,\H,p_k)$. 
Extender can now make a move  $p_{k+1}\in\nobreak M_{\eta_{k}+1}$ 
completely generic over $\{M_\xi:\xi\leq\nobreak\eta_k\}$ with
$x_{p_{k+1}}\subseteq y$ by applying
$\ih{\eta_k-(\eta_{k-1}+1)}$ in the model $M_{\eta_k+1}$ to
the tower $M_{\eta_{k-1}+1}\in\cdots\in M_{\eta_k}$, $p:=p_k$ 
and $y:=y'\setminus\bigcup_{i=0}^{k}s_k$.

The continuity of the $\in$-chain ensures that Extender does not lose the game.
And since Complete does not lose, $p_0\leq p_1\leq\cdots$ 
has  a common extension $q$ with $x_q\subseteq y$ by proposition~\ref{p-18}. 
Since $q$ is completely generic over 
$\{M_\xi:\xi\leq\alpha\}$ the induction is complete. 

To see that $\poset(\F,\H)$ is $\E$-$\aalpha$-proper, 
suppose $M_0\in \cdots\in M_\beta$ is a tower in $\E$ with $\F,\H\in M_0$, 
and take $p\in \poset(\F,\H)\cap M_0$. 
Equation~\eqref{eq:8} gives a $y\subseteq M_\beta$ 
such that $\psi(M_\beta,y,\F,\H,p)$. Then $\ih\beta$ implies the existence of
$q\geq p$ generic over every member of the tower.
\end{proof}

\begin{lem}
\label{l-27}
Let\/ $\F$ be a subfamily of\/ $\powcnt\theta$
closed under finite reductions, 
and let\/ $\H\subseteq\powcnt\theta$ 
be an upwards boundedly order closed beyond $R\subseteq\theta$ that is
$\F$-extendable. 
Suppose that $\psi\to\psimin$ is coherent and respects $\ggen$. 
If $\E\subseteq\E(\theta\setminus R,\theta)$ and Complete has a nonlosing strategy  
for\/ $\ggen(\F,\H)$, $\E$-$\psi$-globally,
then $\poset(\F,\H)$ is\/ $\E$-$\aalpha$-proper.
\end{lem}
\begin{proof}
The proof is the same as the proof of lemma~\ref{l-4}, except that we
use claim~\ref{a-1} in place of claim~\ref{c-2}. This is justified by
the definition of $\E$, because $\H$ is $\F$-extendable. 
\end{proof}

We shall require the following basic observation on an equivalent of
$\E$-$\aupalpha$-properness.
Cf.~equation~\eqref{eq:16} for the notation.

\begin{lem}
\label{l-28}
Let $S\subseteq\theta$ be stationary for some regular cardinal $\theta$. 
If a forcing notion $P$ is 
$\E(S\setminus A,\theta)$-$\aalpha$-proper for some $A\in\NS_\theta$, 
then $P$ is $\E(S,\theta)$-$\aalpha$-proper. 
\end{lem}
\begin{proof}
This is because for any countable $M\prec H_\kappa$ with $P\in M$, we
may assume that $A\in M$, and thus $\sup(\theta\cap M)\notin A$. 
Hence the set $A$ does not interfere with $\E(S)$-$\aupalpha$-properness.
\end{proof}

\subsection{Isomorphic models}
\label{sec:isomorphic-models}

We introduce a new variant of $\ma$ here in definition~\ref{d-24}, 
that is consistent with $\ch$ relative to $\zfc$.
It is based on Shelah's theorem below
(\cite[Ch.~VIII,~Lemma~2.4]{MR1623206}) for obtaining $\aleph_2$-cc
iterations.  Alternatively, we could have used the appropriate axiom
from~\cite[Ch.~VIII,~\Section3]{MR1623206}.

This will only be possible for the case
$\theta=\oone$, and we still need some theorem in this case because in general for
$\F,\H\subseteq[\oone]^{\leq\aleph_0}$, the cardinality of our poset is large:
\begin{equation}
  \label{eq:18}
  |\poset(\F,\H)|=2^{\aleph_1^{\aleph_0}}\geq\aleph_2
\end{equation}
(and equal to $\aleph_2$ under $\gch$).  

\begin{defn}
\label{d-25}
We say that a suitability function $F$ (as in definition~\ref{d-9}) 
\emph{respects isomorphisms} if
for every two isomorphic  models $M$ and $N$ in $\S$, 
for every (first order) isomorphism $h:M\to N$ fixing $M\cap N$,
\begin{multline}
  \label{eq:48}
  x\in F(M)(\vec a)\Iff h[x]\in F(N)(h(\vec a))\\
  \text{for all $x\subseteq M$ and all $\vec a\in\T\cap M$}.
\end{multline}
If $\psi$ is a formula describing a suitability function, then we say
that $\psi$ (provably) \emph{respects isomorphisms} if (it is provable that) 
the function described by $\psi$ respects isomorphisms. 
\end{defn}

\begin{example}
\label{x-1}
In the context of definition~\ref{d-15}, with $\F$ and $\H$ fixed,
$\psi$ respects isomorphisms iff for all $M,N\in\S$, 
and all isomorphisms $h:M\to N$ fixing $M\cap N$,
\begin{multline}
  \label{eq:24}
  \psi(M,y,\F,\H,p)\Iff h[y]\in \psi(N,h[y],h(\F),h(\H),h(p))\\
  \text{for all $p\in\poset(\F,\H)\cap M$}.
\end{multline}
\end{example}

\begin{prop}
\label{p-36}
Provided that we restrict $\T$ to only include families $\F$ of
countable sets of ordinals, $\psimin$ provably respects isomorphisms.
\end{prop}
\begin{proof}
Let $h:M\to N$ be an isomorphism.
Assume $\psimin(M,y,\F,\H,p)$, i.e. $x_p\subseteq y$ and $y\subseteq M$ is a
$\subseteqfnt$-bound of $\F\cap M$. 
Then $h[y]$ is a $\subseteqfnt$-bound of $h(\F)\cap N$ and
$x_{h(p)}=h(x_p)\subseteq h[y]$, because $x\subseteq M$ for all
$x\in\F\cap M$, yielding $\psimin(N,h[y],h(\F),h(\H),h(p))$.
\end{proof}


\begin{lem}
\label{l-8}
Let $\F$ and $\H$ be subfamilies of\/
$[\theta]^{\leq\aleph_0}$ with $\F$ closed under finite reductions.
Suppose that $\psi\to\psimin$, and $\psi$-globally, 
Complete has a nonlosing strategy for $\ggen(\F,\H)$.
Then for any two isomorphic countable models $M,N\prec
H_\kappa$, say $h:M\to N$ is an isomorphism, with $\F,\H\in M$, if
\begin{multline}
  \label{eq:49}
  \psi(M,y,\F,\H,p)\Iff\psi(N,h[y],h(\F),h(\H),h(p))\\
  \textup{for all $y\subseteq M$ and all $p\in\poset(\F,\H)\cap M$},
\end{multline}
then for every $p\in\poset(\F,\H)\cap M$, there exists
$G\in\Genc(M,\poset(\F,\H),p)$ such that $h[G]\in\Genc(N,\poset(h(\F),h(\H)),h(p))$.
\end{lem}
\begin{proof}
Suppose we are given two isomorphic countable models $M,N\prec H_\kappa$, say
$h:M\to N$ is an isomorphism, with $\F,\H\in M$, and $p\in\poset(\F,\H)\cap M$.\linebreak
By equation~\eqref{eq:39}, there exists $y\subseteq M$ such that
$\psi(M,y,\F,\H,p)$. Thus\linebreak $\psi(N,\allowbreak h[y],\allowbreak
h(\F),\allowbreak h(\H),\allowbreak h(p))$ by 
equation~\eqref{eq:49}.
Also note that $\psi$-globally, Complete has a nonlosing
strategy for $\ggen(h(\F),h(\H))$ by elementarity. 

The games $\ggen(M,y,\F,\H,p)$ and 
$\ggen(N,h[y],h(\F),h(\H),h(p))$ are played simultaneously, and 
we let $\Phi$ and $\Phi'$ denote nonlosing strategies for Complete in the respective games.
Extender plays $p_k$ on its $k\Th$ move in the game\linebreak $\ggen(M,y,\F,\H,p)$, according
to a nonlosing strategy which it has by corollary~\ref{o-6} as $\psi\to\psimin$. 
And
\begin{equation}
  \label{eq:21}
  \text{Extender plays $h(p_k)$ in the game $\ggen(N,h[y],h(\F),h(\H),h(p))$};
\end{equation}
the validity of this move is verified below.
On its $k\Th$ move, Complete plays
\begin{equation}
  \label{eq:17}
  s_k\cup h\inv(t_k)\espc\text{where $s_k=\Phi(P_k)$ and $t_k=\Phi'(h(P_k))$}
\end{equation}
in the game $\ggen(M,y,\F,\H,p)$, 
where $P_k$ is the position after Extender's $k\Th$ move,
and plays $h(s_k)\cup t_k$ in the game $\ggen(N,h[y],h(\F),h(\H),h(p))$. 
Note that $t_k\subseteq h[y]\setminus h(x_{p_k})$ by~\eqref{eq:21}, which implies
$h\inv(t_k)\subseteq y\setminus x_{p_k}$, and thus $s_k\cup h\inv(t_k)$ is a valid move for
Complete in the former game, and similarly $h(s_k)\cup t_k$ is a valid move in the latter game.
Also note that Extender's move $h(x_{p_k})$ is valid in the latter game:
\eqref{item:14}~$h(p_k)\geq h(p_{k-1})$ and
\eqref{item:15}~$h(x_{p_k})\subseteq h[y]\setminus\bigcup_{i=0}^{k-1} h(s_i)\cup t_i$ as
$x_{p_k}\subseteq y\setminus\bigcup_{i=0}^{k-1}s_i\cup h\inv(t_i)$.

After the games, let $G$ be the filter of $\poset(\F,\H)\cap M$ generated by
$(p_k:k<\nobreak\omega)$, and let $H$ be the filter of $\poset(h(\F),(\H))\cap N$ generated by
$(h(p_k):k<\nobreak\omega)$, so that $H=h[G]$. 
By proposition~\ref{p-10}, Complete does not lose either games. 
Extender does not lose the former game since it played according to a nonlosing
strategy, and it does not lose the latter game, because for every dense
$D\subseteq\poset(h(\F),h(\H))$ in $N$, $p_k\in h\inv(D)$ for some~$k$, and thus
$h(p_k)\in D$. Therefore,\linebreak both
$G\in\Genc(M,\poset(\F,\H),p)$ and $h[G]\in\Genc(N,\poset(\F,\H),h(p))$ by
proposition~\ref{p-4}. 
\end{proof}

Recall that a poset $P$ satisfies the \emph{properness isomorphism condition} if for
every two isomorphic countable $M,N\prec H_\kappa$, 
for $\kappa$ a sufficiently large regular cardinal, 
via $h:M\to N$ with $P\in M\cap N$ and $h$ fixing $M\cap N$, 
for every $p\in P\cap M$, 
there exists $G\in\Genc(M,P,p)$ such that $h[G]\in\Genc(N,P,h(p))$
and moreover there exists $q\in P$ extending both $G$ and $h[G]$ (this
is the $\aleph_2$-pic from~\cite[Ch.~VIII,~\Section2]{MR1623206}).

\begin{thmo}[Shelah]
\label{shelah:isomorphism}
Assume $\ch$. Let $\vec P=(P_\xi,\dot Q_\xi:\xi<\mu)$ be a countable support
iterated forcing construction of length $\mu\le\omega_2$. 
If each iterand satisfies the properness isomorphism condition then
$P_\mu$ satisfies the $\aleph_2$-chain condition.
\end{thmo}

\begin{defn}
\label{d-24}
We write $\ma(\text{\alphaproper}+\text{\deecmp}+\text{pic}
+\text{$\Delta_0$-$H_{\aleph_2}$-definable})$ 
to be interpreted as in
definition~\ref{d-17}, where \emph{pic} denotes the class of posets
satisfying the properness isomorphism condition, 
and \emph{$\Delta_0$-$H_{\aleph_2}$-definable} denotes those posets $P$
that are $\Delta_0$ definable over $H_{\aleph_2}$, i.e.~there exists a
$\Delta_0$ formula (in the L\'evy hierarchy) 
$\varphi(v_0,\dots,v_n)$ and parameters $a_1,\dots,a_n\in
H_{\aleph_2}$ 
such that $P=\bigl\{x\in
H_{\aleph_2}:H_{\aleph_2}\models\varphi[x,a_1,\dots,a_n]\bigr\}$. 
\end{defn}

\begin{example}
\label{x-6}
Our  class $\poset$ of posets is clearly $\Delta_0$. 
Thus $\poset(\F,\H)$ is $\Delta_0$-$H_{\aleph_2}$-definable whenever
$\F,\H$ are families of subsets of $\theta<\omega_2$. 
\end{example}

\begin{thm}
\label{u-13}
$\ma(\tu{\alphaproper}+\tu{\deecmp}+\tu{pic}
+\tu{$\Delta_0$-$H_{\aleph_2}$-definable})$ is consistent with $\ch$ relative
to $\zfc$. 
\end{thm}
\begin{proof}
Beginning with ground model satisfying $\gch$,
we construct an iteration $(P_\xi,\dot Q_\xi:\xi<\omega_2)$ of forcing
notions of length $\omega_2$. At every stage $\xi$ of the iteration it
is arranged that
\begin{enumerate}[leftmargin=*, label=(\roman*), ref=\roman*, widest=iii]
\item\label{item:58} $P_\xi\forces\gch$,
\item $P_\xi$ has a dense suborder of cardinality at most $\aleph_2$,
\item\label{item:59} $P_\xi$ has the $\aleph_2$-cc.
\end{enumerate}
This entails that there are $\aleph_2$ many $P_\xi$-names for members
of $H_{\aleph_2}$, where $H_{\aleph_2}$ is being interpreted in the
forcing extension. Using standard bookkeeping methods, we can arrange
a sequence $(\varphi_\xi,\dot a_\xi:\xi<\omega_2)$ where $\varphi_\xi$
is a $\Delta_0$ formula and $\dot a_\xi$ is a $P_\xi$-name for a parameter in
$H_{\aleph_2}$, so that regarding $P_\xi$-names as also being
$P_\eta$-names for $\xi\le\eta$, every pair $(\varphi,\dot a)$ where
$\dot a$ is a $P_\xi$-name appears as $(\varphi_\eta,\dot a_\eta)$ for
some $\eta\ge\xi$. By skipping steps if necessary, we may assume that
$P_\xi$ forces that $\dot Q_\xi$ is the object defined by $(\varphi_\xi,\dot a_\xi)$
over $H_{\aleph_2}$, and that $\dot Q_\xi$ is an \alphaproper poset in
the class \deecmp (cf.~definition~\ref{d-17}), with the pair of
formulae describing $\dot Q_\xi$ fixed in $V$,\footnote{This is a
  subtle point. To apply Shelah's Theorem on
  page~\pageref{Shelah:alpha}, we need to to know the pair of formulae
  describing $\dot Q_\xi$ in the ground model, i.e.~a $P_\xi$-name for
  a pair of formulae does not suffice. In practice, however, this does
  not pose any difficulties.} and is in pic.

Since every $P_\xi\forces|\dot Q_\xi|\le\aleph_2$ and $\dot Q_\xi$ is in
pic, it follows from Shelah's Theorem above
that~\eqref{item:58}--\eqref{item:59} hold at every intermediate stage
of the iteration, and moreover that 
$P_{\omega_2}=\varinjlim{}_{\xi<\omega_2}P_\xi$ has the
$\aleph_2$-cc. By Shelah's Theorem on page~\pageref{Shelah:alpha},
$P_{\omega_2}$ does not add any reals, and thus $\ch$ holds in its
forcing extension. 

In the forcing extension by $P_{\omega_2}$:
Suppose that $\varphi$ is a $\Delta_0$ formula and $a\in H_{\aleph_2}$ so that
$(\varphi,a)$ defines an \alphaproper poset $Q$ in \deecmp and pic
over $H_{\aleph_2}$, and suppose that $\D\subseteq Q$ is a family of
dense subsets of cardinality $|\D|=\aleph_1$. 
By the $\aleph_2$-cc, there exists $\xi$ such that $\D$ is in the
intermediate extension and 
$P_{\omega_2}\forces(\varphi_\xi,\dot a_\xi)=(\varphi,\dot a)$. 
In the intermediate extension by $P_\xi$: The interpretation $Q_\xi$
of the iterand $\dot Q_\xi$ is given by 
$Q_\xi=\{x\in H_{\aleph_2}:H_{\aleph_2}\models\varphi[x,a]\}$. 
Since $\Delta_0$ formulae are absolute between transitive models, it
follows that $Q_\xi$ is a suborder of $Q$. Thus $D\cap Q_\xi$ is dense
for all $D\in\D$. And at the $\xi\Th$ stage we force a filter
intersecting $D\cap Q_\xi\subseteq D$ for all $D\in\D$. Notice that
$\aupalpha$-properness is downwards absolute, and by
proposition~\ref{l-26}, $Q_\xi$ is also in \deecmp in the intermediate
model.
\end{proof}

\begin{rem}
\label{r-6}
The same proof shows that we can allow posets whose base set is
$\Sigma_1$-$H_{\aleph_2}$-definable but the ordering still must be
$\Delta_0$-$H_{\aleph_2}$-definable. 
\end{rem}

It sometimes possible to obtain $\poset(\F,H)$ with the properness isomorphism
condition for $\F,\H\subseteq[\oone]^{\leq\aleph_0}$ for the following reason.

\begin{prop}
\label{p-16}
For $M,N\prec H_\kappa$ and $h:M\to N$ an isomorphism, 
$h(\alpha)=\alpha$ for all $\alpha<\oone$. 
Thus $[\oone]^{\leq\aleph_0}\cap M=h\bigl[[\oone]^{\leq\aleph_0}\cap M\bigr]
=[\oone]^{\leq\aleph_0}\cap N$. 
\end{prop}
\begin{proof}
This immediately follows from the fact that $\alpha\subseteq M$ for every countable
ordinal $\alpha\in M$. 
\end{proof}

We will see that for some families $\F,\H\subseteq[\oone]^{\leq\aleph_0}$,
lemma~\ref{l-8} is already enough to establish the properness isomorphism condition.

\begin{cor}[5]{1}
\label{o-7}
Let $\F,\H\subseteq[\oone]^{\leq\aleph_0}$ be subfamilies with $\F$ closed under
finite reductions and $\H$ upwards boundedly order closed.
If $\psi\to\psimin$ respects isomorphisms for the fixed pair $(\F,\H)$
and $\psi$-globally, Complete has a nonlosing strategy for $\ggen(\F,\H)$, 
then $\poset(\F,\H)$ satisfies the properness isomorphism condition.
\end{cor}
\begin{proof}
Suppose $M,N\prec H_\kappa$ are countable, $h:M\to N$ is an isomorphism fixing
$M\cap N$ and $\F,\H\in M\cap N$. Since $h$ fixes $\F$ and $\H$, the
fact that $\psi$ respects isomorphisms entails equation~\eqref{eq:49}.
Take $p\in\poset(\F,\H)\cap M$. 
Then by lemma~\ref{l-8}, there exists $G\in\Genc(M,\poset(\F,\H),p)$ such that
$h[G]\in\Genc(N,\allowbreak \poset(\F,\allowbreak \H),\allowbreak h(p))$. 
By the assumption
on $\H$ and proposition~\ref{p-18}, $G$ has a common extension $q$ such that
$x_q=\bigcup_{\bar p\in G}x_{\bar p}$. Using proposition~\ref{p-18} again, $h[G]$ has a common
extension $q'$ such that
\begin{equation}
  \label{eq:15}
  x_{q'}=\bigcup_{\bar p\in h[G]}x_{\bar p}=\bigcup_{\bar p\in G}h(x_{\bar
    p})=\bigcup_{\bar p\in G}x_{\bar p}=x_q
\end{equation}
by proposition~\ref{p-16}.
Therefore $(x_q,\X_q\cup\X_{q'})$ is a condition in $\poset(\F,\H)$ extending both
of the filters $G$ and $h[G]$.
\end{proof}

\subsection{Preservation of nonspecialness}
\label{sec:pres-nonsp}

The following definition is from~\cite[Ch.~IX,~\Section4]{MR1623206}. 
It was developed by Shelah for his proof that $\sh$ does not imply
that all Aronszajn trees are special. 

\begin{defn}
\label{d-7}
Let $T$ be a tree of height $\oone$. 
A poset $(P,\leq)$ is called~\emph{$(T,R)$-preserving} 
  if for every countable $M\prec H_\kappa$, where
  $\kappa$ is some sufficiently large regular cardinal, with $T,R,P\in
  M$ and
 $\delta_M\notin R$, every $p\in P\cap M$ has an
  $(M,P)$-generic extension $q$ that is \emph{$(M,P,T)$-preserving}, 
  i.e.~the following holds for all $x\in T_\delta$: 
if for all $A\subseteq T$ in $M$,
\begin{equation}
\label{eq:one}
x\in A\impls \exists y\in A\spc y\lln T x,
\end{equation}
then for every $P$-name $\dot A\in M$ for a subset of $T$,
\begin{equation}
\label{eq:two}
q\forces\ulc x\in\dot A\impls \exists y\in\dot A\spc y\lln T x\urc.
\end{equation}
In the case $R=\emptyset$, we just say that the poset is \emph{$T$-preserving}; and
when the poset is $(T,R)$-preserving for every $\oone$-tree $T$, we say that the poset
is \emph{$(\oone$-tree$,R)$-preserving}. 
\end{defn}

The following lemma is straightforward.

\begin{lem}
\label{p-22}
If\/ $T$ is a Souslin tree, $R$ is costationary and\/ $P$ is\/
$(T,R)$-preserving, then\/ $T$ remains nonspecial in the forcing
extension by\/ $P$.
\end{lem}

In~\cite{MR1264960} it\label{page:Sch94} is moreover shown that if $T$ is a Souslin
tree and $R\subseteq\oone$ is costationary, 
then $T$ remains nonspecial in the forcing extension by
a countable support iteration of $(T,R)$-preserving posets. 
In~\cite{NS}, the property itself is preserved:

\begin{lem}[Abraham]
\label{l-30}
Let $R\subseteq\oone$ be costationary and $T$ an $\oone$-tree. 
Suppose $\vec P=(P_\xi,\dot Q_\xi:\xi<\mu)$ is a countable support
iteration of length
$\mu<\omega_2$ such that each iterand is $(T,R)$-preserving. 
Then $P_\mu$ is $(T,R)$-preserving.
\end{lem}

\noindent The preservation is proved for iterations of arbitrary
length, but for a different type of iteration
in~\cite[Ch.~IX]{MR1623206}.

It is here that we require the augmented game. 

\begin{lem}
\label{l-11}
Let\/ $\F$ and $\H$ be subfamilies of\/ $[\theta]^{\leq\aleph_0}$
with $\F$ closed under finite reductions. Let $R\subseteq\oone$.
Suppose that $\psi\to\psimin$, and
 $\E(\oone\setminus R)$-$\psi$-globally,  Extender does not have a
 winning strategy for $\ggen^*(\F,\H)$.
Then the forcing notion $\poset(\F,\H)$ is \tu($\oone$-tree, $R$\tu)-preserving.
\end{lem}
\begin{proof}
Let $M\prec H_\kappa$ be a countable elementary submodel with $T,\F,\H\in M$, 
where $T$ is some $\oone$-tree, and  $\delta_M\notin R$.
Let $Z$ be the set of all $t\in T_{\delta_M}$ satisfying definition~\ref{d-7}\eqref{eq:one},
i.e.~for all $A\subseteq T$ in $M$, $t\in A$ implies $u\lln T t$ for some $u\in A$. 
Then let $(t_k,\dot A_k)$ ($k\in\N$) enumerate all pairs
$(t,\dot A)$ where $t\in Z$ and $\dot A\in M$
is an $\poset(\F,\H)$-name for a subset of $T$.
Given a condition $p\in \poset(\F,\H)\cap M$, we need to produce an
$(M,\poset(\F,\H))$-generic extension that is moreover $(M,\poset(\F,\H),T)$-preserving. 
We enumerate all of the dense subsets of $\poset(\F,\H)$ in $M$ 
as $(D_k:k<\omega)$.

For each $k$ and $\bar p\in\poset(\F,\H)$, a set $A^{\bar p}_k\subseteq T$ is defined 
as the collection of all~$t\in T$, such that every cofinal $X\subseteq\F$,
with $x_{\bar p}\subseteq z$ for all $z\in X$, has an extension $q\geq \bar p$ with
\begin{enumeq}
\item\label{item:24} $x_q\subseteq z$ for some $z\in X$,
\item\label{item:25} $q\forces t\in\dot A_k$.
\end{enumeq}
Note that for each $k$, $\bar p\mapsto A^{\bar p}_k$ is definable in $M$.

$\psi$ gives a $y\subseteq M$ such that $\psi(M,y,\F,\H,p)$ (i.e.~by equation~\eqref{eq:8}). 
Now the game $\ggen^*(M,y,\F,\H,p)$ is played. 
Let $(p_k,X_k)$ denote Extender's $k\Th$ move. We shall describe a strategy for Extender.
On even moves $2k$, Extender plays $(p_{2k},X_{2k})$ such that
\begin{equation}
  \label{eq:42}
  p_{2k}\in D_k,
\end{equation}
by invoking corollary~\ref{o-2} using that fact $\psi\to\psimin$. We
can implicitly use some well ordering so as to obtain a strategy.

After the $2k\Th$ move has been played, 
we consider whether or not $t_k\in\nobreak A^{p_{2k}}_k$. 
First suppose that it is. 
Then there exists $u_k\lln T t_k$ in $A^{p_{2k}}_k$. 
Therefore, since $\psi\to\psimin$, we can apply corollary~\ref{l-9}
with $Q=\{q:q\forces u_k\in\dot A_k\}$, so that
equations~\eqref{item:24} and~\eqref{item:25} 
negate the second alternative~\eqref{item:17},
and thus Extender has a move $(p_{2k+1},X_{2k+1})$ such that
\begin{equation}
  \label{eq:20}
  p_{2k+1}\forces u_k\in \dot A_k.
\end{equation}
Otherwise when $t_k\notin A_k^{p_{2k}}$,
there exists a witness $X_{2k+1}$ with $x_{p_{2k}}\subseteq z$ for all $z\in X_{2k+1}$ 
to the fact that there is no $q\geq p_{2k}$
satisfying~\eqref{item:24} and~\eqref{item:25} with $t:=t_k$. 
Extender then makes the valid move 
$(p_{2k+1},X_{2k+1})$ where $p_{2k+1}=p_{2k}$. Again we can use a well
ordering to obtain an actual strategy. 

Since Extender does not have a winning strategy in this game, there
exists a sequence of plays by  Complete such that Complete does not lose.
And neither does Extender lose, because the described strategy is
nonlosing by~\eqref{eq:42}. Thus the game is drawn, 
and we can find a common extension $q\geq p_k$ for all $k$, 
satisfying $\{X_0,X_1,\dots\}\subseteq\X_q$. Since $q$ is generic by proposition~\ref{p-4}, 
it remains to verify that $q$ is $(M,\poset(\F,\H),T)$-preserving. 

\begin{sublem}
For all\/ $k$, if\/ $t_k\notin A_k^{p_{2k}}$ then\/ $q\forces t_k\notin\dot A_k$.
\end{sublem}
\begin{proof}
Fix $k$ and suppose $t_k\notin A_k^{p_{2k}}$. Given $\bar q\geq q$ we need to
prove that it does not force $t_k\in\dot A_k$. But since $X_{2k+1}\in \X_q$, 
there are cofinally many $z\in X_{2k+1}$ with $x_{\bar q}\subseteq z$. And by the
choice of $X_{2k+1}$ we cannot also have $\bar q\forces t_k\in\dot A_k$. 
\end{proof}

Take $t\in T_{\delta_M}$ satisfying equation~\eqref{eq:one}, and an
$\poset(\F,\H)$-name $\dot A\in M$ for a subset of $T$.
Then $(t,\dot A)=(t_k,\dot A_k)$ for some $k$.  Supposing $\bar q\geq q$ forces
$t_k\in\dot A_k$, then $t_k\in A_k^{p_{2k}}$ by the sublemma, and thus 
by~\eqref{eq:20} the verification is complete.
\end{proof}

Let us also observe the following related technical fact.

\begin{lem}
\label{l-31}
$\sigma$-closed forcing notions are $\oone$-tree-preserving.
\end{lem}
\begin{proof}
This can be proved in a similar, but much simpler, manner to
the proof of lemma~\ref{l-11}. 
\end{proof}

We shall also require the following observation that is completely 
analogous in both its statement and its justification to
lemma~\ref{l-28}.

\begin{lem}
\label{l-29}
Let $R\subseteq\oone$ be costationary. 
If a forcing notion $P$ is $(T,R\cup A)$-preserving for some $A\in\NS$,
then $P$ is $(T,R)$-preserving.
\end{lem}

We shall also need to preserve Aronszajn trees, i.e.~make sure no
uncountable branches are added. We use the following notion
from~\cite[Ch.~IX,~\Section4]{MR1623206}.

\begin{defn}
Let $h:\lim(\oone)\to\oone$ be a function whose domain is the
countable limit ordinals.
A tree $T$ of height $\oone$ is called \emph{$h$-st-special} if there
exists a function $f:\bigcup_{\alpha<\oone}T_{h(\alpha)}\to\oone$
satisfying
\begin{enumerate}[leftmargin=*, label=(\roman*), widest=ii]
\item $x\in T_{h(\alpha)}$ implies $f(x)<\alpha$,
\item for all $\alpha<\beta$, $x\in T_{h(\alpha)}$, $y\in
  T_{h(\beta)}$ and $x\len T y$ imply $f(x)\ne f(y)$. 
\end{enumerate}
\end{defn}

It is easy to see that for any $h$, an $h$-st-special tree is not
Souslin. The consequence we are interested in here is the following
proposition, because $h$-st-specialness is obviously upwards absolute for
$\aleph_1$-preserving extensions.

\begin{prop}
\label{p-38}
If $T$ is an $h$-st-special tree, then $T$ has no uncountable branches.
\end{prop}

\begin{lem}
\label{p-39}
Let\/ $T$ be an Aronszajn tree. 
There is an\/ $\oone$-tree-preserving forcing notion of cardinality\/
$\aleph_1$ forcing that\/ $T$ is $h$-st-special for some $h$. 
\end{lem}
\begin{proof}
We refer to Shelah's book.
A slight modification---we prefer ``$\alpha$'' to
``$\alpha\times\nobreak\omega$''---of the forcing notion $Q(T)$ of 
definition~4.2 is shown to force $T$ is
$h$-st-special for some $h$ in claim~4.4, and is shown to be
$\oone$-tree-preserving in lemma~4.6.
\end{proof}

\subsection[Convex families]{Convex subfamilies of $(\power(\theta),\subseteq)$}
\label{sec:convex-h}

The simplest class of families $\H$ to be used as the second parameter  
are those that form convex subsets of
$(\power(\theta),\subseteq\nobreak)$, i.e.~if $x\subseteq y$ are both
in $\H$ and $x\subseteq z\subseteq y$ then $z\in\H$. 

Since $\subseteq$ is a weaker relation that $\sqsubseteq$, such
families $\H$ are automatically convex in the initial segment
ordering, and hence by corollary~\ref{p-26}:

\begin{prop}
\label{p-23}
If\/ $\H$ is a convex subfamily of\/ $(\power(\theta),\subseteq)$ then
$\H$ is upwards boundedly order closed.
\end{prop}

Convexity gives a relatively simple winning strategy for Complete in the game
$\gcomp(y,\F,\H,p)$ when $y\in\H$. 

\begin{lem}
\label{l-10}
Let $\F$ be a directed subfamily of\/ $(\power(\theta),\subseteqfnt)$ with
$\J(\F,\subseteqfnt)$ a $\lambda$-ideal, and\/ $\H$ be a convex
subfamily of\/ $(\power(\theta),\subseteq)$.
Then Complete has a winning strategy for the game
$\gcomp(y,\F,\H,p)$ whenever $y\in\H$  has cardinality $|y|<\nobreak\lambda$ and $x_p\subseteq y$. 
\end{lem}
\begin{proof}
We know that $\bigcup_{k<\omega}\X_{p_k}$ will be countable, and thus we can arrange a
diagonalization $(Y_k:k<\omega)$ in advance, and since the $\X_{p_k}$'s 
will be increasing with $k$, we can also insist that $Y_k\in\X_{p_k}$ for all $k$. 
After Extender plays $p_k$ on move $k$, we take care of some $Y_k\in\X_{p_k}$
according to the diagonalization. 
\begin{equation}
  \label{eq:9}
  Y'_k=\{x\in Y_k:x_{p_k}\subseteq x\}\textup{ is $\subseteqfnt$-cofinal in $\F$},
\end{equation}
and thus
$Y''_k=\{x\in Y'_k:y\subseteqfnt x\}$ is $\subseteqfnt$-cofinal since
we can assume that 
$y\in\downcl\F$ by remark~\ref{r-5}, and $\F$ is $\subseteqfnt$-directed. Now
\begin{equation}
  \label{eq:25}
  Y''_k=\bigcup_{s\in\Fin(y\setminus x_{p_k})}\{x\in Y'_k:y\setminus s\subseteq x\},
\end{equation}
and hence as $|y|<\lambda$, by the assumption on $\J(\F)$
there exists $s_k\in\Fin(y\setminus x_{p_k})$ such that
\begin{equation}
  \label{eq:26}
  Y'''_k=\{x\in Y'_k:y\setminus s_k\subseteq x\}\text{ is $\subseteqfnt$-cofinal}.
\end{equation}
Complete plays $s_k$ on move $k$. This describes the strategy
for Complete. 

And the end of the game, put $x_q=\bigcup_{k<\omega} x_{p_k}$ and
$\X_q=\bigcup_{k<\omega}\X_{p_k}$.
Then $x_q\in\H$ since $x_p,y\in\H$ and $\H$ is convex; for every $k$, 
$x_q\subseteq x$ for all $x\in Y'''_k$ by proposition~\ref{p-6};
and every $Y\in\bigcup_{k<\omega}\X_{p_k}$ appears as $Y_k$ for some $k$, 
and thus $\{x\in Y:x_{q}\subseteq x\}$ is $\subseteqfnt$-cofinal by~\eqref{eq:26}. 
This proves that $q=(x_q,\X_q)\in\poset(\F,\H)$, 
and thus $q\geq p_k$ for all $k$ and Complete wins the game.
\end{proof}

\begin{cor}{0}
\label{c-4}
Let $\F$ be a directed subfamily of\/ $(\power(\theta),\subseteqfnt)$ with
$\J(\F,\subseteqfnt)$ a $\lambda$-ideal, and\/ $\H$ be a convex
subfamily of\/ $(\power(\theta),\subseteq)$.
Then Complete has a forward winning strategy for the game
$\gcomp(y,\F,\H,p)$ whenever $y\in\H$ 
 has cardinality $|y|<\nobreak\lambda$ and $x_p\subseteq y$. 
\end{cor}
\begin{proof}
We apply lemma~\ref{l-16} with $A=\mathrm{Extender}$ and 
$X=\mathrm{Complete}$
to obtain a forward winning strategy, and
thus we need that Complete has a winning strategy, $\E$-$F$-globally,
for some pair $(\E,F)$. Let $\S=\{V\}$
(cf.~section~\ref{sec:terminology}), $\T=\poset(\F,\H)$ and let $F(V)$
be the function with domain $\T$ where
$F(V)(q)=\{x\in\H:|x|<\lambda$ and $x_q\subseteq x\}$.
Then for every $q\in\poset(\F,\H)$ and every $x\in F(V)(q)$, 
Complete has a winning strategy in the game $\gcomp(x,\F,\H,q)$ by
lemma~\ref{l-10}.
Therefore, $\S$-$F$-globally, we have a winning strategy for
Complete in the parameterized game $\game(V,x,q)$, where
$\game(V,x,q)=\gcomp(x,\F,\H,q)$. 

Now suppose $P=(p_0,s_0),\dots,(p_k,s_k)$ is a position in the game
$\gcomp(x,\allowbreak\F,\allowbreak\H,\allowbreak q)$ with Extender to
play. 
Then letting $x'=x\setminus\bigcup_{i=0}^k s_i$ and
$q'=p_k$, clearly $x'\in F(V)(q')$ and
$\gcomp(x',\F,\H,q')=\gcomp(x,\F,\H,q)\restriction P$ is immediate
from the rules of the game. Thus the
hypothesis of lemma~\ref{l-16} is satisfied, and hence Complete has a
forward winning strategy in the game $\game(V,x,q)$, $\S$-$F$-globally. 
In particular, Complete has
a forward winning strategy in the game $\gcomp(y,\F,\H,p)$, where $y$
and $p$ are from the hypothesis of the corollary.
\end{proof}

\begin{cor}{0}
\label{c-3}
Let $\F$ be a\/ $\sigma$-directed subfamily of\/ $(\powcnt{\theta},\subseteqfnt)$ and let $\H$ be convex. 
Suppose that $\psi$ describes a suitability function 
such that   $\psi(M,y,\F,\H,p)$ implies $y\in\H$ and $x_p\subseteq y$.
Then $\psi$-globally, Complete has a forward winning strategy 
in the game $\gcomp(\F,\H)$.
\end{cor}
\begin{proof}
By corollary~\ref{c-4} with $\lambda=\aleph_1$. See also lemma~\ref{l-1}.
\end{proof}

\begin{rem}
\label{r-9}
Thus, assuming the hypothesis of corollary~\ref{c-3}, $\psi$-globally,
Complete has a nonlosing strategy in the game $\ggen(\F,\H)$, by
propositions~\ref{p-25} and~\ref{p-13}. 
\end{rem}

We shall want to use corollary~\ref{c-3} with $\psi$ satisfying
$\psi\to\psimin$ in order to apply the theory of this section. 

\begin{defn}
\label{d-22}
We let $\psicvx(M,y,\F,\H,p)$ be the conjunction of
$\psimin(M,\allowbreak y,\allowbreak\F,\allowbreak\H,\allowbreak p)$  and $y\in\H$. 
\end{defn}

\begin{rem}
\label{p-24}
When both
\begin{enumerate}[leftmargin=*]
\item\label{item:54} $\F$ is a $\sigma$-directed subfamily of
  $(\powcnt\theta,\subseteqfnt)$,
\item\label{item:55} $\H$ is $\subseteqfnt$-cofinal in $\F$, 
\end{enumerate}
then  $\psicvx$ defines a suitability function for $(\F,\H)$ fixed, as
defined in definition~\ref{d-15}, in that
equation~\eqref{eq:39} holds. This is so because for any countable~$M$, 
there exists $x\subseteq M$ in $\F$ bounding $\F\cap M$ as $\F$ is
$\sigma$-directed and members of $\F$ are countable, 
and hence there exists $y\in\H$ bounding $\F\cap M$.
\end{rem}

\begin{cor}[3]{1}
\label{c-9}
Let $\F$ be a\/ $\sigma$-directed subfamily of\/ $(\powcnt{\theta},\subseteqfnt)$
and let $\H$ be convex. 
Then $\psicvx$-globally, Complete has a forward winning strategy 
in the game $\gcomp(\F,\H)$.
\end{cor}
\begin{proof}
Corollary~\ref{c-3} and remark~\ref{p-24}. 
\end{proof}

Note that since we are using $\psicvx$ in the context of remark~\ref{p-24}, 
we are not losing any generality here, in the crucial extendability
property (definition~\ref{d-8}), by requiring~$\F=\H$.

\begin{lem}
\label{l-13}
Let $\H$ be a directed subset of $(\power(\theta),\subseteqfnt)$, with
$\J(\H,\subseteqfnt)$ a $\sigma$-ideal \tu(e.g.~when $\H$ is $\sigma$-directed\tu).
Suppose $\theta$ has no countable decomposition into sets orthogonal to $\H$, 
and is the least ordinal for which this holds. 
Then $\H$ is extendable.
\end{lem}
\begin{proof}
Take $x\in\H$ and let $\X$ be a countable family of cofinal subsets of
$\H$ with $x\subseteq y$ for all $y\in X$ for all $X\in\X$. 
Applying lemma~\ref{l-2} with $\lambda=\aleph_1$,
there exists $\alpha\geq\xi$ satisfying equation~\eqref{eq:3}. Picking
any $X\in\X$, find $y\in X$ such that $\alpha\in y$. 
Then $x\subseteq
x\cup\{\alpha\}\subseteq y$ implies $x\cup\{\alpha\}\in\H$
by convexity, and thus $x\cup\{\alpha\}$ is a witness to extendability.
\end{proof}

In the following three propositions
(propositions~\ref{p-32}--\ref{p-35}) we are restricting $\T$ only to
include $\F$ that are $\sigma$-directed subfamilies of
$(\powcnt\theta,\subseteqfnt)$. 

\begin{prop}
\label{p-32}
$\psicvx$ is provably coherent \tu(cf.~definition~\tu{\ref{d-5})}.
\end{prop}
\begin{proof}
The exact same proof as for proposition~\ref{p-21} works, because the
element $y\cap y'\cap M$ defined there is in $\H$ by convexity.
\end{proof}

\begin{prop}
\label{p-34}
Provided we restrict ourselves to the class of convex\/~$\H$,
$\psicvx$ provably respects $\ggen$ \tu(cf.~definition~\tu{\ref{d-21})}.
\end{prop}
\begin{proof}
Immediate from the definition and convexity.
\end{proof}

\begin{prop}
\label{p-35}
Restricting to $\theta=\oone$, $\psicvx$ provably
respects isomorphisms \tu(cf.~definition~\tu{\ref{d-25})}.
\end{prop}
\begin{proof}
This follows from equation~\eqref{eq:24}, by first observing that
$h[y]=y$ and $h(x)=x$ for all $x\in\F\cap M$ by proposition~\ref{p-16}.
\end{proof}

\subsection{Closed sets of ordinals}
\label{sec:h-consisting-closed}

We shall only consider one more class of families $\H$ in this
paper, but there are certainly others of interest (see e.g.~\cite{Hir}). 

\begin{notn}
For a family $\F$ of subsets of some ordinal $\theta$, 
and a subset $S\subseteq\theta$, typically stationary, 
we let $\club(\F,S)$ denote the family of all closed subsets of~$S$ (in the ordinal
topology) in $\F$, i.e.~$x\subseteq S$ is closed iff $\delta$ is a limit point of
$x$ and $\delta\in S$ imply $\delta\in x$. 
\end{notn}

\begin{prop}
\label{p-14}
$\club(\downcl{\F},S)$ is upwards boundedly order closed
beyond $S$.
\end{prop}
\begin{proof}
Suppose $\K\subseteq\club(\downcl{\F},S)$ with
$y\in\club(\downcl{\F},S)$ so that $x\sqsubseteq y$ for
all~$x\in\K$ and such that $\sup\bigl(\bigcup\K\bigr)\notin S$. 
Then $\bigcup\K\sqsubseteq y$, and in particular is in $\downcl\F$ since $y$ is, 
and $\bigcup\K$ is relatively closed in $S$ because it is the union of
a $\sqsubseteq$-chain of closed sets and by assumption, 
all of its limit points in $S$ are strictly below its
supremum.
\end{proof}

\begin{rem}
\label{r-14}
In section~\ref{sec:main-theorems} we are going to use the preceding
proposition together with lemma~\ref{l-27} 
to show that under suitable assumptions, the poset
$\poset(\F,\club(\downcl\F,S))$ is $\E(\oone\setminus
S)$-\alphaproper. We remark that by standard counterexamples, it is not
in general \alphaproper.
\end{rem}

\begin{lem}
\label{l-14}
Let\/ $\theta$ be an ordinal of uncountable cofinality.
Let $\F$ be a directed subfamily of $(\power(\theta),\subseteqfnt)$ with
$\J(\F,\subseteqfnt)$ a $\sigma$-ideal, and let $S\subseteq\theta$ be stationary. 
Suppose there is no stationary subset of\/ $S$ orthogonal to $\F$. 
Then\/ $\club(\downcl\F,S)$ is $\F$-extendable.
\end{lem}
\begin{proof}
Take $x\in\club(\downcl\F,S)$, some countable family $\X$ of cofinal subsets
of $\F$ with $x\subseteq y$ for all $y\in X$ for all $X\in\X$, and $\xi<\theta$.
Since $S$ is stationary, there exists a countable $M\prec H_\kappa$ 
with $x,\F,\X,\xi,S\in M$ and $\sup(\theta\cap M)\in S$. 
And by lemma~\ref{l-23}, $\{y\in X:\sup(\theta\cap M)\in y\}$ is
cofinal in $\F$  for all $X\in\X$. 
Since $x\cup\{\sup(\theta\cap M)\}\subseteq S$ is also closed and in $\downcl\F$,
and since $\xi<\sup(\theta\cap M)$, $x\cup\{\sup(\theta\cap M)\}$
witnesses extendability.
\end{proof}

Complete no longer has a winning strategy in the purely combinatorial
game (as it did with $\H$ convex in $\subseteq$), 
but is contented with a nonlosing strategy in~$\ggen^*$, for some models.

\begin{lem}
\label{l-15}
Let $\F$ be a directed subfamily $(\power(\theta),\subseteqfnt)$ with
$\J(\F,\subseteqfnt)$ a $\lambda$-ideal, and let $S\subseteq\theta$ be stationary. 
Suppose $M\prec H_{\theta^+}$ has cardinality $|M|<\cof(\theta)$ with $\F,S\in M$,
\begin{equation}
  \label{eq:31}
  \sup(\theta\cap M)\notin S
\end{equation}
and\/ $p\in\poset(\F,\club(\downcl\F,S))\cap M$. 
Assuming that\/ $\club(\downcl\F,S)$ is $\F$-extendable,
Complete has a nonlosing strategy for the game\/ 
$\ggen^*(M,y,\F,\club(\downcl\F,S),p)$
whenever\/ $y\subseteq M$ is  in $\downcl\F$ and of
cardinality\/ $|y|<\lambda$ with\/ $x_p\subseteq y$. 
\end{lem}
\begin{proof}
Set $\delta=\sup(\theta\cap M)$.
We know that $\bigcup_{k<\omega}\X_{p_k}$ will be countable, 
and thus we can arrange a diagonalization in advance. 
After Extender plays its $k\Th$ move $(p_k,X_k)$, 
we take care of some $Y_k\in\X_{p_k}$
according to the diagonalization.  
Since $Y'_k=\{x\in Y_k:x_{p_k}\subseteq x\}$ is cofinal in $\F$,
$Y''_k=\{x\in Y'_k:y\subseteqfnt x\}$ is $\subseteqfnt$-cofinal since
$y\in\downcl\F$ and $\F$ is $\subseteqfnt$-directed. Now
$Y''_k=\bigcup_{s\in\Fin(y\setminus x_{p_k})}\{x\in Y'_k:y\setminus s\subseteq x\}$,
and hence by the assumption on $\J(\F)$
there exists $s_k\in\Fin(y\setminus x_{p_k})$ such that
\begin{equation}
\label{eq:36}
  Y'''_k=\{x\in Y'_k:y\setminus s_k\subseteq x\}\text{ is $\subseteqfnt$-cofinal}.
\end{equation}
Exactly the same argument shows that there exists $t_k\in\Fin(y\setminus x_{p_k})$
such that
\begin{equation}
  \label{eq:27a}
  X'_k=\{x\in X_k:x_{p_k}\subseteq y\setminus t_k\subseteq x\}
  \text{ is $\subseteqfnt$-cofinal}.
\end{equation}
Complete plays $s_k\cup t_k$ on move $k$. 

If Extender loses then Complete wins. 
Thus we may assume that Extender does not lose.
Put $x_q=\bigcup_{k<\omega} x_{p_k}$ and
$\X_q=\bigcup_{k<\omega}\X_{p_k}\cup\{X_0,X_1,\dots\}$.
By proposition~\ref{p-12}, for every $\xi<\theta$ in $M$, $\D_\xi\in
M$ (cf.~equation~\eqref{eq:2}) is dense, and
thus $x_{p_k}\in\D_\xi$ for some $k$ since Extender did not lose.
Hence $x_q$ is unbounded in~$\delta$, 
and therefore $x_q\in\club(\downcl\F,S)$ for the same reason as in the
proof of proposition~\ref{p-14}.
And for every $k$, by proposition~\ref{p-6},
$x_q\subseteq x$ for all $x\in X'_k$ and all $x\in Y'''_k$ proving that
$q=(x_q,\X_q)\in\poset(\F,\club(\downcl\F,S))$ 
(see~\eqref{eq:27a} and~\eqref{eq:36}).
Thus the game is  drawn.
\end{proof}

When $\sup(\theta\cap M)\in S$, Complete no longer has a nonlosing
strategy in the augmented game, but the unaugmented game is still OK.

\begin{lem}
\label{l-6}
Let\/ $\F$ be a directed subfamily of\/ $(\power(\theta),\subseteqfnt)$ with\/
$\J(\F,\subseteqfnt)$ a $\lambda$-ideal, 
and let\/ $S\subseteq\theta$ be a stationary set.
Suppose\/ $M\prec H_{\theta^+}$ has cardinality\/ $|M|<\cof(\theta)$ 
with\/ $\F,S\in M$, and $p\in\poset(\F,\club(\downcl\F,S))\cap M$.
Assuming that there is no stationary subset of $S$ orthogonal to $\F$,
Complete has a nonlosing strategy for the game\/ 
$\ggen(M,y,\F,\club(\downcl\F,S),p)$
whenever\/ $y\subseteq M$ and 
$y\in\downcl\F$ is of cardinality\/ 
$|y|<\lambda$ with\/ $x_p\subseteq y$.
\end{lem}
\begin{proof}
Set $\delta=\sup(\theta\cap M)$. The case $\delta\notin S$ has been dealt with in
lemma~\ref{l-15}, by lemma~\ref{l-14} and proposition~\ref{p-25}. 
Assume then that $\delta\in S$.
We know that $\bigcup_{k<\omega}\X_{p_k}$ will be countable, and thus we can arrange a
diagonalization in advance. After
Extender plays its $k\Th$ move $p_k$, we take care of some $Y_k\in\X_{p_k}$
according to the diagonalization.  
Since $\{x\in Y_k:x_{p_k}\subseteq x\}\in M$, and is cofinal in $\F$, by lemma~\ref{l-23}, 
\begin{equation}
  \label{eq:7}
  Z_k=\bigl\{x\in Y_k:x_{p_k}\cup\{\delta\}\subseteq x\bigr\}
  \espc\text{is $\subseteqfnt$-cofinal in $\F$}.
\end{equation}
And thus
$Z'_k=\{x\in Z_k:y\subseteqfnt x\}$ is $\subseteqfnt$-cofinal since
$y\in\downcl\F$ and $\F$ is $\subseteqfnt$-directed. Now
$Z'_k=\bigcup_{s\in\Fin(y\setminus x_{p_k})}
  \{x\in Z_k:y\setminus s\subseteq x\}$,
and hence by the assumption on $\J(\F)$
there exists $s_k\in\Fin(y\setminus x_{p_k})$ such that
\begin{equation}
\label{eq:34}
  Z''_k=\{x\in Z_k:y\setminus s_k\subseteq x\}
  \espc\text{is $\subseteqfnt$-cofinal}.
\end{equation}
Complete plays $s_k$ on move $k$. This describes the strategy for Complete.

At the end of the game, assume without loss of generality that
Extender does not lose. Put
\begin{equation}
  \label{eq:1}
  x_q=\bigcup_{k<\omega} x_{p_k}\cup\{\delta\},
\end{equation}
and $\X_q=\bigcup_{k<\omega}\X_{p_k}$.
Every $Y\in\bigcup_{k<\omega}\X_{p_k}$ appears as $Y_k$ for some $k$, 
and thus $\{x\in Y:x_{q}\subseteq x\}$ is $\subseteqfnt$-cofinal 
by~\eqref{eq:34}. Since $\club(\downcl\F,S)$ is $\F$-extendable, 
$x_q\setminus\{\delta\}$ is unbounded in $\delta$, which clearly
implies that $x_q$ is closed in $S$. 
Hence $x_q\in\club(\downcl\F,S)$ as any of the sets $x$
from~\eqref{eq:34} witnesses that $x_q\in\downcl\F$.
This proves that $q=(x_q,\X_q)\in\poset(\F,\club(\downcl\F,S))$, 
and thus Complete does not lose the game.
\end{proof}

\begin{cor}{0}
\label{c-5}
Let $\F$ be a $\sigma$-directed subfamily of\/ $(\powcnt\theta,\subseteqfnt)$,
and let $S\subseteq\theta$ be a stationary set.
Suppose that there is no stationary subset of\/ $S$ orthogonal to $\F$.
Let $\psi$ be a formula describing a suitability function so that
 $\psi(M,y,\F,\H,p)$ implies $y\in\downcl\F$ and $x_p\subseteq y$. 
Then $\psi$-globally, Complete has a forward nonlosing strategy in the game
$\ggen(\F,\club(\downcl\F,S))$.
\end{cor}
\begin{proof}
Fix $M\in\S$ where $\S$ comes from the global
strategy associated with $\psi$, and suppose $p\in\poset(\F,\club(\downcl\F,S))\cap M$
and $\varphi(\F,\club(\downcl\F,S),p,a_3)$.
Then there exists $y\subseteq M$ satisfying $\psi(M,y,\F,\club(\downcl\F,S),p)$. 
It remains to show that Complete has a forward nonlosing strategy in the game
$\ggen(M,y,\F,\club(\downcl\F,S),p)$. Note that (by definition)
$\poset(\F,\C(\downcl\F,S))\in M$ and we may assume that $S$ can be
computed from this poset. 

We proceed as in the proof of corollary~\ref{c-4}, and thus we need 
a pair $(\E,F)$ to be used with lemma~\ref{l-16} 
(of course different than $(\S,\psi)$ above).
Put $\E=\{M\}$, $\T=\poset(\F,\club(\downcl\F,S))$ 
and let $F(M)$ be the function with
domain $\T\cap M$ given by $F(M)(q)=\{x\subseteq M:x\in\downcl\F$ and
$x_q\subseteq x\}$. Then for every
$q\in\poset(\F,\club(\downcl\F,S))\cap M$ and every $x\in F(M)(q)$,
Complete has a nonlosing strategy in the game
$\ggen(M,x,\F,\club(\downcl\F,S),q)$ by lemma~\ref{l-6}. 
Therefore, there is a nonlosing strategy for Complete
in the parameterized game $\game(M,x,q)$, $\E$-$F$-globally,
where $\game(M,x,q)=\ggen(M,x,\F,\allowbreak\club(\downcl\H,S),q)$. 

Supposing $P=(p_0,s_0),\dots,(p_k,s_k)$ is a position in the game
$\ggen(M,x,\allowbreak\F,\allowbreak\club(\downcl\F,S),\allowbreak q)$
with Extender to play, if we put $x'=x\setminus\bigcup_{i=0}^k s_i$ and
$q'=p_k$, clearly $x'\in F(M)(q')$ and
$\ggen(M,x',\F,\club(\downcl\F,S),q')=\ggen(M,x,\F,\club(\downcl\F,S),q)
\restriction P$ is immediate
from the rules of the game. Thus the lemma~\ref{l-16} applies with
$A=\mathrm{Extender}$ and $X=\mathrm{Complete}$, 
and hence Complete has a forward nonlosing strategy, $\E$-$F$-globally. 
In particular, Complete has
a forward nonlosing strategy in the game
$\ggen(M,y,\F,\club(\downcl\F,S),p)$, where $y$
and $p$ are the above fixed parameters.
\end{proof}

\setcounter{@temp}{-1}

\begin{cor}[1]{1}
\label{c-7}
Let $\F$ be a\/ $\sigma$-directed subfamily of\/ $(\powcnt\theta,\subseteqfnt)$,
and let $S\subseteq\theta$ be a stationary set.
Suppose that there is no stationary subset of\/ $S$ orthogonal to $\F$.
Let $\psi$ be a formula describing a suitability function so that 
 $\psi(M,y,\F,\H,p)$ implies $y\in\downcl\F$ and $x_p\subseteq y$. 
Then $\E(\theta\setminus S,\theta)$-$\psi$-globally, 
Complete has a forward nonlosing strategy in the game $\ggen^*(\F,\club(\downcl\F,S))$. 
\end{cor}
\begin{proof}
This can be proved very similarly to corollary~\ref{c-5}. 
The difference is that the game $\ggen^*$ is used in place of $\ggen$,
but we only need to consider $M\in\S$ with $\sup(\theta\cap M)\notin S$.
Thus we can invoke lemma~\ref{l-15} instead of lemma~\ref{l-6}. 
At position
$P=\bigl((p_0,X_0),s_0\bigr),\dots,\bigl((p_k,X_k),s_k\bigr)$ of the
game $\ggen^*(M,x,\F,\club(\downcl\F,S),q)$,
set $\bar x=x\setminus\bigcup_{i=0}^k s_i$ and $q'=p_k$. 
Then as in equation~\eqref{eq:27a}, there is a finite subset $z\subseteq\bar
x\setminus x_{p_k}$ such that putting $x'=\bar x\setminus z$, 
$\{w\in X_i:x'\subseteq w\}$ is $\subseteqfnt$-cofinal in $\F$ for all $i=0,\dots,k$.
Since $z$ is finite we can assume that $z$ is $\subseteq$-minimal with
this property. Now
$\ggen^*(M,x',\F,\club(\downcl\F,S),q')$ is `isomorphic' to
$\ggen^*(M,x,\F,\club(\downcl\F,S),q)\restriction P$: By the
minimality of $z$, the two games are identical for Extender; Complete
has less room to move in the former game, but the difference has no
effect on the outcome of the game. We will not formalize this argument
further, and in any case corollary~\ref{c-7} is not applied in this
paper.
\end{proof}

\begin{defn}
\label{d-23}
We let $\psicls(M,y,\F,\H,p)$  be the conjunction of
$\psimin(M,\allowbreak y,\allowbreak \F,\allowbreak\H,\allowbreak p)$ 
and $y\in\downcl\F$. 
\end{defn}

\begin{rem}
\label{r-15}
When $\F$ is a $\sigma$-directed subfamily of
$(\powcnt\theta,\subseteqfnt)$ then we can make statements ``$\psicls$-globally''
in that equation~\eqref{eq:39} holds.
\end{rem}

\begin{cor}[2]{1}
\label{c-6}
Let $\F$ be a $\sigma$-directed subfamily of
$(\powcnt\theta,\subseteqfnt)$, and let $S\subseteq\theta$ be a
stationary set with no stationary subset orthogonal to $\F$. 
Then  $\psicls$-globally, Complete has a
forward nonlosing strategy in the game $\ggen(\F,\club(\downcl\F,S))$. 
\end{cor}
\begin{proof}
Corollary~\ref{c-5} and remark~\ref{r-15}. 
\end{proof}

\setcounter{@temp}{-1}

\begin{cor}[3]{1}
\label{c-8}
Let $\F$ be a $\sigma$-directed subfamily of
$(\powcnt\theta,\subseteqfnt)$, and let $S\subseteq\theta$ be a
stationary set with no stationary subset orthogonal to $\F$. 
Then $\E(\theta\setminus S,\theta)$-$\psicls$-globally, 
Complete has an nonlosing
strategy in the game $\ggen^*(\F,\club(\downcl\F,S))$. 
\end{cor}
\begin{proof}
Corollary~\ref{c-7} and remark~\ref{r-15}.
\end{proof}

In the following three propositions
(propositions~\ref{p-19}--\ref{p-37}) we are restricting $\T$ only to
include $\F$ that are $\sigma$-directed subfamilies of
$(\powcnt\theta,\subseteqfnt)$. 

\begin{prop}
\label{p-19}
$\psicls$ is provably coherent.
\end{prop}
\begin{proof}
The exact same proof as for proposition~\ref{p-21} works, because the
element $y\cap y'\cap M$ defined there is clearly in $\downcl\F$.
\end{proof}

\begin{prop}
\label{p-27}
$\psicls$ provably respects $\ggen$.
\end{prop}
\begin{proof}
Immediate from the definition.
\end{proof}

\begin{prop}
\label{p-37}
Restricting to $\theta=\oone$, $\psicls$ provably
respects isomorphisms.
\end{prop}
\begin{proof}
Exactly the same as for proposition~\ref{p-35}.
\end{proof}

\section{Hybrid and combinatorial principia}
\label{sec:main-theorems}

The following hybrid principium is the strongest statement of its type considered
here.  

\begin{principle}{$\hybridmax$}
Let $(\F,\H)$ be a pair of subfamilies of $[\theta]^{\leq\aleph_0}$ 
for some ordinal $\theta$, with $\F$ closed under finite reductions. 
If $\H$ is $\F$-extendable and Extender has no winning strategy in the
parameterized game $\ggen(\F,\H)$, $\psi$-globally for some $\psi\to\psimin$, 
then there exists an uncountable
$X\subseteq\theta$ such that every proper initial segment of $X$ 
is in $\downcl(\H,\sqsubseteq)$. 
\end{principle}

\noindent The corresponding combinatorial principium is as follows.

\begin{principle}{$\pcombmax$}
Let $(\F,\H)$ be a pair of subfamilies of $[\theta]^{\leq\aleph_0}$ for some ordinal
$\theta$, with $\F$ closed under finite reductions. 
If $\H$ is $\F$-extendable and Complete has a winning strategy in the
parameterized game $\gcomp(\F,\H)$, $\psi$-globally for some $\psi\to\psimin$, 
then there exists an uncountable
$X\subseteq\theta$ such that every proper initial segment of $X$ is in
$\downcl(\H,\sqsubseteq)$. 
\end{principle}
These principles are mentioned because they have enough constraints to be
consistent. 

\begin{thm}
\label{u-9}
$\pfa$ implies $\hybridmax$.
\end{thm}
\begin{proof}
Let $\F$ and $\H$ be as specified in the principle. 
By corollary~\ref{o-5}, $\poset(\F,\H)$ is completely proper. And by
proposition~\ref{p-12}, $\D_\xi$ is dense for each $\xi<\theta$.

We still need a density argument to produce the desired uncountable
$X\subseteq\theta$. For each $\alpha<\oone$, put
\begin{equation}
  \label{eq:56}
  \E_\xi=\{p\in\poset(\F,\H):\otp(x_p)>\xi\}.
\end{equation}
Observe that each $\E_\xi$ is dense: Given $\xi<\oone$ and
$p\in\poset(\F,\H)$, take a countable elementary
submodel $M\prec H_\kappa$ with $\xi,p,\F,\H\in M$. By complete properness,
there exists $G\in\Genc(M,P,p)$. Since the $\D_\xi$'s are dense,
\begin{equation}
  \label{eq:57}
  M[G]\models\bigcup_{q\in G}x_q\text{ is cofinal in }\theta.
\end{equation}
Then taking the transitive collapse,
$\overbar M[\overbarg G]\models\bigcup_{q\in\overbarg
  G}x_q\text{ is cofinal in }\overbarg\theta$.
As $\xi\in M$ is fixed under the collapse, $\xi<\overbarg\theta$,
and thus there exists $q\ge p$ in $G$ with
$\otp(x_q)>\xi$. Hence $q\in\E_\xi$.

Now applying $\pfa$ to the
proper poset $\poset(\F,\H)$, it has a filter
$G$ intersecting $\E_\xi$ for all $\xi<\oone$.
Therefore $X_G=\bigcup_{p\in G}x_p\subseteq\theta$ is uncountable.
And every proper initial segment 
 $y\sqsubset X$ is in $\downcl{(\H,\sqsubseteq)}$, as required
(proposition~\ref{p-31}).
\end{proof}

A slight formal weakening of these principles, 
namely requiring a \emph{forward} strategy of Complete,
allows us to significantly weaken $\pfa$
in the hypothesis. The hybrid version $\hybrid$ has already been
presented in the introduction (page~\pageref{principle:hybrid}),
and following is the corresponding combinatorial principle.

\begin{principle}{$\pcomb$}
Let $(\F,\H)$ be a pair of subfamilies of $[\theta]^{\leq\aleph_0}$ 
for some ordinal $\theta$, with $\F$ closed under finite reductions. 
If $\H$ is $\F$-extendable and Complete has a forward winning strategy in the
parameterized game $\gcomp(\F,\H)$, $\psi$-globally for some $\psi\to\psimin$,
then there exists an uncountable
$X\subseteq\theta$ such that every proper initial segment of $X$ is in~$\downcl(\H,\sqsubseteq)$. 
\end{principle}

Recall that $\mathbb D$-completeness implies complete properness
(proposition~\ref{p-2}). 

\begin{thm}
\label{u-2}
$\ma($\tu{$\mathbb D$-complete}$)$ implies $\hybrid$.
\end{thm}
\begin{proof}
The additional requirement that $\psi$-globally, 
Complete has a forward
nonlosing strategy in the game $\ggen(\F,\H)$, guarantees that
$\poset(\F,\H)$ is $\mathbb D$-complete by lemma~\ref{l-7} and
remark~\ref{r-10}. The
rest of the proof is the same as the proof of theorem~\ref{u-9}.
\end{proof}

Just in case one decides to do a more in depth study of such
principles, we might ask whether this weakening is purely formal. 

\begin{question}
\label{q-1}
Does either $\hybrid\to\hybridmax$ or $\pcomb\to\pcombmax$\tu?
\end{question}

Our goal here is to obtain principles compatible with $\ch$. Although the
``medicine'' against (destroying) 
weak diamond has been taken for, e.g., $\pcomb$, no ``medicine''
has been taken for the so called disjoint clubs (cf.~\cite{math.LO/0003115}). 
Thus we would be surprised if
$\pcomb$ is consistent with $\ch$. On the other hand, we do not know of a
counterexample.

\begin{question}
\label{q-2}
Is $\pcomb$ compatible with $\ch$\tu?\footnote{We did not use the word
  ``consistent'' because we want to avoid large cardinal considerations for the moment.}
\end{question}

The only way we know of to take the latter ``medicine'' is to make `geometrical'
restrictions on the the second family $\H$. First the hybrid principle:

\begin{principle}{$\hybridboc$}
Let $\F$ be a subfamily of $[\theta]^{\leq\aleph_0}$ for some ordinal $\theta$,
closed under finite reductions. 
Suppose $\H$ is an upwards boundedly order closed subfamily of
$([\theta]^{\leq\aleph_0},\sqsubseteq)$.
If $\H$ is $\F$-extendable and Complete has a forward 
nonlosing strategy in the
parameterized game $\ggen(\F,\H)$, $\psi$-globally for some $\psi\to\psimin$
that is coherent and respects $\ggen$, 
then there exists an uncountable
$X\subseteq\theta$ such that every proper initial segment of $X$ is in~$\downcl(\H,\sqsubseteq)$. 
\end{principle}

\noindent And the combinatorial principle:

\begin{principle}{$\pcombboc$}
Let $\F$ be a subfamily of $[\theta]^{\leq\aleph_0}$ for some ordinal $\theta$,
closed under finite reductions. 
Suppose $\H$ is an upwards boundedly order closed subfamily of
$([\theta]^{\leq\aleph_0},\sqsubseteq)$.
If $\H$ is $\F$-extendable and Complete has a
forward winning strategy in the
parameterized game $\gcomp(\F,\H)$, $\psi$-globally for some $\psi\to\psimin$ that is
coherent and respects $\ggen$, then there exists an uncountable
$X\subseteq\theta$ such that every proper initial segment of $X$ 
is in~$\downcl(\H,\sqsubseteq)$. 
\end{principle}

\begin{thm}
\label{u-3}
$\ma($\tu{\alphaproper and $\mathbb D$-complete}$)$ implies $\hybridboc$.
\end{thm}
\begin{proof}
The only thing that needs to be added to the proof of
theorem~\ref{u-2} is that $\poset(\F,\H)$ is \alphaproper.
And this is by lemma~\ref{l-4} because $\H$ is upwards boundedly order
closed and by the additional requirements on $\psi$.
\end{proof}

\begin{corthm}
\label{o-8}
$\hybridboc$ is compatible with $\ch$. More precisely, $\hybridboc$ is consistent
with $\ch$ relative to a supercompact cardinal. 
\end{corthm}
\begin{proof}
Theorem~\ref{u-3} and Shelah's Corollary on page~\pageref{Shelah:ma-alpha}.
\end{proof}

As we shall see below, $\pcombboc$ implies $\pstar$ and thus has considerable large
cardinal strength (cf.~page~\pageref{large-cardinal}). 
However, if we want to restrict to
$\oone$ then we expect that no large cardinals are necessary, as is the case 
with~$\pstar_\oone$. We prove this for a slight weakening of $\hybridboc$.

\begin{principle}{$\hybridbociso$}
This is the same principle as $\hybridboc$ except that we add the requirement that
$\psi$ respects isomorphisms.
\end{principle}

\noindent The combinatorial principle $\pcombbociso$ is exactly analogous.

\begin{thm}
\label{u-12}
$\ma(\tu{\alphaproper}+\tu{\deecmp}+\tu{pic}
+\tu{$\Delta_0$-$H_{\aleph_2}$-definable})$ implies $\hybridbociso_\oone$.
\end{thm}
\begin{proof}
Letting $\F$ and $\H$ be subsets of $[\oone]^{\le\aleph_0}$ 
satisfying the requirements of the principle, we need in
addition to the proof of theorem~\ref{u-3} to show that
$\poset(\F,\H)$ satisfies the properness isomorphism condition and is
$\Delta_0$-definable over $H_{\aleph_2}$. These are true by
corollary~\ref{o-7} and example~\ref{x-6}, respectively. 
\end{proof}

\begin{corthm}
\label{u-4}
$\hybridbociso_{\oone}$ is consistent with $\ch$ relative to the consistency of\/~$\zfc$. 
\end{corthm}
\begin{proof}
By theorem~\ref{u-13}.
\end{proof}

This raises an interesting question. ``Naturally occurring'' 
combinatorial principles on $\oone$ generally
(always?) have no large cardinal strength. The only explanation that the author
knows of is that they can be forced without collapsing cardinals. 
Perhaps there is a more satisfactory explanation? (It is quite possible
that there is a well known explanation that the author is simply
unaware of.) If all naturally occurring
combinatorial principles on $\oone$---notice that this notion has not
been clearly defined, and this may well be the essential point---can
be decided without large cardinals, then we should be able to prove
the consistency of $\pcombmax_\oone$ without large cardinal assumptions.

\begin{question}
\label{q-4}
Is $\pcombmax_\oone$ consistent relative to $\zfc$\tu?
\end{question}

Even if question~\ref{q-4} has a positive answer, it is conceivable that
adding $\ch$ requires large cardinals.

\begin{question}
\label{q-5}
Is the conjunction of\/ $\pcombboc_\oone$ and\/~$\ch$ relatively
consistent with\/ $\zfc$\tu?
If question~\tu{\ref{q-2}} has a positive answer, is $\pcomb_\oone$ and\/ $\ch$
relatively consistent with $\zfc$\tu?
\end{question}

\noindent Note that we can obtain the consistency of $\pcombmax_\oone$ 
or $\pcombboc_\oone+\ch$
by assuming the existence of an inaccessible cardinal. 
E.g.~one can produce a model of
$\pcombboc_\oone$ and $\ch$ by iterating up to an inaccessible cardinal instead of
using the properness isomorphism property to guarantee a suitable 
chain condition. 

Now we observe that the combinatorial version $\pcombboc$ is already
weak enough to be consistent with the existence of a nonspecial
Aronszajn tree. 

\begin{thm}
\label{u-5}
The conjunction of\/ $\pcombboc$, $\ch$ and the existence of a 
nonspecial Aronszajn
tree is consistent relative to a supercompact cardinal.
\end{thm}
\begin{proof}
By the argument in the proof of theorem~\ref{u-9}, 
to obtain a model of $\pcombboc$ it suffices to ensure that for every
pair $(\F,\H)$ as specified in the principle, for every family of size
$\aleph_1$ consisting of dense subsets of the poset $\poset(\F,\H)$,
there exists a filter intersecting every member of the family, i.e.~we
prove $\ma(\poset(\F,\H))$ for every such $(\F,\H)$. 

Using the argumentation of the proof of the consistency of $\pfa$ (see
e.g. \cite[\Section1]{MR924672}), we can obtain such a model by extending
by a countable support iteration $(P_\xi,\dot
Q_\xi:\xi<\kappa)$ where $\kappa$ is a supercompact cardinal and each iterand $\dot Q_\xi$ is either a $P_\xi$-name for a poset of the form
$\poset(\dot\F_\xi,\dot\H_\xi)$ with $(\dot\F_\xi,\dot\H_\xi)$ as
above, or else a $P_\xi$-name for a L\'evy collapse of the form
$\coll(\aleph_1,\theta)$ which in particular is a $\sigma$-closed poset. 
As argued in the proof of theorem~\ref{u-3}, each $\dot Q_\xi$ is \alphaproper, 
and in the case $\dot Q_\xi=\poset(\dot\F_\xi,\dot\H_\xi)$, $\poset$
has the $\sigma$-complete completeness system from remark~\ref{r-10},
and the collapsing poset has the trivial completeness system
(i.e.~$\Gen(M,P)=\Genc(M,P)$ whenever $P$ is $\sigma$-closed).
Therefore, by Shelah's Theorem on 
page~\pageref{Shelah:alpha}, the iteration does not add new reals and
thus if the ground model satisfies $\ch$ then so does the extension.

Now it is specified (in particular) that $\psi$-globally, Complete has a 
winning strategy in the game $\gcomp(\F,\H)$. 
Hence Complete has, $\psi$-globally, a nonlosing strategy
in the game $\ggen^*(\F,\H)$ by proposition~\ref{p-13}. 
Thus by lemma~\ref{l-11}, the poset $\poset(\F,\H)$ is
$\oone$-tree-preserving. And $\sigma$-closed posets are also
$\oone$-tree preserving by lemma~\ref{l-31}. Therefore, every iterand of our
iteration is $\oone$-tree-preserving, and thus
$\varinjlim{}_{\xi<\kappa}P_\xi$ is $\oone$-tree-preserving 
by Schlindwein's theorem~\cite{MR1264960} discussed on
page~\pageref{page:Sch94}.
Thus by lemma~\ref{p-22}, if we begin with a ground model satisfying $\ch$
containing a Souslin tree $T$ then our extension satisfies $\pcombboc$,
$\ch$ and $T$ remains nonspecial. To ensure that $T$ moreover remains
Aronszajn in our extension, we can further assume that $T$ is
$h$-st-special in the ground model (this is better explained in the
proof of theorem~\ref{u-14}). 
\end{proof}

\begin{thm}
\label{u-6}
The conjunction of\/ $\pcombbociso_\oone$, $\ch$ 
and the existence of a nonspecial Aronszajn
tree is consistent relative to\/ $\zfc$.
\end{thm}
\begin{proof}
This is done very similarly to the proof of theorem~\ref{u-13}, except
that we only include iterands of the form $\poset(\dot\F,\dot\H)$,
where $(\dot\F,\dot\H)$ names a pair of subfamilies of~$\powcnt\oone$ 
as in the specification of
$\pcombbociso$. Thus, as in the preceding theorem, each iterand is
$\oone$-tree-preserving in addition to being \alphaproper and has a
fixed $\sigma$-complete completeness system and is in pic. 
Therefore, by starting out with a ground model
containing an $h$-st-special Souslin tree and satisfying $\gch$, we end up with a
model satisfying $\pcombbociso$, $\ch$ and the existence of a
nonspecial Aronszajn tree.
\end{proof}

\begin{lem}
\label{l-25}
$\pcombboc$ implies $\pstar$.
\end{lem}
\begin{proof}
Let $\H$ be a $\sigma$-directed subfamily of
$(\powcnt\theta,\subseteqfnt)$. We assume that the second
alternative~\eqref{item:51} of $\pstar$ fails, and prove that the
first alternative~\eqref{item:50} holds.
Next we verify that the pair
$(\downcl\H,\downcl\H)$ satisfies the hypotheses of the principle~$\pcombboc$.
Since $\downcl\H$ is convex, it is
upwards boundedly order closed by proposition~\ref{p-23}. 
By restricting to a smaller ordinal if necessary, we can assume that
$\theta$ is least ordinal that has no countable decomposition into
pieces orthogonal to $\H$. 
Then by lemma~\ref{l-13}, $\downcl\H$ is extendable (i.e.~$\downcl\H$ is $\downcl\H$-extendable). 
We use $\psicvx$ from definition~\ref{d-22}. 
Complete has a forward winning strategy in the game
$\gcomp(\downcl\H,\downcl\H)$, $\psicvx$-globally, by corollary~\ref{c-3}. And
$\psicvx\to\psimin$ is coherent and respects $\ggen$ by
propositions~\ref{p-32} and~\ref{p-34}. Therefore, $\pcombboc$ gives
an uncountable $X\subseteq\theta$ with every proper initial segment in
$\downcl\H$. In particular, $X$ is locally in~$\downcl\H$ establishing
the first alternative.
\end{proof}

\begin{lem}
\label{l-3}
$\pcombbociso_\oone$ implies $\pstar_\oone$.
\end{lem}
\begin{proof}
Additionally to the proof of lemma~\ref{l-25}, restricting
$\F,\H\subseteq\powcnt{\oone}$ implies that $\psicvx$ respects
isomorphisms by proposition~\ref{p-35}. 
\end{proof}

Now we have answered the Abraham--Todor\v cevi\'c question
(cf.~page~\pageref{q-3}), by establishing that $\pstar$ does not imply
that all Aronszajn trees are special (cf.~theorem~\ref{u-1}).

\begin{proof}[Proof of theorem~\tu{\ref{u-1}}]
Theorem~\ref{u-5} and lemma~\ref{l-25} for the unrestricted principle, 
and theorem~\ref{u-6} and lemma~\ref{l-3} for $\pstar_\oone$.
\end{proof}

\begin{thm}
\label{u-14}
The conjunction of $\pstarsplus_\oone$ and $\pstar_\oone$ is
consistent with $\ch$ and the existence of a nonspecial Aronszajn tree
relative to $\zfc$. 
\end{thm}
\begin{proof}
First we describe the ground model $V$.
Assume, by going to a forcing extension if necessary, that $\ch$ holds and
$2^{\aleph_1}=\aleph_2$. By further forcing if necessary, 
we may assume moreover that the stationary coideal $\NS^+$ does not satisfy
the $\aleph_2$-chain condition (see e.g.~\cite{MR1905156}). 
Thus we can fix a maximal antichain $\A\subseteq\NS^+$ of
cardinality~$\aleph_2$ (maximality is unimportant here).
By forcing yet again if necessary,
we may assume that in addition there exists a Souslin tree~$T$. 

Next we let $W$ denote the forcing extension of $V$ by the
$\oone$-tree-preserving (and thus proper) 
forcing notion from lemma~\ref{p-39}, so that in $W$, $T$ is
$h$-st-special for some $h$. 

In $W$: We construct a forcing notion $P=\varinjlim{}_{\xi<\omega_2}P_\xi$ 
using an iterated forcing construction 
$(P_\xi,\dot Q_\xi:\xi<\omega_2)$ of length $\omega_2$ 
with countable supports and inverse limits. Just as in the proof of
theorem~\ref{u-13}, and as is standard,
each iterand $\dot Q_\xi$ is forced to be of size
at most $\aleph_2$, and satisfy the pic. Thus, 
as each $P_\xi$ has the the $\aleph_2$-cc and a dense suborder of
cardinality at most $\aleph_2$, we
can use standard bookkeeping to obtain an enumeration 
$(\dot R_\xi,\dot\H_\xi:\xi<\omega_2)$ in advance such that, viewing
$P_\xi$-names as $P_\eta$-names for $\xi\le\eta$, 
every pair of $P_\xi$-names $(\dot R,\dot\H)$ for a subset $\dot R$ of
$\oone$ and a subfamily $\dot\H$ of
 $[\oone]^{\leq\aleph_0}$, respectively, appears as 
$(\dot R_\eta,\dot\H_\eta)$ for cofinally many $\xi\leq\eta<\omega_2$. 
By skipping steps, we may assume
that for all $\xi$,
\begin{enumerate}[leftmargin=*, label=(\roman*), ref=\roman*, widest=ii]
\item $P_\xi\forces\dot R_\xi\subseteq\oone$ is stationary,
\item $P_\xi\forces\dot\H_\xi$ is a $\sigma$-directed subset of
  $([\oone]^{\le\aleph_0},\subseteqfnt)$. 
\end{enumerate}

We also recursively choose $P_\xi$-names $\dot S_\xi$ for subsets of
$\oone$, according to which one of the following mutually exclusive cases holds:
\begin{enumerate}[leftmargin=*, label=(\alph*), ref=\alph*, widest=b]
\item\label{item:37} $P_\xi\forces\ulc\dot R_\xi=\check O_\xi\urc$ for
  some $O_\xi\in\A$,
\item\label{item:39} $p\forces\ulc\dot R_\xi\cap S\in\NS$ 
for all $S\in\A\cup\{\dot S_\gamma:\gamma<\xi\}\urc$ for some $p\in P_\xi$,
\item\label{item:38} $p\forces\ulc\dot R_\xi\in\{\dot
  S_\gamma:\gamma<\xi\}\setminus\A\urc$ for some $p\in P_\xi$ and
  case~\eqref{item:39} fails, 
\item\label{item:40} otherwise.
\end{enumerate}
In case~\eqref{item:37}, we set $\dot S_\xi=\check O_\xi$;  
in case~\eqref{item:39}, we define $\dot S_\xi$ so that $p\forces \dot
S_\xi=\dot R_\xi$ 
whenever $p\in P_\xi$ is as specified there, and
$q\forces\dot S_\xi=\check\emptyset$ for $q$ incompatible with such a
$p$; 
in case~\eqref{item:38}, we define $\dot S_\xi$ so that 
$p\forces\dot S_\xi=\dot R_\xi$ whenever $p\in P_\xi$ is as specified
in~\eqref{item:38}, and $q\forces\dot S_\xi=\check\emptyset$ if $q$ is
incompatible with any such $p$; in case~\eqref{item:40}, we set $\dot S_\xi=\check\emptyset$. 
Thus we have that the nonempty members of 
$\A\cup\{\dot S_\gamma:\gamma\leq\xi\}$ are forced to 
form an antichain of $\NS^+$. Put
  $\Gamma_\xi=\{\gamma\le\xi:P_\gamma\forces\dot S_\gamma=\check
  O_\gamma\}$. 
Then observe that
\begin{equation}
  \label{eq:55}
  P_\gamma\forces \check S\cap\dot S_\gamma\in\NS\espc\text{for all
    $S\in\A\setminus\Gamma_\xi$ and all $\gamma\le\xi$}.
\end{equation}
We shall
also ensure that at each stage $\xi$,
\begin{alignat}{1}
  \label{eq:50}
  P_\xi&\forces\dot Q_\xi\text{ is $\E(\oone\setminus\dot S_\xi)$-\alphaproper}\\
  \label{eq:52}
  P_\xi&\forces\dot Q_\xi\text{ is
    $(\omega_1\text{-tree},\dot S_\xi)$-preserving}
\end{alignat}

We now start working in the forcing extension by $P_\xi$, in order to
specify $\dot Q_\xi$. Let  $\H_\xi$ and  $S_\xi$ be the
interpretations of $\dot\H_\xi$ and  $\dot S_\xi$, respectively.
If $S_\xi=\emptyset$ we let $Q_\xi$ be trivial. 
Otherwise, we now specify $Q_\xi$, determined in order by the following cases.

\smallskip
\paragraph{\emph{Case} 1.} $\H_\xi$ has a stationary orthogonal set, 
but $\oone$ has no countable decomposition into sets orthogonal to
$\H_\xi$. Here we force with $Q_\xi=\poset(\downcl\H_\xi)$. 

First we verify that it forces the desired object.
Let $X_{\dot G}$ be a $Q_\xi$-name for the generic object. Clearly it
is forced to be locally in $\downcl\H_\xi$ (cf.~proposition~\ref{p-31}). 
Since $\oone$ is in fact the least ordinal that cannot be decomposed
in to countable many pieces orthogonal to $\H_\xi$
(proposition~\ref{p-3}), the fact that $\H_\xi$ is $\sigma$-directed
implies  that it is extendable by lemma~\ref{l-13}. 
Therefore $X_{\dot G}$ is forced to be uncountable (proposition~\ref{p-12}). 

By corollary~\ref{c-9}, $\psicvx$-globally, Complete has a forward winning
strategy in the game $\gcomp(\downcl\H_\xi,\downcl\H_\xi)$. 
Now $\downcl\H_\xi$ is upwards boundedly order closed
(proposition~\ref{p-23}), $\psicvx\to\psimin$ is coherent and respects
$\ggen$ (propositions~\ref{p-32} and~\ref{p-34}) and
$\psicvx$-globally, Complete has a
forward nonlosing strategy for $\ggen(\downcl\H_\xi,\downcl\H_\xi)$
(propositions~\ref{p-25} and~\ref{p-13}).
Therefore, $Q_\xi$ is in \deecmp by lemma~\ref{l-7} as witnessed by the
fixed pair of formulae given in remark~\ref{r-10}. And $Q_\xi$ is
$\aupalpha$-proper by lemma~\ref{l-4}, and in particular~\eqref{eq:50}
is satisfied. Also, since $\psicvx$ respects
isomorphisms (proposition~\ref{p-35}), $Q_\xi$ satisfies
the properness isomorphism condition by corollary~\ref{o-7}.

We verify that $Q_\xi$ is $\oone$-tree-preserving. 
Indeed by proposition~\ref{p-13}, $\psicvx$-globally,
Complete has a nonlosing
strategy in the game $\ggen^*(\downcl\H_\xi,\downcl\H_\xi)$; 
hence, we obtain the
$\oone$-tree-preserving property from lem\-ma~\ref{l-11}.

\smallskip
\paragraph{\emph{Case} 2.}
$\H_\xi$ has no stationary set orthogonal to it. In this case we force
with $Q_\xi=\poset(\downcl\H_\xi,\club(\downcl\H_\xi,S_\xi))$. 

By lemma~\ref{l-14}, $\club(\downcl\H_\xi,S_\xi)$ is $\H_\xi$-extendable, and
thus $X_{\dot G}$ is forced to be uncountable.
And $Q_\xi$ forces that every proper initial segment of $X_{\dot G}$ is in
$\downcl{(\club(\downcl\H_\xi,S_\xi),\sqsubseteq)}$ 
by proposition~\ref{p-31}.

By corollary~\ref{c-6}, $\psicls$-globally, 
Complete has a forward nonlosing strategy in the game
$\ggen(\downcl\H_\xi,\club(\downcl\H_\xi,S_\xi))$.
Thus $Q_\xi$ is in \deecmp by lemma~\ref{l-7}. 
Since $\club(\downcl\H_\xi,S_\xi)$ is upwards boundedly order closed
beyond $S_\xi$ (proposition~\ref{p-14}) and $\psicls\to\psimin$ is
coherent and respects $\ggen$ (propositions~\ref{p-19}
and~\ref{p-27}), Complete's
nonlosing strategy $\psi$-globally
ensures that $Q_\xi$ is $\E(\oone\setminus S_\xi)$-\alphaproper 
by lemma~\ref{l-27}, and thus~\eqref{eq:50} holds.  And since
$\psicls$ respects isomorphisms (proposition~\ref{p-37}), $Q_\xi$ is in
pic by corollary~\ref{o-7}. 

As for equation~\eqref{eq:52}, Complete has a
nonlosing strategy in the game $\ggen^*(\downcl\H_\xi,
\allowbreak\club(\downcl\H_\xi,S_\xi))$,  
$\E(\oone\setminus S_\xi)$-$\psicls$-globally, by corollary~\ref{c-8}. Hence $Q_\xi$ is
$(\omega_1$-tree$,S_\xi)$-preserving by lemma~\ref{l-11}.

\smallskip
\paragraph{\emph{Case} 3.} There is a countable decomposition of
$\oone$ into sets
orthogonal to $\H_\xi$. Then $Q_\xi$ is trivial.

\smallskip
Having defined the iteration, we verify that it has the desired
properties. In the final forcing extension $W[G]$ of $W$ by $P$: 
We know that the nonempty
members $\B$ of $\A\cup\{S_\xi:\xi<\omega_2\}$ form an antichain of
$\NS^+$. And $\B$ is in fact a
maximal antichain. For suppose to the contrary that there is a
stationary set $R$ with $R\cap S\in\NS$ for all $S\in\B$.
Since $P$ is proper and thus does
not collapse~$\aleph_1$ and since it has the $\aleph_2$-cc,
$R$ appears at some intermediate stage, i.e.~there exists
$\eta<\omega_2$ such that $R\in W[G_\eta]$, where
$G_\eta=G\restriction P_\eta$ is a generic filter on $P_\eta$.
Therefore by our bookkeeping, we can find $\xi<\omega_2$ such that 
$\dot R_\xi[G]=R$, and thus $p$ forces $\dot R_\xi\cap
S\in\NS$ for all $S\in\A\cup\{\dot S_\gamma:\gamma<\xi\}$, for some
$p\in G_\xi$. Thus we are in case~\eqref{item:39}, and hence 
$p\forces \dot S_\xi=\dot R_\xi$, and we
arrive at the contradiction~$R\in\B$.

Let us verify that $\B$ serves to instantiate
$\pstarsplus_\oone$. Suppose $\H$ is a $\sigma$-directed subfamily of
$(\powcnt\oone,\subseteqfnt)$, and that the second
alternative~\eqref{item:53} of $\pstarsplus$ fails, and hence there is
no stationary set orthogonal to $\H$. Take $S\in\B$. Now there exists
$\xi<\omega_2$ such that
\begin{equation}
(\H,S)=(\dot\H_\xi[G],\dot R_\xi[G]).\label{eq:53}
\end{equation} 
If $S\in\A$, then $p\forces\dot R_\xi=\check S$ for some $p\in G$, 
and hence by our bookkeeping, we can find such an $\xi$ so
that we are furthermore in case~\eqref{item:37}. 

Otherwise, $S\in\{S_\gamma:\gamma<\omega_2\}\setminus\A$, and thus
$p\forces\dot R_\xi=\dot S_\gamma\notin\A$ for some
$\gamma<\omega_2$ and  $p\in G_\xi$.
By our bookkeeping, we can find such an $\xi$ so that $\xi>\gamma$ 
and so that moreover we are not in case~\eqref{item:39}. Therefore, we
are in case~\eqref{item:38}. 

In either of these two situations, $p\forces\dot S_\xi=\dot R_\xi$ for
some $p\in G$. And we arrive at case~2 of the construction of
$\dot Q_\xi$. Hence $Q_\xi=\dot Q_\xi[G]$ forces an uncountable 
$X\subseteq \oone$ with all of its initial segments in
$\downcl(\club(\downcl\H,S),\sqsubseteq)$. In particular, $X$
relatively closed in $S$ and is locally in $\downcl\H$, as wanted.

The verification that $\pstar_\oone$ holds is similar. Here we will
arrive in either case~1 or~2 depending on whether there is a
stationary set orthogonal to $\H$, and in either case we obtain an
uncountable $X\subseteq\oone$ locally in $\downcl\H$.

Next we prove that $W[G]\models\ch$, by showing that $P$ does not add
any new reals. Suppose towards a contradiction that there is a real
number $r\in\reals$ (in~$W[G]$) that is not in the model $W$.
Then at some intermediate stage~$\eta$, $r\in W[G_\eta]$, and thus the
initial segment $P_\eta$ of the iteration has added a new real.
Since $|\A|=\aleph_2$ while $|\Gamma_\eta|\le\aleph_1$, there exists
$S\in\A\setminus\Gamma_\eta$. But for all $\xi<\eta$, $P_\xi\forces
\dot Q_\xi$ is $\E(S\setminus\dot S_\xi)$-\alphaproper, by
equation~\eqref{eq:50}, and $P_\xi\forces \check S\cap\dot
S_\xi\in\NS$ by equation~\eqref{eq:55}. 
Therefore, by lemma~\ref{l-28}, we in fact have that $P_\xi\forces \dot Q_\xi$ is
$\E(S)$-\alphaproper. Since we have also demonstrated that each
iterand is in \deecmp, by Shelah's Theorem on
page~\pageref{Shelah:alpha}, $P_\eta$ does not add new reals, a contradiction.

To conclude the proof, we need to show that the Souslin tree $T$ 
in the ground model $V$, is a nonspecial Aronszajn tree 
in the forcing extension $W[G]$. Since $W\models\ulc T$ is
$h$-st-special$\urc$, we know that $T$ is Aronszajn by proposition~\ref{p-38}.
Hence it remains to establish nonspecialness.
However, supposing towards a contradiction that it is special, the
specializing function appears in some intermediate model, and thus
$W[G_\eta]\models\ulc T$ is special$\urc$ for some $\eta<\omega_2$. 

By equation~\eqref{eq:52}, $P_\xi\forces\dot
Q_\xi$ is $(T,\dot S_\xi)$-preserving, for all $\xi<\eta$. 
And taking any $S\in\A\setminus\Gamma_\eta$, 
we have $P_\xi\forces S\cap\dot S_\xi\in\NS$ for all $\xi<\eta$ by~\eqref{eq:55}.
Therefore, $P_\xi\forces(\oone\setminus S)\cup\dot S_\xi
=(\oone\setminus S)\cup\dot A_\xi$,  where $P_\xi\forces\dot A_\xi\in\NS$.
Since every $P_\xi$ forces that $\dot Q_\xi$ is
$(T,(\oone\setminus S)\cup\dot S_\xi)$-preserving, 
by lemma~\ref{l-29}, $P_\xi\forces\dot Q_\xi$ is 
$(T,(\oone\setminus S))$-preserving. 
But then $P_\eta$ is $(T,(\oone\setminus S))$-preserving by lemma~\ref{l-30}. 
And therefore, since $W$ is an $\oone$-tree-preserving forcing extension of
$V$, $W[G_\eta]$ is also a $(T,(\oone\setminus S))$-preserving
extension of $V$. 
Hence, by lemma~\ref{p-22}, $W[G_\eta]\models\ulc\forces T$ is
nonspecial$\urc$, a contradiction.
\end{proof}

\begin{question}
\label{q-8}
Is the conjunction of $\pstarsplus$ and $\pstar$ consistent with $\ch$
and the existence of a nonspecial Aronszajn tree relative to a
supercompact cardinal\tu?
\end{question}

\begin{rem}
\label{r-3}
We expect that the proof of theorem~\ref{u-14} can be readily
generalized to prove the consistency of $\pstarsplus$ with $\ch$;
indeed, this should probably be stated and proved as a separate
theorem. In the proposed proof one would construct maximal antichains
in $(\NS_\theta^+,\subseteq)$ for each regular $\theta<\kappa$
simultaneously, where $\kappa$ is a supercompact cardinal. Then at any
intermediate stage $\xi<\kappa$ of the iteration, we would have that $P_\xi$ is
$\E(S,\theta)$-\alphaproper for some regular cardinal $\theta<\kappa$ 
and some stationary $S\subseteq\theta$; hence, the iteration would not
add new reals. Forcing $\pstar$ to hold simultaneously should not pose
any difficulties.

However, preserving the nonspecialness of some Aronszajn tree does
seem problematic, and seems to require a new approach even if
question~\ref{q-8} has a positive answer. The difficulty
is that if we try to preserve the $(T,R)$-preserving property for some
costationary $R\subseteq\oone$, we run out of possibilities for
$R$ after $2^{\aleph_1}$ many steps; unlike for
$\E(S)$-$\aupalpha$-properness, 
we cannot use $(T,R)$-preserving with $R\subseteq\theta>\oone$.
\end{rem}
\section{Applications}
\label{sec:applications}

\subsection{The principle $\princA$}
\label{sec:p_aleph_1}

\begin{defn}
\label{d-18}
We say that a lower subset $\L$ of some poset $(P,\leq)$ is
\emph{$\lambda$-generated} if there exists a subset $\G\subseteq\L$ of cardinality
at most $\lambda$ such that every member $x$ of $\L$ satisfies $x\leq y$ for some
member $y$ of $\G$, i.e.~$\L\subseteq\downcl{(\G,\le)}$.

We say that $\L$ is \emph{locally $\lambda$-generated} if for every
$a\in P$, $\downcl{\L}\cap\downcl a$ is $\lambda$-generated.
\end{defn}

\begin{prop}
\label{p-20}
Let $\lambda$ be an infinite cardinal.
An ideal is $\lambda$-generated in the standard sense iff it is $\lambda$-generated
as a lower set. 
\end{prop}

We shall consider lower subsets of posets of the form
$(\powcnt{\theta},\subseteq)$. 
For a family~$\F$ of sets and a set $A$, 
write $\F\cap A=\{x\cap A:x\in\F\}$. 

\begin{prop}
\label{p-28}
Let $\H\subseteq\powcnt\theta$. Then $\downcl\H$ is locally
$\lambda$-generated in $(\powcnt\theta,\subseteq\nobreak)$ 
iff $\downcl{\H\cap A}$ is $\lambda$-generated for every countable $A\subseteq\theta$. 
\end{prop}

The next lemma generalizes~\cite[Theorem~2.6]{EN}. We refer to two of
the standard cardinal characteristics of the continuum $\p$ and $\b$
(see~e.g.~\cite{Bl}). $\p$ is the smallest cardinality of a subfamily
$\A\subseteq[\omega]^{\aleph_0}$ 
with the \emph{finite intersection property},
meaning that $\bigcap F$ is infinite for every nonempty finite
$F\subseteq\A$, and such that $\F$ has no infinite
$\subseteqfnt$-lower bound. 

\begin{lem}[$\p>\lambda$]
\label{l-18}
Let $\H$ be a directed subset of $(\powcnt\theta,\subseteqfnt)$.
If $\downcl\H$ is locally $\lambda$-generated then
\begin{equation}
  \label{eq:29}
  \H^{\perp\perp}=\downcl(\H,\subseteqfnt).
\end{equation}
\end{lem}
\begin{proof}
We need to show that
$\H^{\perp\perp}\subseteq\downcl{(\H,\subseteqfnt)}$, 
because  the opposite inclusion holds for any family $\H$.
Suppose $y\in\powcnt\theta$ is not in $\downcl(\H,\subseteqfnt)$. 
By assumption, $\downcl\H\cap y$ is $\lambda$-generated, say by some family of
generators $\G\subseteq\H\cap y$ with $|\G|\leq\lambda$. Since $\H$ is directed, there can
be no finite $F\subseteq\G$ with $y\subseteqfnt\bigcup F$. Therefore, 
$y\setminus\G=\{y\setminus x:x\in\G\}$ is a family of cardinality at most $\lambda$ 
with the finite intersection property. 
Thus $\p>\lambda$ yields an infinite $z\subseteq y$ that is a $\subseteqfnt$-lower
bound of $y\setminus\G$. Since $\G$ generates $\H\cap y$, $z\in\H^\perp$; 
and since it is infinite, it witnesses that $y\notin\H^{\perp\perp}$. 
\end{proof}

The following is essentially the characterization of $\b$ as the least cardinal~$\kappa$ 
for which there exists an $(\omega,\kappa)$ gap in $\power(\omega)\div\Fin$. 

\begin{lem}[$\b>\lambda$]
\label{l-17}
Let $\H\subseteq\powcnt\theta$. If\/ $\H$ is $\lambda$-generated 
then $(\H^\perp,\subseteqfnt)$ is $\sigma$-directed.
\end{lem}

The following should be compared to the principle $\princA$.
\begin{principle}{$\aprinc$}
Every directed subfamily  $\H$ of $(\powcnt\theta,\subseteqfnt)$ for which
$\downcl\H$ is generated by a subset $\B\subseteq\H$ of cardinality
$|\B|<\b$, and is locally generated by fewer than $\p$ elements, has either
\begin{enumerate}[leftmargin=*, label=(\arabic*), ref=\arabic*,
  widest=2, labelsep=\label@sep]
\item\label{item:42} a countable decomposition of $\theta$ into
  singletons and sets locally in $\downcl\H$,
\item\label{item:43} an uncountable subset of $\theta$ orthogonal to $\H$.
\end{enumerate}
\end{principle}

\begin{prop}[$\p>\aleph_1$]
\label{p-29}
$\aprinc$ implies $\princA$.
\end{prop}
\begin{proof}
Note that $\p\le\b$.
\end{proof}

\begin{thm}
\label{u-8}
$\pstar$ implies $\aprinc$.
\end{thm}
\begin{proof}
Let $\H$ be as in the hypothesis of $\aprinc$. 
By lemma~\ref{l-18}, equation~\eqref{eq:29} is satisfied. 
And by lemma~\ref{l-17},
$\H^\perp$ is $\sigma$-$\subseteqfnt$-directed. Therefore, $\pstar$ can be applied
to $\H^\perp$. 

First suppose that the second alternative~\eqref{item:51} of $\pstar$ holds. 
Then there is a
countable decomposition of $\theta$ into pieces orthogonal to $\H^\perp$, and
therefore the first alternative of $\aprinc$ holds by
lemmas~\ref{l-22} and~\ref{l-18} since
$\downcl(\H,\subseteqfnt)=\H^{\perp\perp}$. 

Otherwise, $\pstar$ implies that the first alternative~\eqref{item:50} of $\pstar$
holds, and thus there is an uncountable $X\subseteq\theta$ locally in
$\downcl{(\H^\perp)}=\H^\perp$. This means that $X$ is orthogonal to $\H$.
\end{proof}

\begin{corthm}
\label{u-7}
$\p>\aleph_1$ and\/ $\pstar$ together imply\/ $\princA$. 
\end{corthm}

\begin{corthm}[$\p>\aleph_1$]
\label{o-10}
$\pcombmax\to\pcomb\to\pcombboc\to\pstar\to\princA$.
\end{corthm}

In the article~\cite{Hir} we will establish that one of our combinatorial principles,
in conjunction with $\p>\aleph_1$, implies $\princAstar$. Note that it is already
known (cf.~\cite{MR1441232}) that $\ma_{\aleph_1}$ and $\princA$ imply $\princAstar$.

Of course, corollary~\ref{u-7} does not replace the principle $\princA$. On the
other hand, if it is the case that $\princA$ implies $\p>\aleph_1$, then this would
render $\princA$ obsolete. Hence we are interested in the following question (recall
from section~\ref{sec:pstar-does-not} that $\princA\to2^{\aleph_0}>\aleph_1$).

\begin{question}
\label{q-6}
Does $\princA$ imply $\p>\aleph_1$\tu? How about $\b>\aleph_1$\tu?
\end{question}
\subsection{Nearly special Aronszajn trees}
\label{sec:nearly-spec-aronsz}

In~\cite{NS} notions of `almost specialness' for Aronszajn trees were considered. 
The main goal was to show that all Aronszajn trees being `almost special' 
does not necessarily imply that all Aronszajn trees are special in the usual sense. 
The following notion was proposed and appears 
to be optimal with respect to being nearly special.

\begin{defn}
\label{d-20}
We say that an $\oone$-tree $T$ is \emph{$(\NS^+,\NS^*)$-special} if there exists a
maximal antichain $\A$ of $(\NS^+,\subseteq)$ such that every $S\in\A$ has a
relatively closed $C\subseteq S$ on which $T\restriction C$ is special 
(cf.~\Section\ref{sec:terminology}). 
\end{defn}

The following is the main result of~\cite{NS}.

\begin{thmo}[Abraham--Hirschorn]
\label{thmo:AH}
The existence of a nonspecial Aronszajn tree is simultaneously consistent 
with $\ch$ and every Aronszajn tree being $(\NS^+,\NS^*)$-special. 
\end{thmo}

We now show that the results of this paper do indeed generalize this theorem.
Given a tree $T$, we assume without loss of generality that it
resides on some ordinal $\theta$. The associated family, defined in~\cite{MR1441232},
is
\begin{equation}
  \label{eq:30}
  \ideal_T=\bigl\{x\in\powcnt\theta:x\perp\{\pred_T(\xi):\xi\in\theta\}\bigr\},
\end{equation}
in other words, $x\in\powcnt\theta$ is in $\ideal_T$ iff every node of the tree has
only finitely many predecessors in $x$. Clearly $\ideal_T$ is
$\subseteqfnt$-directed and is moreover a lower subset of~$\powcnt\theta$; indeed it is an ideal. 
As is well known (e.g.~\cite{MR1441232}), in the case $\theta=\oone$, if $T$ is an
$\oone$-tree then $\ideal_T$ is $\sigma$-$\subseteqfnt$-directed and thus a
$P$-ideal. More generally, the following lemma~\ref{l-20} follows. 
This embedability property was shown to be equivalent to what 
we called \emph{local countability}
in~\cite{H2}.

\begin{lem}
\label{l-20}
If\/ $T$ is a tree that embeds into an\/ $\oone$-tree then\/ $\ideal_T$ is\/ $\sigma$-directed.
\end{lem}

The following is observed in~\cite{MR1441232}.
Lemma~\ref{l-19} implies that when $X$ is locally in $\ideal_T$ it is a countable
union of antichains (i.e.~special in our terminology). 

\begin{lem}
\label{l-19}
If $X\subseteq\theta$ is locally in $\ideal_T$ then as a subtree of $T$, 
$X$ is of height at most $\omega$. 
\end{lem}

It is noted in~\cite{MR1809418} that if $X$ is orthogonal to
$\ideal_T$, then $X$ is a finite union of chains. We provide a proof
here, using the following basic tree lemma, which to the best of our
knowledge is from the folklore.

\begin{lem}
\label{l-24}
Let\/ $T$ be a tree. If\/ $T$ has no infinite antichains then it is a
finite union of chains.
\end{lem}
\begin{proof}
Let $\C$ be the collection of all maximal chains in $T$. We may assume
that all levels of $T$ are finite. 
Supposing that $\C$ is infinite, $T$ must have infinitely many
splitting nodes.
We choose $s_n,t_n\in T$ by recursion on $n<\omega$ so that each
$s_n$ is a splitting node and $s_n\lln T t_{n+1},s_{n+1}$ with
$s_{n+1}\ne t_{n+1}$ on the same level.
This is possible, because $s_n$ is finitely splitting and thus we can 
recursively guarantee that $\{C\in\C:s_n\in C\}$ is infinite, 
which implies that there is another
splitting node strictly above $s_n$. Then $\{t_1,t_2,\dots\}$ forms an
infinite antichain of~$T$. 
\end{proof}

\begin{cor}{0}
\label{l-21}
If $X\subseteq\theta$ is orthogonal to $\ideal_T$ then it is a finite union of chains.
\end{cor}
\begin{proof}
Note that $\ideal_T$ includes the countable antichains of $T$. 
Thus $X$ contains no infinite antichains, and lemma~\ref{l-24} is
applied to the subtree $(X,\len T)$.
\end{proof}

\begin{thm}
\label{u-10}
$\pstarsplus$ implies that every Aronszajn tree is $(\NS^+,\NS^*)$-special.
\end{thm}
\begin{proof}
Let $\A$ be the maximal antichain of $\NS^+$ given by $\pstarsplus$. 
Let $T$ be an Aronszajn tree, and enumerate each level $T_\alpha$ ($\alpha<\oone$) as
$\xi_{\alpha,0},\xi_{\alpha,1},\dots$. For each~$n$, define a tree
$(U,\len U)$ on $\oone$ by
\begin{equation}
  \label{eq:43}
  \alpha\len{U_n}\beta\If \xi_{\alpha,n}\len T\xi_{\beta,n}.
\end{equation}
Thus each $U_n$ embeds into $T$ via $\alpha\mapsto\xi_{\alpha,n}$. 
And thus each $\ideal_{U_n}$ is $\sigma$-directed by lemma~\ref{l-20}.
By corollary~\ref{l-21}, since $T$ is Aronszajn, 
no $\ideal_{U_n}$ can have an uncountable set orthogonal to it (let
alone stationary). Therefore the second alternative~\eqref{item:53} of
$\pstarsplus$ gives for each $S\in\A$, and
each $n$, a relatively closed $C_{S,n}\subseteq S$ locally in $\ideal_{U_n}$. Then
each $C_{S,n}$ is special as a subtree of $U_n$ by lemma~\ref{l-19},
which means that $\{\xi_{\alpha,n}:\alpha\in C_{S,n}\}$ is special as
a subtree of $T$.
For each $S\in\A$, put $C_S=\bigcap_{n=0}^\infty C_{S,n}$. 
It is now clear that $T\restriction C_S$ is special, completing the proof.
\end{proof}

\vspace*{-50pt}
\bibliographystyle{alpha}
\bibliography{database}

\addtocontents{toc}{\protect\setcounter{unbalance}{0}\protect
\end{multicols}}
\end{document}

%% file: macros.tex
\newcommand{\ch}{\mathrm{CH}}
\newcommand{\gch}{\mathrm{GCH}}
\newcommand{\ma}{\mathrm{MA}}
\newcommand{\zfc}{\mathrm{ZFC}}

\newcommand{\pfa}{\mathrm{PFA}}
\newcommand{\mm}{\mathrm{MM}}

\newcommand{\sh}{\textrm{\rm SH}}

\newcommand{\pstar}{(*)}
\newcommand{\pstarc}{(\gdot_c)}

\newcommand{\length}{\operatorname{len}}

\newcommand{\game}{\Game}

\newcommand{\downcl}[1]{{\downarrow\!#1}}

\newcommand{\gdot}{\star}
\newcommand{\power}{\P}

\newcommand{\inv}{^{-1}}

\newcommand{\otp}{\operatorname{otp}}
\newcommand{\pred}{\operatorname{pred}}

\newcommand\bigext{^\frown}

\newcommand{\lands}{\mathrel{\land}}

\renewcommand{\div}{\mathbin{/}}

\newcommand{\con}{\operatorname{Con}}

\newcommand{\Gen}{\operatorname{Gen}}
\newcommand{\Genc}{\operatorname{Gen}^+}

\newcommand{\spr}{\bigvee}

\newcommand{\@@eq}[2]{=_{#1}^{#2}}
\newcommand{\@@neq}[2]{\ne_{#1}^{#2}}
\newcommand{\@@l}[2]{<_{#1}^{#2}}
\newcommand{\@@nl}[2]{\nless_{#1}^{#2}}
\newcommand{\@@le}[2]{\le_{#1}^{#2}}
\newcommand{\@@nle}[2]{\nleq_{#1}^{#2}}
\newcommand{\@@equiv}[2]{\equiv_{#1}^{#2}}

\newcommand{\@in}[1]{\in_#1}
\newcommand{\@eq}[1]{=_#1}
\newcommand{\@le}[1]{\le_{#1}}
\newcommand{\@l}[1]{<_#1}
\newcommand{\@nle}[1]{\nleq_#1}
\newcommand{\@ge}[1]{\ge_#1}

\newcommand{\u@le}[1]{\le^#1}
\newcommand{\u@l}[1]{<^#1}
\newcommand{\u@eq}[1]{=^#1}

\newcommand{\@lefnt}[1]{\@@le#1*}
\newcommand{\@nlefnt}[1]{\@@nle#1*}
\newcommand{\@eqfnt}[1]{\@@eq#1*}

\newcommand{\subseteqfnt}{\subseteq^*}

\newcommand{\supseteqfnt}{\mathrel{\supseteq^*}}

\newcommand{\proves}{\vdash}
\newcommand{\len}[1]{\@le{#1}}
\newcommand{\lln}[1]{\@l#1}
\newcommand{\gen}[1]{\@ge{{#1}}}
\newcommand{\altle}{\lesssim}

\newcommand{\str}{\mathrm{str}}

\newcommand{\asym}{\mathrm{asym}}
\newcommand{\trn}{\mathrm{tran}}
\newcommand{\negstr}{-\str}

\newcommand{\ggen}[1]{\gg_{#1}}

\newcommand{\strlen}[1]{\@@le{#1}\str}
\newcommand{\nstrlen}[1]{\@@nle{#1}\str}
\newcommand{\strln}[1]{\@@l{#1}\str}
\newcommand{\nstrln}[1]{\@@nl{#1}\str}
\newcommand{\strequiv}[1]{\@@equiv{#1}\str}

\newcommand{\lenn}[2]{\@@le{#1}{#2}}
\newcommand{\nlenn}[2]{\@@nle{#1}{#2}}
\newcommand{\lnn}[2]{\@@l{#1}{#2}}
\newcommand{\nlnn}[2]{\@@nl{#1}{#2}}

\newcommand{\eqnn}[2]{\@@eq{#1}{#2}}
\newcommand{\neqnn}[2]{\@@neq{#1}{#2}}
\newcommand{\asymle}{\u@le\asym}
\newcommand{\asyml}{\u@l\asym}

\newcommand{\strtrnlen}[1]{\@@le{#1}{\str,\trn}}
\newcommand{\strtrnln}[1]{\@@l{#1}{\str,\trn}}
\newcommand{\trnle}{\u@le{\trn}}
\newcommand{\trnl}{\u@l{\trn}}
\newcommand{\trneq}{\u@eq{\trn}}
\newcommand{\negstrlen}[1]{\@@le{#1}\negstr}
\newcommand{\nnegstrlen}[1]{\@@nle{#1}\negstr}
\newcommand{\nnegstrln}[1]{\@@nl{#1}\negstr}
\newcommand{\negstreq}[1]{\@@eq{#1}\negstr}

\renewcommand{\aa}{\mathrm{aa}}
\newcommand{\eqaa}{\@eq\aa}
\newcommand{\inaa}{\@in\aa}

\newcommand{\leaa}{\@le\aa}
\newcommand{\nleaa}{\@nle\aa}
\newcommand{\geaa}{\@ge\aa}
\newcommand{\lefntaa}{\@lefnt\aa}
\newcommand{\nlefntaa}{\@nlefnt\aa}
\newcommand{\eqfntaa}{\@eqfnt\aa}

\newcommand{\naa}{\mathrm{naa}}
\newcommand{\eqnaa}{\@eq\naa}
\newcommand{\eqfntnaa}{\@eqfnt\naa}
\newcommand{\lenaa}{\@le\naa}
\newcommand{\lefntnaa}{\@lefnt\naa}

\newcommand{\eqae}{\@eq\ae}
\newcommand{\leae}{\@le\ae}

\newcommand{\forces}{\mskip 5mu plus 5mu\|\hspace{-2.5pt}{\textstyle \frac{\hspace{4.5pt}}{\hspace{4.5pt}}}\mskip 5mu plus 5mu}

\newcommand{\Iff}{\espc\mathrm{iff}\espc}
\newcommand{\If}{\espc\textrm{if}\espc}
\newcommand{\impls}{\espc\mathrm{implies}\espc}

\renewcommand{\and}{\hespc\mathrm{and}\hespc}

\renewcommand{\b}{\mathfrak b}

\newcommand{\p}{\mathfrak p}

\newcommand{\cof}{\operatorname{cof}}

\newcommand{\Fin}{\mathrm{Fin}}
\newcommand{\pN}{\power(\N)}
\newcommand{\reals}{\mathbb R}

\newcommand{\N}{\mathbb N}

\newcommand{\fin}{\mathrm{fin}}

\newcommand{\NS}{\mathcal{NS}}

\newcommand{\coll}{\mathrm{Col}}

\newcommand{\pmax}{\mathbb P_{\mathrm{max}}}

\DeclareFontFamily{U}{cmsy}{}
\DeclareFontShape{U}{cmsy}{m}{n}{<12> sfixed * [10] cmsy10 
<10> <9> <8> <7> <6> <5> sfixed * [10] cmsy10}{}
\DeclareSymbolFont{customtwo}{U}{cmsy}{m}{n} 
\DeclareMathSymbol{\sctn}{\mathord}{customtwo}{"78}

\DeclareFontFamily{U}{cmmi}{}
\DeclareFontShape{U}{cmmi}{m}{n}{<20> sfixed * [11] cmmib10 <12> sfixed * [10]
cmmi10 <10> <9> <8> sfixed * [6] cmmi6 <5> <6> <7> sfixed * [5] cmmi5}{}
\DeclareSymbolFont{custom}{U}{cmmi}{m}{n}
\DeclareMathSymbol{\rharpoon}{\mathord}{custom}{"2A}
\newlength{\widt}
\newlength{\widttwo}
\newlength{\hgt}

\newcommand{\oone}{{\omega_1}}

\newcommand{\<}{\langle}
\renewcommand{\>}{\rangle}

\newcommand{\espc}{\quad}
\newlength{\@@hespc}
\newcommand{\hespc}{\hspace{\@@hespc}}
\newcommand{\overbar}[1]{\,\overline{\!{#1}}}

\newcommand{\overbarg}[1]{\hsp\overline{\hspb{#1}}}
\newcommand{\Th}{{^{\mathrm{th}}}}

\renewcommand{\ae}{\mathrm{ae}}

\newcommand{\ulc}{\ulcorner}
\newcommand{\urc}{\urcorner}

\newcommand{\tu}{\textup}
\newcommand{\ih}[1]{(\mathrm{IH}_{#1})}

\newcommand{\A}{\mathcal A}
\newcommand{\B}{\mathcal B}
\newcommand{\C}{\mathcal C}
\newcommand{\D}{\mathcal D}
\newcommand{\E}{\mathcal E}
\newcommand{\F}{\mathcal F}

\newcommand{\ideal}{\mathcal I}
\newcommand{\J}{\mathcal J}
\newcommand{\K}{\mathcal K}
\newcommand{\G}{\mathcal G}
\renewcommand{\H}{\mathcal H}
\newcommand{\M}{\mathcal M}
\renewcommand{\O}{\mathcal O}
\renewcommand{\P}{\Pcal}
\newcommand{\Pcal}{\mathcal P}

\newcommand{\R}{\mathcal R}
\newcommand{\Section}{{\char'237}}
\renewcommand{\S}{\mathcal S}
\newcommand{\T}{\mathcal T}

\newcommand{\X}{\mathcal X}

\newcommand{\hsp}{\mspace{1.5mu}}
\newcommand{\hspb}{\mspace{-1.5mu}}
\newcommand{\spc}{\,\,\,}

\newcounter{saveenumi}
\newcommand{\save}{\setcounter{saveenumi}{\value{enumi}}}
\newcommand{\restore}{\setcounter{enumi}{\value{saveenumi}}}

